\documentclass[aps,preprint,groupedaddress,showkeys]{revtex4-1}
\usepackage{graphicx}
\usepackage{dcolumn}
\usepackage{bm}
\usepackage{amssymb}
\usepackage{amsthm}
\usepackage{amsmath}
\usepackage{graphicx}
\usepackage{epstopdf}

\newtheorem{theorem}{\hskip\parindent\bf Theorem}
\newtheorem{remark}{\hskip\parindent\bf Remark}
\begin{document}

\preprint{}

\title[Double Hopf bifurcation  in delayed reaction-diffusion systems]{Double Hopf bifurcation in delayed reaction-diffusion systems}

\author{Yanfei Du}

\affiliation{Shaanxi University of Science and Technology, Xi'an 710021,  China.}

\author{Ben Niu*, Yuxiao Guo, Junjie Wei}

\affiliation{Department of Mathematics, Harbin Institute of Technology, Weihai 264209,  China.\\*Corresponding author, niu@hit.edu.cn}

\date{\today}

\begin{abstract}
Double Hopf bifurcation analysis can be used to reveal some complicated dynamical behavior in a dynamical system, such as the existence or coexistence of  periodic orbits, quasi-periodic orbits, or even chaos. In this paper, an algorithm for deriving the normal form  near a codimension-two double Hopf bifurcation of a reaction-diffusion system with time delay and Neumann boundary condition is rigorously established, by employing the center manifold reduction technique and the normal form method. We find that the dynamical behavior near   bifurcation points are proved to be governed by  twelve distinct unfolding systems.   Two examples are performed to illustrate our results: for a stage-structured epidemic model, we find that double Hopf bifurcation appears when varying the diffusion rate and time delay, and two stable spatially inhomogeneous periodic oscillations are proved to coexist near the bifurcation point; in a diffusive  predator-prey system, we theoretically proved that quasi-periodic orbits exist on two- or three-torus near a double Hopf bifurcation point, which will break down after slight perturbation, leaving the system a strange attractor.

\end{abstract}

                           \maketitle

\section{Introduction}
\label{section1}
A Hopf bifurcation refers to the phenomenon that steady states lose their stability and give rise to periodic solutions, as a parameter crosses  critical values \cite{Andronov,Hale,Hopf,Poincare,Wiggins}. Since the condition of existence of Hopf bifurcation can be easily verified, and both the direction of Hopf bifurcation and the stability of the bifurcating periodic solution can be determined by formulas derived from center manifold reduction technique, Hopf bifurcation analysis has been an effective method of the study on the existence of periodic solutions of differential equations. In recent years, Hopf bifurcation has been intensely studied to investigate the dynamics of differential equations in various fields \cite{J Belair,HW Hethcote,DVR Reddy,Y SongHopf}.  Normal forms theory provides powerful tools to bifurcation analysis, whose  basic idea is to employ  near-identity nonlinear transformations that lead  the original system to a qualitatively equivalent differential equation  with the simplest form.

In a reaction-diffusion system, especially reaction-diffusion system with delays, there may exist several Hopf bifurcation critical points with   parameters' varying. Thus, like the case in ordinary differential equations \cite{Kuznetsov}, double Hopf bifurcations will be characterized by  intersections of two Hopf bifurcation curves in a two-parameter plane.
In fact, a double Hopf bifurcation occurs from a critical point at which the Jacobian evaluated there involves two conjugate pairs of pure imaginary eigenvalues.  Guckenheimer and Holmes \cite{Guckenheimer} gave the dynamics in the neighborhood of a  codimension-two double Hopf point in an ordinary differential equation. When bifurcation parameters are closed to the double Hopf bifurcation point, the system may exhibit rich dynamics, such as periodic oscillations, quasi-periodic oscillations,  coexisting of several oscillations, two- or three-dimensional invariant torus, and even chaos \cite{Guckenheimer,Kuznetsov,Yu}.
For delayed systems governed by delay differential equations, due to the fact that delay differential equations  are of infinite dimensional, the center manifold method was adopted to reduce delay differential equations into finite-dimensional systems, and then normal forms can be calculated on the center manifold \cite{P.L. Buono,S.A. Campell,C. Elphick,FariaJDE,T. FariaBT,J. Halefde}.
Correspondingly, the multiple timescale method has also been successfully applied to analyze
codimension-two  non-resonant double Hopf bifurcations \cite{Luongo,PEI YU}.
 Research on the double Hopf bifurcation and complex dynamics of delay differential equations is also carried out by these methods we mentioned above \cite{Ping BiTumor,S.A. Campellfeedback,Suqi Ma,P. Yu,Y.Y. Zhang}.

In recent years, the analysis of Hopf bifurcation in a reaction-diffusion systems has attracted much attention. Based on the method provided in \cite{B. Hassard,JWu}, Yi et al. \cite{F. Yipredator-prey} carried out Hopf bifurcation analysis in a reaction-diffusion system  without delays. By decomposing the equation with the orthonormal Fourier basis corresponding to the eigenvalue problem, an explicit expression is given for several key parameters in determining the properties of bifurcation according to the method given in \cite{B. Hassard}.
Faria \cite{T. Fariapredator-preydiffusion}
 investigated the effect of both time delay and diffusion on Hopf bifurcation, and improved the method given in \cite{Faria W. Huang,FariaJDE} to study normal forms with perturbation parameters of Hopf bifurcation  in delayed reaction-diffusion equations with Neumann boundary condition.
 The results in \cite{T. Fariapredator-preydiffusion} have been further applied to equations with various practical backgrounds, such as the predator-prey model \cite{S. Chen} and so on.
Recently, Hopf bifurcation analysis of delayed reaction-diffusion equations with  Dirichlet boundary condition has also made considerable progress, such as in  \cite{S. GuoDirichlet,Y. SuDirichlet,X.-P. YanDirichlet}, the existence of Hopf bifurcation and the stability and direction of bifurcating periodic solutions for some population models with Dirichlet boundary condition were considered; Furthermore, Guo \cite{S. Guononlocal},  Chen and Yu \cite{S. Chennonlocal},  combining the non-local phenomenon, considered  Hopf bifurcation of the delayed reaction-diffusion equations with Dirichlet boundary condition. Comparing the research on Dirichlet boundary conditions with  the research on Neumann boundary conditions, the main difficulty when investigating such models is the existence of positive steady-state solutions and the complexity of the corresponding eigenvalue problem, so practical analysis is usually carried out by combining  the method of phase space analysis, topological degree theory or Liapunov-Schmidt reduction theory \cite{H. Kielhfer}. In this paper, only homogeneous Neumann boundary condition is considered, as we can   assume that zero is an equilibrium of the system.

For analysis on high codimensional bifurcations of reaction-diffusion equations, so far as we know, there are very few published results. Only a few results have been conducted mainly on the research of Turing-Hopf bifurcation. Research on interaction between the Turing bifurcation and Hopf bifurcation in  reaction-diffusion equations can be traced back to the work of De Wit et al. \cite{A. De Wit}.
Later, Meixner et al. \cite{M. Meixner} conducted a preliminary analysis of the problem through codimension-two bifurcation analysis.  Parallel research  was also seen in literature \cite{M. Baurmann} on the predator-prey system. Song et al. \cite {YongliSongTuringHopf} investigated the Turing-Hopf bifurcation from the point of view of  analysis of bifurcation and normal forms, in which they extended the method given in \cite{T. Fariapredator-preydiffusion,F. Yipredator-prey}, and deduced the normal form with parameters of  Turing-Hopf bifurcation.
 Xu and Wei \cite{X. Xu} applied this method to study the Turing-Hopf bifurcation of a predator-prey model.
 Turing-Hopf bifurcation can be seen as  a special case of  Hopf-zero bifurcation  of dynamic system from the point of view of classification of bifurcation. Results about Hopf-zero bifurcations and normal form analysis  in delayed reaction-diffusion equations can be found in An and Jiang's preprint \cite{Q. An}. For double Hopf bifurcation, Lewis and Nagata \cite{GREGORY} studied the transitions from axisymmetric steady solutions to nonaxisymmetric waves in a Navier-Stokes model of the differentially
 heated rotating annulus experiment, and an analytical-numerical center manifold reduction
 is used to analyze the double Hopf bifurcation points that occur at this transition.
Besides \cite{GREGORY}, to our best knowledge, there have not been any results  on  double Hopf bifurcation analysis and the corresponding normal form calculation   in the  reaction-diffusion equations with or without delays.

Based on the general method of normal form simplification in ordinary differential equation or delay differential equation, the ratio of two pairs of imaginary roots appearing at the double Hopf bifurcation point is important to determine the final simplest normal form. When  the ratio is not a rational number, the normal form has a universal form as those provided in
 \cite{Ping BiTumor,P.L. Buono,S.A. Campellfeedback,Luongo,Suqi Ma,PEI YU,P. Yu,Y.Y. Zhang}. However, if the ratio is rational, there must exist some additional terms which couldn't be eliminated because of the resonance. Precisely, if we are about to obtain the topological classification of the dynamical behavior near the double Hopf points by analyzing the third order normal form, this only  require  that  the ratio is not $m:n$ for $m,n=1,2,3,4$, namely, weak resonance case. For all weakly resonant double Hopf bifurcation, we can treat them in the same method as that used for non-resonant cases. Otherwise, if the ratio is
 $m:n$ for $m,n=1,2,3,4$, we say that this is a strongly resonance case, which has also been investigated in ordinary differential equations and delay differential equations \cite{AK Bajaj,S. A. CampbellResonant,S. VAN GILS,W. GOVAERTS,Ji,Revel,Revel1:2,PH Steen}.

In this paper, we aim at the development of a method of computing normal forms on center manifold near an equilibrium of a system of a general delayed reaction-diffusion equations with non-resonant or weakly resonant double Hopf singularity.
Consider a general reaction-diffusion model
 \begin{equation}\label{general rd}
 \left\lbrace \begin{array}{l}
 \frac{\partial u_1(x,t)}{\partial t}=d_1(\mu)\Delta u_1(x,t)+f_1(\mu,u_1^t,u_2^t,\cdots,u_n^t),\\
 \frac{\partial u_2(x,t)}{\partial t}=d_2(\mu)\Delta u_2(x,t)+f_2(\mu,u_1^t,u_2^t,\cdots,u_n^t),\\
~~~~ \vdots\\
 \frac{\partial u_n(x,t)}{\partial t}=d_n(\mu)\Delta u_n(x,t)+f_n(\mu,u_1^t,u_2^t,\cdots,u_n^t),\\
 \end{array}\right.~x\in(0,l\pi),
 \end{equation}
 with homogeneous Neumann boundary conditions and $u_i^t(\theta)(x)=u_i(t+\theta,x)$. Define $f_i(\mu,0,\cdots,0)=0$, $i=1,2,\cdots,n$, such that
 (\ref{general rd}) has zero equilibrium, and assume that every $f_i$  is at least $C^3$.

We will extend the normal form method given by Faria and Magalh\~aes \cite{Faria,FariaJDE}   to investigate non-resonant or weakly resonant double Hopf bifurcations in reaction-diffusion systems with delays.  The first step of this method is to decompose the phase space, and then to decompose the equations we investigate. Then we can calculate the second and third order normal form on the center manifold,  and get the normal form for double Hopf bifurcation, which is a four-dimensional system up to the third order with unfolding parameters  restricted to the center manifold. By polar coordinates transformation, the four-dimensional system can be reduced to a two-dimensional amplitude system, in which the unfolding parameters can be expressed by the parameters in the original system. Due to the unfolding analysis given in \cite{Guckenheimer}, we can determine the type of twelve distinct kinds of unfoldings by the coefficients we calculated in the normal form, and get  corresponding  bifurcation sets. By detailed analysis of bifurcation sets, all the dynamical behaviors near the double Hopf bifurcation point can be figured out. In this paper, a  series of  explicit calculation formulas  of the normal form for non-resonant or weakly resonant double Hopf bifurcation is given, with three different cases of wave number: $n_1=0,n_2=0$, $n_1=0,n_2\neq 0$, and  $n_1\neq 0,n_2\neq 0$. The wave number is related to the spatial profile of bifurcating periodic solutions.

 To illustrate the calculation process of the normal forms, we perform two examples. The first one is a diffusive, stage-structured  epidemic model with two delays.  Xiao and Chen \cite{YXiao} proposed an SIS epidemic model with stage structure and time delays in the following form
 \begin{equation}
 \label{xiaoyanni}
 \left\{
 \begin{array}{lll}
 \dfrac{dS}{dt} &=& \alpha y(t)-dS(t)-\alpha e^{-d\tau}y(t-\tau)-\mu S(t-\omega)I(t)+\gamma I(t), \\
 \dfrac{dI}{dt}&=&\mu S(t-\omega)I(t)-dI(t)-\gamma I(t),\\
 \dfrac{dy}{dt}&=&\alpha e^{-d\tau}y(t-\tau)-\beta y^2(t),\\
 \end{array}
 \right.
 \end{equation}
 where $S(t)$ and $I(t)$ are the population of susceptible and infected immature individuals,  $y(t)$ represents the population of mature individuals,  $\tau$ is the maturity delay, and $\omega$ is the freely moving delay. Hopf bifurcation results for such a model are given in Du et al. \cite{duguo}.  Adding random  diffusion of individuals into (\ref{xiaoyanni}),  we have a diffusive epidemic model, which is proposed in section \ref{section epidemic}.
   We take the time delay $\omega$ and a diffusion coefficient as bifurcation parameters. We show that double Hopf bifurcation can appear with two different cases of wave numbers: $n_1=0,n_2\neq 0$, or  $n_1\neq 0,n_2\neq 0$.  Thus, the spatio-temporal dynamics turn out to be very complicated near the double Hopf bifurcation point. In some regions, there are two stable and spatially nonhomogeneous periodic solutions or a homogeneous and a nonhomogeneous periodic solution coexisting.

 The second example we consider is a predator-prey system with time delay.    Predator-prey systems have been widely investigated \cite{S. B. HSUpredator,S. RUANpredator,D. XIAOpredator}.  Here, we choose a simple system which was studied in \cite{Y SongHopf} to illustrate the procedure of normal form calculation
  \begin{equation}\label{predatorintroduction}
  \begin{aligned}
  &\frac{{\rm d} u(t)}{{\rm d} t}=u(t)[r_1-a_{11}u(t-\tau)-a_{12}v(t)],\\
  &\frac{{\rm d} v(t)}{{\rm d} t}=v(t)[-r_2+a_{21}u(t)-a_{22}v(t)].
  \end{aligned}
  \end{equation}
Incorporating diffusion terms into (\ref{predatorintroduction}), the system becomes a reaction-diffusion system, which is proposed in section \ref{section predatorprey}. Taking $r_1$ and $\tau$ as bifurcation parameters, we will show that such a simple system will exhibit complicated dynamical behavior near the double Hopf bifurcation point such as the existence of quasi-periodic solution on a 2-torus and quasi-periodic solution on a 3-torus. Generally, a vanishing 3-torus might accompany strange attractors, and lead a system into chaos, which is called the ``Ruelle-Takens-Newhouse" scenario \cite{P. Battelino,J.P. Eckmann,D. Ruelle}. Thus a strange attractor is found near the double Hopf bifurcation point in this predator-prey system.

The paper is organized as follows. In section \ref{PFDE}, we discuss eigenfunctions and decomposition of the phase space. Section \ref{Center manifold reduction and normal form} is devoted to calculating the normal forms of non-resonant or weakly resonant double Hopf bifurcation of (\ref{general rd}), and  explicit formula of normal form truncated to the third order is obtained in three cases of different wave numbers. As applications of
the method, a diffusive epidemic model with two delays and a diffusive predator-prey system with a delay are considered in sections  \ref{section epidemic} and \ref{section predatorprey}, respectively, where the conditions for the existence of double Hopf bifurcation are obtained, the normal forms are calculated using the method and formulas in section \ref{Center manifold reduction and normal form}. Two-parameter unfoldings and bifurcation diagrams near the critical point are analyzed, and  some numerical simulations about the rich dynamics near these bifurcation points are demonstrated. Finally, we close the paper with some conclusions.
\section{Eigenfunctions and Decomposition of the Phase Space}
\label{PFDE}
Normal form analysis usually depends on a decomposition of the corresponding phase space \cite{Kuznetsov,Wiggins}. In case of reaction-diffusion model with time delay, it also requires the decomposition with respect to spatial variables which can be accomplished by finding the eigenfunctions of Laplacian operator. Thus, we first rewrite system (\ref{general rd}) into an abstract ordinary differential equation in an appropriate phase space.
\subsection{Eigenfunctions}
Using the general approach to put a partial differential equation into an abstract form introduced in \cite{JWu}, we
define the real-valued Hilbert space \[X:=\left\lbrace (u_1,u_2,\cdots,u_n)^T\in(H^2(0,l\pi))^n:\frac{\partial u_i}{\partial x}(0,t)=\frac{\partial u_i}{\partial x}(l\pi,t)=0,i=1,2,\cdots,n\right\rbrace.\]
Since we are about to deal with the double Hopf bifurcations with two pairs of purely imaginary eigenvalues, we usually extend the space $X$ into the
 corresponding complexification space of $X$ in a natural way by
\[X_{\mathbb{C}}:=X\oplus iX=\{U_1+iU_2:U_1,U_2\in X\}\]
with  the general complex-value   $L^2$ inner product
\[\langle u,v\rangle=\int_0^{l\pi}(\overline{u}_1v_1+\overline{u}_2v_2+\cdots+\overline{u}_nv_n)dx,\]
for $u=(u_1,u_2,\cdots,u_n)^T$, $v=(v_1,v_2,\cdots,v_n)^T\in X_{\mathbb{C}}$. Let $\mathcal{C}:=C([-r,0],X_{\mathbb{C}})$ denote the phase space with the sup norm. We write $u^t\in\mathcal{C}$ for $u^t(\theta)=u(t+\theta)$, $-r\leq \theta\leq 0$.

Now we can rewrite  system (\ref{general rd}) into an abstract form as
  \begin{equation}
  \label{general}
  \frac{{\rm{d}}U(t)}{{\rm{d}} t}=D(\mu)\Delta U(t)+L(\mu)(U^t)+F(\mu,U^t),~~ t>0,
  \end{equation}
  where $U(t)=\left( \begin{array}{c}
   u_1(t)\\u_2(t)\\\vdots\\u_n(t)
  \end{array} \right)\in X_{\mathbb{C}}$,  $u^t=\left( \begin{array}{c}
     u_1^t\\u_2^t\\\vdots\\u_n^t
    \end{array} \right)\in \mathcal{C}$,
     $D(\mu)=\left(\begin{array}{cccc}
     d_1(\mu)&0&\cdots&0\\0&d_2(\mu)&\cdots&0\\\vdots&\vdots&\ddots&\vdots\\0&0&\cdots&d_n(\mu)
     \end{array} \right)$, $F(\mu,U^t)=\left( \begin{array}{c}
     F^{(1)}(\mu,U^t)\\F^{(2)}(\mu,U^t)\\\vdots\\F^{(n)}(\mu,U^t)
       \end{array}\right)$,
      and $d_i(\mu)>0$ for $i=1,2,\cdots, n$, $\mu\in\mathbb{R}^2$. $L:\mathbb{R}^2\times \mathcal{C}\rightarrow X_{\mathbb{C}}$ is a bounded linear operator.
  $F(\mu,\phi)\in X_{\mathbb{C}}$ for $(\mu,\phi)\in\mathbb{R}^2\times \mathcal{C}$, and $F$ is a $C^k$  $(k\geq 3)$ function such that $F(\mu,0)=0$, $D_{\phi} F(\mu,0)=0$, where $D_{\phi} F(\mu,0)$ stands for the Fr\'{e}chet derivative of $F(\mu,\phi)$ with respect to $\phi$ at $\phi=0$.

 Linearizing system  (\ref{general}) at $(0,0,\cdots,0)$, we have
   \begin{equation}
    \label{general linearized}
    \frac{{\rm{d}}U(t)}{{\rm{d}} t}=D(\mu)\Delta U(t)+L(\mu)(U^t).
    \end{equation}
It is well known that the eigenvalue problem
\begin{equation*}
\Delta \varphi=\lambda\varphi,~~x\in(0,l\pi),~\varphi_x\mid_{x=0,l\pi}=0,
\end{equation*}
has eigenvalues
\begin{equation*}
\lambda_m=-\frac{m^2}{l^2},~m\in \mathbb{N}_0= \mathbb{N}\cup\{0\},
\end{equation*}
with corresponding normalized eigenfunctions
\begin{equation*}
\gamma_m(x)=\dfrac{\cos\frac{m}{l}x}{\parallel\cos\frac{m}{l}x\parallel_{L^2}}=\left\lbrace \begin{array}{ll}
\sqrt{\frac{1}{l\pi}},& m=0,\\
\sqrt{\frac{2}{l\pi}}\cos\frac{m}{l}x, & m\geq 1.
\end{array}
\right.\end{equation*}
Let $\beta_m^{(j)}(x)=\gamma_m(x)e_j$, where $e_j$ is the $j$-th unit coordinate vector of $\mathbb{R}^n$. Then $\{\beta_m^{(j)}\}_{m\geq 0}$ are eigenfunctions of $D(\mu)\Delta$ with corresponding eigenvalues $-d_i(\mu)\frac{m^2}{l^2}$ ($i=1,2,\cdots n$). Applying the general theory about elliptic operators, we know   $\{\beta_m^{(j)}\}_{m\geq 0}$ form an orthnormal basis of $X$.

Define $\mathcal{B}_m$ the subspace of $\mathcal{C}$ by
\[\mathcal{B}_m:={\rm span} \left\lbrace \langle v(\cdot),\beta_m^{(j)}\rangle\beta_m^{(j)}~\arrowvert  ~v\in \mathcal{C},j=1,2,\cdots,n\right\rbrace. \]
For simplification of notations, in the following we write
 \[\left\langle v(\cdot),\beta_m \right\rangle=\left( \begin{array}{c}
\langle v(\cdot),\beta_m^{(1)} \rangle \\
\langle v(\cdot),\beta_m^{(2)} \rangle \\
\vdots\\
\langle v(\cdot),\beta_m^{(n)} \rangle \\
\end{array}\right). \]
The characteristic equation associated with  (\ref{general linearized}) is
\begin{equation}\label{characteristic1}
\lambda y-D(\mu)\Delta y-L(\mu)(e^{\lambda\cdot}y)=0,
\end{equation}
where $e^{\lambda\cdot}(\theta)y=e^{\lambda\theta}y$, for $\theta\in[-r,0]$.
By using the Fourier expansion
\begin{equation*}
y=\sum_{m=0}^\infty\left(\begin{array}{c}
a_{1m}\\a_{2m}\\\vdots\\ a_{nm}
\end{array} \right) \gamma_m,
\end{equation*}
we find that the characteristic equation (\ref{characteristic1}) is equivalent to a sequence of characteristic equations
\begin{equation}
\label{characteristic2}
{\rm det}\left[ \lambda I+\frac{m^2}{l^2}D(\mu)-L(\mu)(e^{\lambda\cdot}I)\right]=0,~m\in\mathbb{N}_0.
\end{equation}
In order to consider  double Hopf bifurcation, we  assume  that  the following conditions hold  for some $\mu_0\in \mathbb{R}^2$ \cite{Kuznetsov}:

$(H_1)$ There exists a neighborhood $\textit{O}$ of $\mu_0$ such that for $\mu\in\textit{O}$, (\ref{general linearized}) has two pairs of complex simple conjugate eigenvalues $\alpha_1(\mu)\pm i\omega_1(\mu)$ and $\alpha_2(\mu)\pm i \omega_2(\mu)$, all continuously differential in $\mu$ with $\alpha_1(\mu_0)=0$ , $\omega_1(\mu_0)=\omega_1>0$, $\alpha_2(\mu_0)=0$ , $\omega_2(\mu_0)=\omega_2>0$,
  and the remaining eigenvalues of $(\ref{general linearized})$ have non-zero real part for $\mu\in \textit{O}$.

$(H_2)$ Assume that $\omega_1<\omega_2$ and  $\omega_1:\omega_2\neq i:j$ for $i,j\in\mathbb{N}$ and $1\leq i\leq j \leq 4$, i.e., we do not consider the strongly resonant cases.

$(H_3)$ The conjugate eigenvalues $\alpha_k(\mu)\pm i\omega_k(\mu)$ are obtained by $(\ref{characteristic2}_{n_k})$ and the corresponding eigenvalues belong to $\mathcal{B}_{n_k}$, $k=1,2$, and $n_1, n_2 \in\mathbb{N}_0 $.  Assume that $n_1\leq n_2$.

To figure out the spatiotemporal dynamical behavior near the critical point $\mu=\mu_0$, let $\mu=\alpha+\mu_0$, where $\alpha\in  \mathbb{R}^2$. Now the system (\ref{general linearized}) is equivalent to
\begin{equation}
\label{generalalpha}
\frac{{\rm {d}} U(t)}{\rm {d}t}=D_0\Delta U(t)+L_0(U^t)+\widetilde{F}(\alpha,U^t),
\end{equation}
where $D_0=D(\mu_0)$, $L_0(\cdot)=L(\mu_0)(\cdot)$ is a linear operator from $\mathcal{C}$ to $X_{\mathbb{C}}$, and $\widetilde{F}(\alpha,\varphi)=[D(\alpha+\mu_0)-D_0]\Delta\varphi(0)+[L(\alpha+\mu_0)-L_0](\varphi)+F(\alpha+\mu_0,\varphi)$.
\subsection{Decomposition of the Phase Space}
In order to adapt the  center manifold reduction technique, we have to operate on an equation about $U^t$, which requires  an enlarged   phase space $\mathcal{BC}$ defined by
\[\mathcal{BC}:=\{\psi:[-r,0]\rightarrow X_{\mathbb{C}}:\psi {\rm ~is~ continuous~ on~} [-r,0), \exists \lim_{\theta\rightarrow 0^-}\psi(\theta)\in X_{\mathbb{C} }\}.\]
Eq. (\ref{generalalpha}) can be rewritten as an abstract ordinary differential equation on $\mathcal{BC}$ \cite{Faria}:
\begin{equation}\label{ode}
\dfrac{{\rm d}U^t}{{\rm d}t}=AU^t+X_0\widetilde{F}(\alpha,U^t),
\end{equation}
where $A$ is the infinitesimal generator of the $C_0$-semigroup of solution maps of  the linear equation (\ref{general linearized}), defined by
\begin{equation}
\label{A}
A:\mathcal{C}_0^1\cap\mathcal{BC}\rightarrow \mathcal{BC}, ~A\varphi=\dot{\varphi}+X_0[D_0\Delta\varphi(0)+L_0(\varphi)-\dot{\varphi}(0)],
\end{equation}
with ${\rm dom}(A)=\{\varphi\in\mathcal{C}:\dot{\varphi}\in\mathcal{C},\varphi(0)\in {\rm dom}(\Delta)\}$ and $X_0$ is given by
\[X_0(\theta)=\left\lbrace \begin{array}{cc}
0,&  -r\leq\theta< 0,\\
I,&\theta=0.
\end{array}\right.\]
Then on $\mathcal{B}_m$, the linear equation $$\frac{{\rm d}U(t)}{{\rm d}t}=D_0\Delta U(t)+L_0(U^t)$$ is equivalent to the retarded functional differential equation on $\mathbb{C}^n$:
\begin{equation}
\label{RFDE}
\dot{y}(t)=-\frac{m^2}{l^2}D_0y(t)+L_0y^t.
\end{equation}
Define functions of bounded variation $\eta_k\in BC0([-r,0],\mathbb{R})$ for  $n_k$ $(k=1,2)$ mentioned in $(H_3)$, such that
\begin{equation*}
-\frac{n_k^2}{l^2}D_0\varphi(0)+L_0(\varphi)=\int_{-r}^0d\eta_k(\theta)\varphi(\theta), \varphi\in \mathcal{C}.
\end{equation*}
Let $A_k$ ($k=1,2$) denote the infinitesimal generator of the semigroup defined by (\ref{RFDE}) with $\mu =\mu_0, m = n_1,n_2$,  and $A_k^*$ denotes the formal adjoint of $A_k$ under the bilinear form
\begin{equation*}
(\psi,\phi)_k=\psi(0)\phi(0)-\int_{-r}^0\int_0^\theta\psi(\xi-\theta)d\eta_k(\theta)\phi(\xi)d\xi.
\end{equation*}
Let
\begin{equation*}
\begin{array}{ll}
\Phi_1(\theta)=(\phi_1(\theta),\phi_2(\theta)),&
\Phi_2(\theta)=(\phi_3(\theta),\phi_4(\theta)),\\
\Psi_1(s)=\left( \begin{array}{c}
\psi_1(s)\\\psi_2(s)
\end{array}\right),&
\Psi_2(s)= \left( \begin{array}{c}
\psi_3(s)\\\psi_4(s)
\end{array}\right),
\end{array}
\end{equation*}
be the basis of the generalized eigenspace of $A_k$, $A_k^*$ corresponding to the eigenvalues  $i\omega_1$ and $i\omega_2$, respectively, and satisfying
\[A_k\Phi_k=\Phi_kB_k,~A_k^*\Psi_k=B_k \Psi_k,~(\Psi_k,\Phi_k)_k=I, k=1,2,\]
with $B_1={ \rm diag} (i\omega_1,-i\omega_1)$, $B_2={ \rm diag} (i\omega_2,-i\omega_2)$.
Denote $\Phi(\theta)=(\Phi_1(\theta),\Phi_2(\theta))$, and $ \Psi(s)={ (\Psi_1(s),\Psi_2(s))^T}.$ Now, we can decompose $\mathcal{BC}$ into a center subspace and its orthocomplement, i.e.,
\begin{equation}\label{ker}
\mathcal{BC}=\mathcal{P}\oplus{\rm Ker}\pi,
\end{equation}
where $\pi:\mathcal{BC}\rightarrow\mathcal{P}$ is the projection defined by
\begin{equation*}
\pi(\varphi)=\sum_{k=1}^2\Phi_k(\Psi_k,\langle\varphi(\cdot),\beta_{n_k}\rangle)_k\cdot\beta_{n_k},
\end{equation*}
with $\beta_{n_k}=\left( \beta_{n_k}^{(1)}, \beta_{n_k}^{(2)},\cdots \beta_{n_k}^{(n)}\right) $,  $\langle \varphi(\cdot),\beta_{n_k}\rangle=\left( \begin{array}{c}
\langle \varphi(\cdot),\beta_{n_k}^{(1)} \rangle \\
\langle \varphi(\cdot),\beta_{n_k}^{(2)} \rangle \\
\vdots\\
\langle \varphi(\cdot),\beta_{n_k}^{(n)} \rangle \\
\end{array}\right)$, $c\cdot \beta_{n_k}=c_1 \beta_{n_k}^{(1)}+c_2 \beta_{n_k}^{(2)}+\cdots +c_n \beta_{n_k}^{(n)}$ for $c=(c_1,c_2,\cdots,c_n)^T\in \mathcal{C}$.

Decompose $U^t\in \mathcal{C}_0^1$ according  to (\ref{ker}) as
\begin{equation}\label{Uttheta}
U^t(\theta)=\sum_{k=1}^2\Phi_k(\theta)\widetilde{z}_k(t)\cdot\beta_{n_k}+y_t(\theta),
\end{equation}
where $\widetilde{z}_k(t)=(\Psi_k,\langle U^t,\beta_{n_k}\rangle)_k$, and $y_t \in \mathcal{C}_0^1\bigcap{\rm Ker}\pi:=\mathcal{Q}^1$ for any $t$.
Then in $\mathcal{BC}$ the system (\ref{ode}) is equivalent to the system
\begin{equation}\label{zdoty}
\begin{array}{l}
\dot{z}=Bz+\Psi(0)\left(\begin{array}{c}
\langle\widetilde{F}(\alpha,\sum\limits_{k=1}^2(\Phi_k\widetilde{z}_k)\cdot\beta_{n_k}+y),\beta_{n_1}\rangle\\
\langle\widetilde{F}(\alpha,\sum\limits_{k=1}^2(\Phi_k\widetilde{z}_k)\cdot\beta_{n_k}+y),\beta_{n_2}\rangle\\
\end{array} \right),\\
\frac{{\rm d}y}{{\rm d}t} =A_1y+(I-\pi)X_0\widetilde{F}(\alpha,\sum\limits_{k=1}^2(\Phi_k\widetilde{z}_k)\cdot\beta_{n_k}+y),
\end{array}
\end{equation}
where $z=(\widetilde{z}_1,\widetilde{z}_2)^T:=(z_1,z_2,z_3,z_4)^T$, $B={\rm diag (i\omega_1,-i\omega_1,i\omega_2,-i\omega_2)}$, and $A_1$ is  the restriction of $A$ on  $\mathcal{Q}^1\subset{\rm Ker \pi\rightarrow{\rm Ker}\pi}$, $A_1\varphi=A\varphi$ for $\varphi\in \mathcal{Q}^1$.

\section{Center Manifold Reduction and Normal Form}
\label{section3}
In the previous section, we have calculated the basis of the center subspace and defined a projection onto it.  Applying the center manifold theory \cite{Lin,JWu}, we know there exists an invariant local center manifold of the equilibrium, which will be calculated approximately and be used to obtain the normal form.
\label{Center manifold reduction and normal form}
\subsection{Center Manifold Reduction}
Consider the formal Taylor expansion
\begin{equation*}
\widetilde{F}(\alpha,\varphi)=\sum_{j\geq 2}\frac{1}{j!}\widetilde{F}_j(\alpha,\varphi),
\end{equation*}
where $\widetilde{F}_j$ is the $j{\rm th}$ Fr\'{e}chet derivation of $\widetilde{F}$. Then (\ref{zdoty}) can be written  as
\begin{equation}\label{zdotTaylor}
\begin{array}{l}
\dot{z}=Bz+\sum\limits_{j\geq 2}\frac{1}{j!}f_j^1(z,y,\alpha),\\
\dfrac{{\rm d}y}{{\rm d}t} =A_1y+\sum\limits_{j\geq 2}\frac{1}{j!}f_j^2(z,y,\alpha),
\end{array}
\end{equation}
where $z=(z_1,z_2,z_3,z_4)^T\in \mathbb{C}^4, y\in \mathcal{Q}^1$, and $f_j=(f_j^1,f_j^2), j\geq 2$, are defined by
\begin{equation}
\label{fj1fj2}
\begin{array}{l}
f_j^1(z,y,\alpha)=\Psi(0)\left(\begin{array}{c}
\langle\widetilde{F}_j(\alpha,(\sum\limits_{k=1}^2\Phi_k\widetilde{z}_k)\cdot\beta_{n_k}+y),\beta_{n_1}\rangle\\
\langle\widetilde{F}_j(\alpha,(\sum\limits_{k=1}^2\Phi_k\widetilde{z}_k)\cdot\beta_{n_k}+y),\beta_{n_2}\rangle\\
\end{array} \right),\\
f_j^2=(I-\pi)X_0\widetilde{F}_j(\alpha,(\sum\limits_{k=1}^2\Phi_k\widetilde{z}_k)\cdot\beta_{n_k}+y).
\end{array}
\end{equation}

Let us introduce the following notations as those in \cite{Faria}: for a normed space $Y$, $V_j^6(Y)$ denotes the space of homogeneous polynomials of degree $j$ in $6$ variables $z=(z_1,z_2,z_3,z_4)^T$, $\alpha=(\alpha_1,\alpha_2)^T$ with coefficients in $Y$,
\[V_j^{6}(Y)=\left\lbrace\sum_{|(q,l)|=j}c_{(q,l)}z^q\alpha^l:(q,l)\in \mathbb{N}_0^{6},c_{(q,l)}\in Y \right\rbrace \]
and the  norm $|\sum\limits_{|(q,l)|=j}c_{(q,l)}z^q\alpha^l|=\sum\limits_{|(q,l)|=j}|c_{(q,l)}|_Y$.
Define the operator $M_j=(M_j^1,M_j^2)$, $j\geq 2$ by
\begin{equation}\label{mjU}
\begin{array}{ll}
M_j^1:V_j^{6}(\mathbb{C}^4)\rightarrow V_j^{6}(\mathbb{C}^4),\\(M_j^1p)(z,\alpha)=D_zp(z,\alpha)Bz-Bp(z,\alpha),\\
M_j^2:V_j^{6}(\mathcal{Q}_1)\subset V_j^{6}({\rm Ker\pi})\rightarrow  V_j^{6}({\rm Ker\pi}),\\(M_j^2h)(z,\alpha)=D_zh(z,\alpha)Bz-A_1h(z,\alpha).\\
\end{array}
\end{equation}

According to \cite{Faria}, the normal forms are obtained by a recursive transformations of variables of the form
\begin{equation}
\label{zyalpha}
(z,y,\alpha)=(\widehat{z},\widehat{y},\alpha)+\frac{1}{j!}(U_j^1(\widehat{z},\alpha),U_j^2(\widehat{z},\alpha),0),
\end{equation}
with $U_j=(U_j^1,U_j^2)\in V_j^{6}(\mathbb{C}^4)\times V_j^{6}(\mathcal{Q}_1)$.
This recursive process transforms (\ref{zdotTaylor}) into the normal form
\begin{equation*}\label{zdotnormalform}
\begin{array}{l}
\dot{z}=Bz+\sum\limits_{j\geq 2}\frac{1}{j!}g_j^1(z,y,\alpha),\\
\frac{dy}{dt} =A_1y+\sum\limits_{j\geq 2}\frac{1}{j!}g_j^2(z,y,\alpha),
\end{array}
\end{equation*}
where $g_j=(g_j^1,g_j^2), j\geq 2$, are the new terms of order $j$, given by
\[g_j(z,y,\alpha)=\overline{f}_j(z,y,\alpha)-M_jU_j(z,\alpha),\]
and $U_j\in V_j^{6}(\mathbb{C}^4)\times V_j^{6}(\mathcal{Q}_1)$
are given by
\begin{equation}\label{Ujzalpha}
U_j(z,\alpha)=(M_j)^{-1}\textbf{P}_{{\rm Im}(M_j^1)\times{\rm Im}(M_j^2)}\circ\overline{f}_j(z,0,\alpha),
\end{equation}
where $\textbf{P}$ is the projection operator, and $\overline{f}_j=(\overline{f}_j^1,\overline{f}_j^2)$ denote the terms of order $j$ in $(z,y)$ obtained after the computation of normal forms up to order $j-1$.

  \subsection{Second Order Normal Form}
For simplification of notations, we write $z=(z_1,z_2,z_3,z_4)^T$, $\Phi(\theta)=(\Phi_1(\theta),\Phi_2(\theta))$, $z_x=(z_1\gamma_{n_1},z_2\gamma_{n_1},z_3\gamma_{n_2},z_4\gamma_{n_2})^T$. (\ref{Uttheta}) can be written  as
\begin{equation}
\label{Ut}\begin{array}{l}
U^t(\theta)\\
=\Phi_1(\theta)\left( \begin{array}{c}
z_1(t)\\z_2(t)
\end{array}\right)\cdot\beta_{n_1} +\Phi_2(\theta)\left( \begin{array}{c}
z_3(t)\\z_4(t)
\end{array}\right)\cdot\beta_{n_2}+y_t(\theta)\\
=\Phi_1(\theta)\left( \begin{array}{c}
z_1(t)\\z_2(t)
\end{array}\right)\gamma_{n_1} +\Phi_2(\theta)\left( \begin{array}{c}
z_3(t)\\z_4(t)
\end{array}\right)\gamma_{n_2}+y_t(\theta)\\
=\Phi(\theta)\left(\begin{array}{c}
z_1(t)\gamma_{n_1}\\z_2(t)\gamma_{n_1}\\z_3(t)\gamma_{n_2}\\z_4(t)\gamma_{n_2}
\end{array} \right) +\left(\begin{array}{c}
 y_{1t}(\theta)\\
 \vdots\\
  y_{nt}(\theta)
\end{array}\right)\\
=\phi_1(\theta)z_1(t)\gamma_{n_1}+\phi_2(\theta)z_2(t)\gamma_{n_1}+\phi_3(\theta)z_3(t)\gamma_{n_2}+\phi_4(\theta)z_4(t)\gamma_{n_2}+y_t(\theta)\\ \stackrel{\bigtriangleup}{=}\Phi(\theta) z_x(t)+y_t(\theta).
\end{array}
\end{equation}

Then, by (\ref{zdoty}), (\ref{general}) is equivalent to the system
\begin{equation*}\label{zdotA01}
\begin{array}{l}
\dot{z}=Bz+\Psi(0)\left(\begin{array}{c}
\langle\widetilde{F}(\alpha,\Phi z_x+y),\beta_{n_1}\rangle\\
\langle\widetilde{F}(\alpha,\Phi z_x+y),\beta_{n_2}\rangle\\
\end{array} \right),\\
\dfrac{{\rm d}y}{{\rm d}t} =A_1y+(I-\pi_0)X_0\widetilde{F}(\alpha,\Phi z_x+y).
\end{array}
\end{equation*}

When  $\omega_1:\omega_2\neq m:n$ for $m,n\in\mathbb{N}$ and $1\leq m,n\leq 3$, by (\ref{mjU}), it is easy to verify that
\begin{equation}
\label{Mj1}
\begin{array}{l}
M_j^1(z^p\alpha^\iota e_{\xi})=D_z(z^p\alpha^\iota e _{\xi})Bz-Bz^p\alpha^\iota e_{\xi}\\
=\left\lbrace \begin{array}{l}
(i\omega_1p_1-i\omega_1p_2+i\omega_2p_3-i\omega_2p_4+(-1)^{\xi} i\omega_1)z^p\alpha^\iota e_{\xi},\xi=1,2,\\
(i\omega_1p_1-i\omega_1p_2+i\omega_2p_3-i\omega_2p_4 +(-1)^{\xi} i\omega_2)z^p\alpha^\iota e_{\xi},\xi=3,4,
\end{array}
\right.
\end{array}
\end{equation}
where $\xi=1,2,3,4$,  $e_1=(1,0,0,0)^T$, $e_2=(0,1,0,0)^T$, $e_3=(0,0,1,0)^T$, $e_4=(0,0,0,$ $1)^T$,  $z^p=z_1^{p_1}z_2^{p_2}z_3^{p_3}z_4^{p_4}$,  $\alpha^\iota=a_1^{\iota_1}a_2^{\iota_2}$, $p_1,p_2,p_3,p_4,\iota_1,\iota_2 \in \mathbb{N}_0$, $p_1+p_2+p_3+p_4+\iota_1+\iota_2=j$.
Therefore,
\begin{equation}
\label{M21}
{\rm Im}(M_2^1)^c={\rm span}\left\lbrace \left(\begin{array}{c}\alpha_iz_1\\0\\0\\0\end{array} \right) ,
\left(\begin{array}{c}0\\\alpha_iz_2\\0\\0\end{array} \right) ,
\left(\begin{array}{c}0\\0\\\alpha_iz_3\\0\end{array} \right) ,
\left(\begin{array}{c}0\\0\\0\\\alpha_iz_4\end{array} \right)\right\rbrace ,
 i=1,2.
\end{equation}

Thus, the second order normal form of (\ref{general}) on the center manifold of the origin near $\mu=\mu_0$ has the following form
\begin{equation*}
\dot{z}=Bz+\frac{1}{2!}g_2^1(z,0,\alpha)+\cdots,
\end{equation*}
with $g_2^1(z,0,\alpha)={\rm Proj}_{{\rm Im}(M_2^1)^c}f_2^1(z,0,\alpha)$.

Write the Taylor expansions of $L(\mu_0+\alpha)$ and $D(\mu_0+\alpha)$ as follows
\begin{equation*}
L(\alpha+\mu_0)=L_0+\alpha_1L_1^{(1,0)}+\alpha_2L_1^{(0,1)}+\frac{1}{2}\left( \alpha_1^2L_2^{(2,0)}+2\alpha_1\alpha_2L_2^{(1,1)}+\alpha_2^2L_2^{(0,2)}\right) +\cdot\cdot\cdot,
\end{equation*}
\begin{equation*}
D(\alpha+\mu_0)=D_0+\alpha_1D_1^{(1,0)}+\alpha_2D_1^{(0,1)}+\frac{1}{2}\left( \alpha_1^2D_2^{(2,0)}+2\alpha_1\alpha_2D_2^{(1,1)}+\alpha_2^2D_2^{(0,2)}\right) +\cdot\cdot\cdot.
\end{equation*}

Thus, the second order term of $\widetilde{F}(\alpha,U^t)$ is
\begin{equation*}\label{F2Ut}
\begin{aligned}
\widetilde{F}_2(\alpha,U^t)=2\left( \alpha_1D_1^{(1,0)}+\alpha_2D_1^{(0,1)}\right)\Delta U^t(0)+ 2\left( \alpha_1L_1^{(1,0)}+\alpha_2L_1^{(0,1)}\right)U^t+F_2(\alpha,U^t).
\end{aligned}
\end{equation*}

By (\ref{Ut}),  we can write
\begin{equation}
\label{F2}
\begin{aligned}
\widetilde{F}_2(z,y,\alpha)
&=\widetilde{F}_2(\alpha,\Phi z_x+y)\\
&=2\left( \alpha_1D_1^{(1,0)}+\alpha_2D_1^{(0,1)}\right)\Delta(\Phi(0)z_x+y(0))\\&+2\left( \alpha_1L_1^{(1,0)}+\alpha_2L_1^{(0,1)}\right)(\Phi z_x+y)+F_2(\alpha,\Phi z_x+y).\\
\end{aligned}
\end{equation}
Recalling assumptions $F(\mu,0)=0$, $D_{\phi}F(\mu,0)=0$, we have $F_2(\alpha,\Phi z_x+y)=F_2(0,\Phi z_x+y)$, which can be expressed by
\begin{equation}
\label{F2Phi}
\begin{aligned}
F_2(0,\Phi z_x+y)&=\sum_{q_1+q_2+q_3+q_4=2}F_{q_1q_2q_3q_4}\gamma_{n_1}^{q_1+q_2}(x)\gamma_{n_2}^{q_3+q_4}(x)z_1^{q_1}z_2^{q_2}z_3^{q_3}z_4^{q_4}\\
& +S_2(\Phi z_x,y)+o(|y|^2).
\end{aligned}
\end{equation}
Here $S_2(\Phi z_x,y)$  represents the linear terms of $y$, which can be calculated by $D_\phi F_2(0,\Phi z_x+y)\mid_{y=0}(y)$.

From (\ref{fj1fj2}), (\ref{Ut}), and (\ref{F2}), we have
\begin{equation*}
\frac{1}{2!}f_2^1(z,0,\alpha)=\frac{1}{2!}\Psi(0)\left( \begin{array}{l}
\langle \widetilde{F}_2(z,0,\alpha),\beta_{n_1}\rangle \\
\langle \widetilde{F}_2(z,0,\alpha),\beta_{n_2}\rangle
\end{array}
\right) .
\end{equation*}
According to (\ref{M21}) and the fact that
\begin{equation*}
\int_0^{l\pi}\gamma_{n_1}^2dx=\int_0^{l\pi}\gamma_{n_2}^2dx=1,
\end{equation*}
we obtain
\begin{equation}
\label{2!g21zalpha}
\frac{1}{2!}g_2^1(z,0,\alpha)=\frac{1}{2!}{\rm Proj}_{{\rm Im}(M_2^1)^c}f_2^1(z,0,\alpha)=\left(\begin{array}{c}
B_{11}\alpha_1z_1+B_{21}\alpha_2z_1\\
\overline{B_{11}}\alpha_1z_2+\overline{B_{21}}\alpha_2z_2\\
B_{13}\alpha_1z_3+B_{23}\alpha_2z_3\\
\overline{B_{13}}\alpha_1z_4+\overline{B_{23}}\alpha_2z_4
\end{array} \right) ,
\end{equation}
with
\begin{equation}\label{B11}
\begin{array}{l}
B_{11}=\psi_1(0)\left(-\frac{n_1^2}{l^2}D_1^{(1,0)} \phi_1(0)+L_1^{(1,0)}\phi_1\right), \\
B_{21}=\psi_1(0)\left(-\frac{n_1^2}{l^2}D_1^{(0,1)} \phi_1(0)+L_1^{(0,1)}\phi_1\right),\\
B_{13}=\psi_3(0)\left(-\frac{n_2^2}{l^2}D_1^{(1,0)} \phi_3(0)+L_1^{(1,0)}\phi_3\right), \\
B_{23}=\psi_3(0)\left(-\frac{n_2^2}{l^2}D_1^{(0,1)} \phi_3(0)+L_1^{(0,1)}\phi_3\right).\\
\end{array}
\end{equation}

\begin{remark}\label{2orderdiscrete}
For many systems with discrete time delays, if there are $m$ discrete  delays in system (\ref{general}), we can write
\[L(\alpha,\varphi)=A(\alpha)\varphi(0)+G_1(\alpha)\varphi(-r_1)+G_2(\alpha)\varphi(-r_2)+\cdots+G_m(\alpha)\varphi(-r_m),\]
and thus  (\ref{B11}) can be calculated explicitly as follows
\[
\begin{array}{l}
B_{11}=\psi_1(0)\left(-\frac{n_1^2}{l^2}D_1^{(1,0)} \phi_1(0)+A_1^{(1,0)}\phi_1(0)+G_{11}^{(1,0)}\phi_1(-r_1)+\cdots+G_{m1}^{(1,0)}\phi_1(-r_m)\right), \\
B_{21}=\psi_1(0)\left(-\frac{n_1^2}{l^2}D_1^{(0,1)} \phi_1(0)+A_1^{(0,1)}\phi_1(0)+G_{11}^{(0,1)}\phi_1(-r_1)+\cdots+G_{m1}^{(0,1)}\phi_1(-r_m)\right),\\
B_{13}=\psi_3(0)\left(-\frac{n_2^2}{l^2}D_1^{(1,0)} \phi_3(0)+A_1^{(1,0)}\phi_3(0)+G_{11}^{(1,0)}\phi_3(-r_1)+\cdots+G_{m1}^{(1,0)}\phi_3(-r_m)\right), \\
B_{23}=\psi_3(0)\left(-\frac{n_2^2}{l^2}D_1^{(0,1)} \phi_3(0)+A_1^{(0,1)}\phi_3(0)+G_{11}^{(0,1)}\phi_3(-r_1)+\cdots+G_{m1}^{(0,1)}\phi_3(-r_m)\right).\\
\end{array}\]
\end{remark}

\subsection{Third Order Normal Form}
To find the third order normal form of the double Hopf singularity, neglecting the high order terms of the perturbation parameters and noticing the assumption $(H_2)$, we   have
\begin{equation*}\begin{array}{r}
{\rm Im}(M_3^1)^c={\rm span}\left\{
\left(\begin{array}{c}z_1^2z_2\\0\\0\\0\end{array} \right) ,
\left(\begin{array}{c}z_1z_3z_4\\0\\0\\0\end{array} \right) ,
\left(\begin{array}{c}0\\z_1z_2^2\\0\\0\end{array} \right) ,
\left(\begin{array}{c}0\\z_2z_3z_4\\0\\0\end{array} \right),
\left(\begin{array}{c}0\\0\\z_3^2z_4\\0\end{array} \right),\right. \\\left.
\left(\begin{array}{c}0\\0\\z_1z_2z_3\\0\end{array} \right) ,
\left(\begin{array}{c}0\\0\\0\\z_3z_4^2\end{array} \right),
\left(\begin{array}{c}0\\0\\0\\z_1z_2z_4\end{array} \right)
\right\}. \end{array}
\end{equation*}
According to \cite{Faria}, the normal form up to the third order is given by
\begin{equation*}
\dot{z}=Bz+\frac{1}{2!}g_2^1(z,0,\alpha)+\frac{1}{3!}g_3^1(z,0,0)+\cdots,
\end{equation*}
where $g_3^1(z,0,0)={\rm Proj}_{{\rm Ker}(M_3^1)}\overline{f}_3^1(z,0,0)$, with
\begin{equation}
\label{f31bar}\begin{array}{rl}
\overline{f}_3^1(z,0,0)=&f_3^1(z,0,0)+\frac{3}{2}[D_zf_2^1(z,0,0)U_2^1(z,0)\\&+D_yf_2^1(z,0,0)U_2^2(z,0)-D_zU_2^1(z,0)g_2^1(z,0,0)],\end{array}
\end{equation}
with $(U_2^1(z,\alpha),U_2^2(z,\alpha))\in V_j^{6}(\mathbb{C}^4)\times V_j^{6}(\mathcal{Q}_1)$ given by (\ref{Ujzalpha}).
After  the calculation of the second  order normal forms, we obtain the third order $\overline{f}_3^1(z,0,0)$ given by (\ref{f31bar}).
 Noticing that $g_2^1(z,0,0)=0$, we divide the computation of  the third term $g_3^1(z,0,0)={\rm Proj}_{{\rm Im}(M_3^1)^c}\overline{f}_3^1(z,0,0)$  into the following three parts
 \[{\rm Proj}_{{\rm Im}(M_3^1)^c}(f_3^1(z,0,0)),\]  \[{\rm Proj}_{{\rm Im}(M_3^1)^c}(D_zf_2^1(z,0,0)U_2^1(z,0)),\]  and \[{\rm Proj}_{{\rm Im}(M_3^1)^c}(D_yf_2^1(z,0,0)U_2^2(z,0)),\]
which will be calculated in the following part.

\textbf{Step 1} The calculation of ${\rm Proj}_{{\rm Im}(M_3^1)^c}(f_3^1(z,0,0))$.

Calculating the third order Fr\'{e}chet derivative  of $\widetilde{F}(0,\Phi z_x)$ as follows
\begin{equation*}
\widetilde{F}_3(0,\Phi z_x)=\sum_{q_1+q_2+q_3+q_4=3}F_{q_1q_2q_3q_4}\gamma_{n_1}^{q_1+q_2}(x)\gamma_{n_2}^{q_3+q_4}(x)z_1^{q_1}z_2^{q_2}z_3^{q_3}z_4^{q_4},
\end{equation*}
we have
\begin{equation*}
\begin{array}{l}
f_3^1(z,0,0)=\Psi(0)\left( \begin{array}{l}
\langle\widetilde{F}_3(0,\Phi z_x),\beta_{n_1}\rangle\\
\langle\widetilde{F}_3(0,\Phi z_x),\beta_{n_2}\rangle
\end{array}
\right) \\
=\Psi(0)\left(\begin{array}{c}
\sum\limits_{q_1+q_2+q_3+q_4=3}F_{q_1q_2q_3q_4}\int_0^{l\pi}\gamma_{n_1}^{q_1+q_2+1}(x)\gamma_{n_2}^{q_3+q_4}(x)dxz_1^{q_1}z_2^{q_2}z_3^{q_3}z_4^{q_4}\\
\sum\limits_{q_1+q_2+q_3+q_4=3}F_{q_1q_2q_3q_4}\int_0^{l\pi}\gamma_{n_1}^{q_1+q_2}(x)\gamma_{n_2}^{q_3+q_4+1}(x)dxz_1^{q_1}z_2^{q_2}z_3^{q_3}z_4^{q_4}\\
\end{array} \right). \\
\end{array}
\end{equation*}

Thus,
\begin{equation*}
\label{3!f31z00}
\frac{1}{3!}{\rm Proj}_{{\rm Im}(M_3^1)^c}(f_3^1(z,0,0))=\left(\begin{array}{c}
C_{2100}z_1^2z_2+C_{1011}z_1z_3z_4\\
\overline{C_{2100}}z_1z_2^2+\overline{C_{1011}}z_2z_3z_4\\
C_{0021}z_3^2z_4+C_{1110}z_1z_2z_3\\
\overline{C_{0021}}z_3z_4^2+\overline{C_{1110}}z_1z_2z_4\\
\end{array} \right) ,
\end{equation*}
where
\begin{equation}\label{C2100}
\begin{aligned}
C_{2100}=\frac{1}{6}\psi_1(0)F_{2100}\gamma_{40},~& C_{1011}=\frac{1}{6}\psi_1(0)F_{1011}\gamma_{22},\\C_{0021}=\frac{1}{6}\psi_3(0)F_{0021}\gamma_{04},~& C_{1110}=\frac{1}{6}\psi_3(0)F_{1110}\gamma_{22}.\\
\end{aligned}
\end{equation}
Here $\gamma_{ij}=\int_0^{l\pi}\gamma_{n_1}^i(x)\gamma_{n_2}^j(x)dx$, and
\[\begin{aligned}
&\int_0^{l\pi}\gamma_{n_k}^4(x)dx=\left\lbrace \begin{array}{cc}
\frac{1}{l\pi},n_k=0,\\
\frac{3}{2 l\pi},n_k\neq 0,
\end{array}   \right. (k=1,2),\\
&\int_0^{l\pi}\gamma_{n_1}^2(x)\gamma_{n_2}^2(x)dx=\left\lbrace \begin{array}{ll}
\frac{1}{l\pi},&n_1=0, n_2=0,\\
\frac{1}{ l\pi},&n_1=0, n_2\neq 0,\\
\frac{1}{ l\pi},&n_1\neq 0,n_2\neq 0, n_1\neq n_2,\\
\frac{3}{2 l\pi},&n_1=n_2\neq 0.
\end{array} \right.
\end{aligned}\]

\textbf{Step 2} ${\rm Proj}_{{\rm Im}(M_3^1)^c}(D_zf_2^1(z,0,0)U_2^1(z,0)). $

From (\ref{fj1fj2}), (\ref{Ut}), (\ref{F2}) and (\ref{F2Phi}), we can write $f_2^1(z,0,0)=(f_2^{1(1)},f_2^{1(2)},f_2^{1(3)},f_2^{1(4)})^T$ as
\begin{equation}
\label{f21z00}
\begin{array}{l}
f_2^1(z,0,0)=\Psi(0)\left( \begin{array}{l}
\langle F_2(0,\Phi z_x),\beta_{n_1}\rangle\\
\langle F_2(0,\Phi z_x),\beta_{n_2}\rangle
\end{array}
\right) \\
=\Psi(0)\left(\begin{array}{c}
\sum\limits_{q_1+q_2+q_3+q_4=2}F_{q_1q_2q_3q_4}\int_0^{l\pi}\gamma_{n_1}^{q_1+q_2+1}(x)\gamma_{n_2}^{q_3+q_4}(x)dxz_1^{q_1}z_2^{q_2}z_3^{q_3}z_4^{q_4}\\
\sum\limits_{q_1+q_2+q_3+q_4=2}F_{q_1q_2q_3q_4}\int_0^{l\pi}\gamma_{n_1}^{q_1+q_2}(x)\gamma_{n_2}^{q_3+q_4+1}(x)dxz_1^{q_1}z_2^{q_2}z_3^{q_3}z_4^{q_4}\\
\end{array} \right), \\
\end{array}
\end{equation}
where
\[\begin{aligned}
f_2^{1(1)}= \psi_1(0)(F_{2000}z_1^2\beta_{30}+F_{1100}z_1z_2\beta_{30}+F_{1010}z_1z_3\beta_{21}+F_{1001}z_1z_4\beta_{21}
+F_{0200}z_2^2\beta_{30}\\+F_{0110}z_2z_3\beta_{21}+F_{0101}z_2z_4\beta_{21}+F_{0020}z_3^2\beta_{12}+F_{0011}z_3z_4\beta_{12}+F_{0002}z_4^2\beta_{12}),
\end{aligned} \]
\[\begin{aligned}
f_2^{1(2)}= \psi_2(0)(F_{2000}z_1^2\beta_{30}+F_{1100}z_1z_2\beta_{30}+F_{1010}z_1z_3\beta_{21}+F_{1001}z_1z_4\beta_{21}
+F_{0200}z_2^2\beta_{30}\\+F_{0110}z_2z_3\beta_{21}+F_{0101}z_2z_4\beta_{21}+F_{0020}z_3^2\beta_{12}+F_{0011}z_3z_4\beta_{12}+F_{0002}z_4^2\beta_{12}),
\end{aligned} \]
\[\begin{aligned}
f_2^{1(3)}=\psi_3(0) (F_{2000}z_1^2\beta_{21}+F_{1100}z_1z_2\beta_{21}+F_{1010}z_1z_3\beta_{12}+F_{1001}z_1z_4\beta_{12}+F_{0200}z_2^2\beta_{21}\\+F_{0110}z_2z_3\beta_{12}+F_{0101}z_2z_4\beta_{12}+F_{0020}z_3^2\beta_{03}+F_{0011}z_3z_4\beta_{03}+F_{0002}z_4^2\beta_{03}),
\end{aligned}\]
\[\begin{aligned}
f_2^{1(4)}=\psi_4(0) (F_{2000}z_1^2\beta_{21}+F_{1100}z_1z_2\beta_{21}+F_{1010}z_1z_3\beta_{12}+F_{1001}z_1z_4\beta_{12}+F_{0200}z_2^2\beta_{21}\\+F_{0110}z_2z_3\beta_{12}+F_{0101}z_2z_4\beta_{12}+F_{0020}z_3^2\beta_{03}+F_{0011}z_3z_4\beta_{03}+F_{0002}z_4^2\beta_{03}),
\end{aligned}\]
and
\begin{equation}\label{f200011}
\begin{array}{cc}
f_{2000}^{1(1)}=\psi_1(0)F_{2000}\beta_{30},& f_{2000}^{1(2)}=\psi_2(0)F_{2000}\beta_{30},\\f_{2000}^{1(3)}=\psi_3(0)F_{2000}\beta_{21},&f_{2000}^{1(4)}=\psi_4(0)F_{2000}\beta_{21},\\
\end{array}
\end{equation}
\begin{equation}\label{f110011}
\begin{array}{cc}
f_{1100}^{1(1)}=\psi_1(0)F_{1100}\beta_{30},& f_{1100}^{1(2)}=\psi_2(0)F_{1100}\beta_{30},\\f_{1100}^{1(3)}=\psi_3(0)F_{1100}\beta_{21},&f_{1100}^{1(4)}=\psi_4(0)F_{1100}\beta_{21},\\
\end{array}
\end{equation}
\begin{equation}\label{f101011}
\begin{array}{cc}
f_{1010}^{1(1)}=\psi_1(0)F_{1010}\beta_{21},& f_{1010}^{1(2)}=\psi_2(0)F_{1010}\beta_{21},\\f_{1010}^{1(3)}=\psi_3(0)F_{1010}\beta_{12},&f_{1010}^{1(4)}=\psi_4(0)F_{1010}\beta_{12},\\
\end{array}
\end{equation}
\begin{equation}\label{f100111}
\begin{array}{cc}
f_{1001}^{1(1)}=\psi_1(0)F_{1001}\beta_{21},& f_{1001}^{1(2)}=\psi_2(0)F_{1001}\beta_{21},\\f_{1001}^{1(3)}=\psi_3(0)F_{1001}\beta_{12},&f_{1001}^{1(4)}=\psi_4(0)F_{1001}\beta_{12},\\
\end{array}
\end{equation}
\begin{equation}\label{f020011}
\begin{array}{cc}
f_{0200}^{1(1)}=\psi_1(0)F_{0200}\beta_{30},& f_{0200}^{1(2)}=\psi_2(0)F_{0200}\beta_{30},\\f_{0200}^{1(3)}=\psi_3(0)F_{0200}\beta_{21},&f_{0200}^{1(4)}=\psi_4(0)F_{0200}\beta_{21},\\
\end{array}
\end{equation}
\begin{equation}\label{f011011}
\begin{array}{cc}
f_{0110}^{1(1)}=\psi_1(0)F_{0110}\beta_{21},& f_{0110}^{1(2)}=\psi_2(0)F_{0110}\beta_{21},\\f_{0110}^{1(3)}=\psi_3(0)F_{0110}\beta_{12},&f_{0110}^{1(4)}=\psi_4(0)F_{0110}\beta_{12},\\
\end{array}
\end{equation}
\begin{equation}\label{f010111}
\begin{array}{cc}
f_{0101}^{1(1)}=\psi_1(0)F_{0101}\beta_{21},& f_{0101}^{1(2)}=\psi_2(0)F_{0101}\beta_{21},\\f_{0101}^{1(3)}=\psi_3(0)F_{0101}\beta_{12},&f_{0101}^{1(4)}=\psi_4(0)F_{0101}\beta_{12},\\
\end{array}
\end{equation}
\begin{equation}\label{f002011}
\begin{array}{cc}
f_{0020}^{1(1)}=\psi_1(0)F_{0020}\beta_{12},& f_{0020}^{1(2)}=\psi_2(0)F_{0020}\beta_{12},\\f_{0020}^{1(3)}=\psi_3(0)F_{0020}\beta_{03},&f_{0020}^{1(4)}=\psi_4(0)F_{0020}\beta_{03},\\
\end{array}
\end{equation}
\begin{equation}\label{f001111}
\begin{array}{cc}
f_{0011}^{1(1)}=\psi_1(0)F_{0011}\beta_{12},& f_{0011}^{1(2)}=\psi_2(0)F_{0011}\beta_{12},\\f_{0011}^{1(3)}=\psi_3(0)F_{0011}\beta_{03},&f_{0011}^{1(4)}=\psi_4(0)F_{0011}\beta_{03},\\
\end{array}
\end{equation}
\begin{equation}\label{f000211}
\begin{array}{cc}
f_{0002}^{1(1)}=\psi_1(0)F_{0002}\beta_{12},& f_{0002}^{1(2)}=\psi_2(0)F_{0002}\beta_{12},\\f_{0002}^{1(3)}=\psi_3(0)F_{0002}\beta_{03},&f_{0002}^{1(4)}=\psi_4(0)F_{0002}\beta_{03},\\

\end{array}
\end{equation}
with $\beta_{ij}=\int_0^{l\pi}\gamma_{n_1}^i(x)\gamma_{n_2}^j(x)dx$, and
\[\begin{aligned}
&\int_0^{l\pi}\gamma_{n_k}^3(x)dx=\left\lbrace \begin{array}{ll}
\sqrt{\frac{1}{l\pi}},& n_k=0,\\
0,& n_k\neq 0,
\end{array}   ~(k=1,2),\right.\\
&\int_0^{l\pi}\gamma_{n_1}^2(x)\gamma_{n_2}(x)dx=\left\lbrace \begin{array}{cc}
\sqrt{\frac{1}{l\pi}},& n_1=0,n_2=0,\\
0,& n_1=0,  n_2\neq 0,\\
\sqrt{\frac{1}{ 2l\pi}},& n_1\neq 0,n_2\neq 0, n_2=2n_1,\\
0,& n_1\neq 0,n_2\neq 0, n_2\neq 2n_1,
\end{array} \right.\\
&\int_0^{l\pi}\gamma_{n_1}(x)\gamma_{n_2}^2(x)dx=\left\lbrace \begin{array}{cc}
\sqrt{\frac{1}{l\pi}},& n_1=0,n_2=0,\\
\sqrt{\frac{1}{l\pi}},& n_1=0,  n_2\neq 0,\\
0,& n_1\neq 0,n_2\neq 0.
\end{array}\right.
\end{aligned}\]
Combining with (\ref{Mj1}) and (\ref{f21z00}), we can calculate $U_2^1(z,0)=(U_2^{1(1)},U_2^{1(2)},U_2^{1(3)},U_2^{1(4)})^T$ from the following formula
\begin{equation*}
U_2^1(z,0)=(M_2^1)^{-1}{\rm Proj}_{{\rm Im}(M_2^1)}f_2^1(z,0,0).
\end{equation*}
Then, we have
\[\begin{aligned}
U_2^{1(1)}=&\frac{1}{i\omega_1}f_{2000}^{1(1)}z_1^2+\frac{1}{-i\omega_1}f_{1100}^{1(1)}z_1z_2+\frac{1}{i\omega_2}f_{1010}^{1(1)}z_1z_3+\frac{1}{-i\omega_2}f_{1001}^{1(1)}z_1z_4\\&+\frac{1}{-3i\omega_1}f_{0200}^{1(1)}z_2^2+\frac{1}{-2i\omega_1+i\omega_2}f_{0110}^{1(1)}z_2z_3+\frac{1}{-2i\omega_1-i\omega_2}f_{0101}^{1(1)}z_2z_4\\&+\frac{1}{2i\omega_2-i\omega_1}f_{0020}^{1(1)}z_3^2+\frac{1}{-i\omega_1}f_{0011}^{1(1)}z_3z_4+\frac{1}{-2i\omega_2-i\omega_1}f_{0002}^{1(1)}z_4^2 ,
\end{aligned}\]
\[\begin{aligned}
U_2^{1(2)}=&\frac{1}{3i\omega_1}f_{2000}^{1(2)}z_1^2+\frac{1}{i\omega_1}f_{1100}^{1(2)}z_1z_2+\frac{1}{2i\omega_1+i\omega_2}f_{1010}^{1(2)}z_1z_3\\&+\frac{1}{2i\omega_1-i\omega_2}f_{1001}^{1(2)}z_1z_4+\frac{1}{-i\omega_1}f_{0200}^{1(2)}z_2^2+\frac{1}{i\omega_2}f_{0110}^{1(2)}z_2z_3+\frac{1}{-i\omega_2}f_{0101}^{1(2)}z_2z_4\\&+\frac{1}{2i\omega_2+i\omega_1}f_{0020}^{1(2)}z_3^2+\frac{1}{i\omega_1}f_{0011}^{1(2)}z_3z_4+\frac{1}{-2i\omega_2+i\omega_1}f_{0002}^{1(2)}z_4^2,
\end{aligned}\]
\[\begin{aligned}
U_2^{1(3)}&=\frac{1}{2i\omega_1-i\omega_2}f_{2000}^{1(3)}z_1^2+\frac{1}{-i\omega_2}f_{1100}^{1(3)}z_1z_2+\frac{1}{i\omega_1}f_{1010}^{1(3)}z_1z_3+\frac{1}{i\omega_1-2i\omega_2}f_{1001}^{1(3)}z_1z_4\\&+\frac{1}{-2i\omega_1-i\omega_2}f_{0200}^{1(3)}z_2^2+\frac{1}{-i\omega_1}f_{0110}^{1(3)}z_2z_3+\frac{1}{-i\omega_1-2i\omega_2}f_{0101}^{1(3)}z_2z_4\\&+\frac{1}{i\omega_2}f_{0020}^{1(3)}z_3^2+\frac{1}{-i\omega_2}f_{0011}^{1(3)}z_3z_4+\frac{1}{-3i\omega_2}f_{0002}^{1(3)}z_4^2,
\end{aligned}\]
\[\begin{aligned}
U_2^{1(4)}&=\frac{1}{2i\omega_1+i\omega_2}f_{2000}^{1(4)}z_1^2+\frac{1}{i\omega_2}f_{1100}^{1(4)}z_1z_2+\frac{1}{i\omega_1+2i\omega_2}f_{1010}^{1(4)}z_1z_3+\frac{1}{i\omega_1}f_{1001}^{1(4)}z_1z_4\\&+\frac{1}{-2i\omega_1+i\omega_2}f_{0200}^{1(4)}z_2^2+\frac{1}{-i\omega_1+2i\omega_2}f_{0110}^{1(4)}z_2z_3+\frac{1}{-i\omega_1}f_{0101}^{1(4)}z_2z_4\\&+\frac{1}{3i\omega_2}f_{0020}^{1(4)}z_3^2+\frac{1}{i\omega_2}f_{0011}^{1(4)}z_3z_4+\frac{1}{-i\omega_2}f_{0002}^{1(4)}z_4^2.
\end{aligned}\]
Thus,
\begin{equation}
\label{3!Dzf21z00U21}
\frac{1}{3!}{\rm Proj}_{{\rm Im}(M_3^1)^c}(D_zf_2^1(z,0,0)U_2^1(z,0))=\left(\begin{array}{c}
D_{2100}z_1^2z_2+D_{1011}z_1z_3z_4\\
\overline{D_{2100}}z_1z_2^2+\overline{D_{1011}}z_2z_3z_4\\
D_{0021}z_3^2z_4+D_{1110}z_1z_2z_3\\
\overline{D_{0021}}z_3z_4^2+\overline{D_{1110}}z_1z_2z_4\\
\end{array} \right) ,
\end{equation}
where
\begin{equation}\label{D2100}
\begin{aligned}
D_{2100}=&\frac{1}{6}(2f_{2000}^{1(1)}\frac{1}{-i\omega_1}f_{1100}^{1(1)}+f_{1100}^{1(1)}\frac{1}{i\omega_1}f_{2000}^{1(1)}+f_{1100}^{1(1)}\frac{1}{i\omega_1}f_{1100}^{1(2)}\\&+2f_{0200}^{1(1)}\frac{1}{3i\omega_1}f_{2000}^{1(2)}+f_{1010}^{1(1)}\frac{1}{-i\omega_2}f_{1100}^{1(3)}+f_{0110}^{1(1)}\frac{1}{2i\omega_1-i\omega_2}f_{2000}^{1(3)}\\&+f_{1001}^{1(1)}\frac{1}{i\omega_2}f_{1100}^{1(4)}+f_{0101}^{1(1)}\frac{1}{2i\omega_1+i\omega_2}f_{2000}^{1(4)}),\\
\end{aligned}
\end{equation}
\begin{equation}\label{D1011}
\begin{aligned}
D_{1011}=&\frac{1}{6}(2f_{2000}^{1(1)}\frac{1}{-i\omega_1}f_{0011}^{1(1)}+f_{1010}^{1(1)}\frac{1}{-i\omega_2}f_{1001}^{1(1)}+f_{1001}^{1(1)}\frac{1}{i\omega_2}f_{1010}^{1(1)}\\&+f_{1100}^{1(1)}\frac{1}{i\omega_1}f_{0011}^{1(2)}+f_{0110}^{1(1)}\frac{1}{2i\omega_1-i\omega_2}f_{1001}^{1(2)}+f_{0101}^{1(1)}\frac{1}{2i\omega_1+i\omega_2}f_{1010}^{1(2)}\\&+f_{1010}^{1(1)}\frac{1}{-i\omega_2}f_{0011}^{1(3)}+2f_{0020}^{1(1)}\frac{1}{i\omega_1-2i\omega_2}f_{1001}^{1(3)}+f_{0011}^{1(1)}\frac{1}{i\omega_1}f_{1010}^{1(3)}\\&+f_{1001}^{1(1)}\frac{1}{i\omega_2}f_{0011}^{1(4)}+f_{0011}^{1(1)}\frac{1}{i\omega_1}f_{1001}^{1(4)}+2f_{0002}^{1(1)}\frac{1}{i\omega_1+2i\omega_2}f_{1010}^{1(4)}),\\
\end{aligned}
\end{equation}
\begin{equation}\label{D0021}
\begin{aligned}
D_{0021}&=\frac{1}{6}(f_{1010}^{1(3)}\frac{1}{-i\omega_1}f_{0011}^{1(1)}+f_{1001}^{1(3)}\frac{1}{2i\omega_2-i\omega_1}f_{0020}^{1(1)}+f_{0110}^{1(3)}\frac{1}{i\omega_1}f_{0011}^{1(2)}\\&+f_{0101}^{1(3)}\frac{1}{2i\omega_2+i\omega_1}f_{0020}^{1(2)}+2f_{0020}^{1(3)}\frac{1}{-i\omega_2}f_{0011}^{1(3)}+f_{0011}^{1(3)}\frac{1}{i\omega_2}f_{0020}^{1(3)}\\&+f_{0011}^{1(3)}\frac{1}{i\omega_2}f_{0011}^{1(4)}+2f_{0002}^{1(3)}\frac{1}{3i\omega_2}f_{0020}^{1(4)}),\\
\end{aligned}
\end{equation}
\begin{equation}\label{D1110}
\begin{aligned}
D_{1110}&=\frac{1}{6}(2f_{2000}^{1(3)}\frac{1}{-2i\omega_1+i\omega_2}f_{0110}^{1(1)}+f_{1100}^{1(3)}\frac{1}{i\omega_2}f_{1010}^{1(1)}+f_{1010}^{1(3)}\frac{1}{-i\omega_1}f_{1100}^{1(1)}\\&+f_{1100}^{1(3)}\frac{1}{i\omega_2}f_{0110}^{1(2)}+2f_{0200}^{1(3)}\frac{1}{2i\omega_1+i\omega_2}f_{1010}^{1(2)}+f_{0110}^{1(3)}\frac{1}{i\omega_1}f_{1100}^{1(2)}\\&+f_{1010}^{1(3)}\frac{1}{-i\omega_1}f_{0110}^{1(3)}+f_{0110}^{1(3)}\frac{1}{i\omega_1}f_{1010}^{1(3)}+2f_{0020}^{1(3)}\frac{1}{-i\omega_2}f_{1100}^{1(3)}\\&+f_{1001}^{1(3)}\frac{1}{-i\omega_1+2i\omega_2}f_{0110}^{1(4)}+f_{0101}^{1(3)}\frac{1}{i\omega_1+2i\omega_2}f_{1010}^{1(4)}+f_{0011}^{1(3)}\frac{1}{i\omega_2}f_{1100}^{1(4)}).\\
\end{aligned}
\end{equation}

\textbf{Step 3.} The calculation of  $ {\rm Proj}_{{\rm Im}(M_3^1)^c}(D_yf_2^1(z,0,0)U_2^2(z,0))$.

First, we will calculate the Fr\'{e}chet derivative $D_yf_2^1(z,0,0):\mathcal{Q}_1\rightarrow X_{\mathbb{C}}$. From (\ref{F2}) and (\ref{F2Phi}), $\widetilde{F}_2(z,y,0)$ can be written as
\begin{equation}\label{S2}
\begin{array}{rl}&\widetilde{F}_2(z,y,0)=S_2(\Phi z_x,y)+\textit{o}(z^2,y^2)\\&=S_{yz_1}(y)z_1\gamma_{n_1}+S_{yz_2}(y)z_2\gamma_{n_1}+S_{yz_3}(y)z_3\gamma_{n_2}+S_{yz_4}(y)z_4\gamma_{n_2}+\textit{o}(z^2,y^2),\end{array}
\end{equation}
where
$S_{yz_i} (i=1,2,3,4)$ are linear operators from $\mathcal{Q}_1$ to $X_{\mathbb{C}}$.
\begin{remark}\label{S}
Again if there are $m$ discrete  delays in the system (\ref{general}), we can get the explicit formula
\begin{equation*}
\begin{aligned}
S_{yz_i}(\varphi)&=(F_{y_1(0)z_i},F_{y_2(0)z_i},\cdots,F_{y_n(0)z_i})\varphi(0)\\&+(F_{y_1(-r_1)z_i},F_{y_2(-r_1)z_i},\cdots,F_{y_n(-r_1)z_i})\varphi(-r_1)+\cdots\\&+(F_{y_1(-r_m)z_i},F_{y_2(-r_m)z_i},\cdots,F_{y_n(-r_m)z_i})\varphi(-r_m).
\end{aligned}
\end{equation*}
\end{remark}

Now, we have
\[D_y\widetilde{F}_2(z,0,0)(\varphi)=S_{yz_1}(\varphi)z_1\gamma_{n_1}+S_{yz_2}(\varphi)z_2\gamma_{n_1}+S_{yz_3}(\varphi)z_3\gamma_{n_2}+S_{yz_4}(\varphi)z_4\gamma_{n_2}.\]

Let
$U_2^2(z,0)=h(z)=\sum\limits_{j\geq 0}h_j(z)\cdot\beta_j(x)=\sum\limits_{j\geq 0}h_j(z)\gamma_j(x)$
with
\[h_j(z)=\left(\begin{array}{c}
h_j^{(1)}(z)\\h_j^{(2)}(z)\\\vdots\\h_j^{(n)}(z)\\
\end{array} \right)=\sum_{q_1+q_2+q_3+q_4=2}\left(\begin{array}{c}
h_{jq_1q_2q_3q_4}^{(1)}(z)\\h_{jq_1q_2q_3q_4}^{(2)}(z)\\\vdots\\h_{jq_1q_2q_3q_4}^{(n)}(z)\\
\end{array} \right)z_1^{q_1}z_2^{q_2}z_3^{q_3}z_4^{q_4} .\]
Thus,
\[\begin{array}{l}
D_yf_2^1(z,0,0)(U_2^2(z,0))=\Psi(0)\left(\begin{array}{c}
\langle D_y\widetilde{F}_2(z,0,0)(U_2^2(z,0)),\beta_{n_1}\rangle\\
\langle D_y\widetilde{F}_2(z,0,0)(U_2^2(z,0)),\beta_{n_2}\rangle\\
\end{array} \right) \\
=\Psi(0)\left( \begin{array}{c}
\langle S_{yz_1} (\sum\limits_{j\geq 0}h_j\gamma_j)\gamma_{n_1},\beta_{n_1}\rangle z_1+\langle S_{yz_2} (\sum\limits_{j\geq 0}h_j\gamma_j)\gamma_{n_1},\beta_{n_1}\rangle z_2\\+\langle S_{yz_3} (\sum\limits_{j\geq 0}h_j\gamma_j)\gamma_{n_2},\beta_{n_1}\rangle z_3+\langle S_{yz_4} (\sum\limits_{j\geq 0}h_j\gamma_j)\gamma_{n_2},\beta_{n_1}\rangle z_4 \\
\langle S_{yz_1} (\sum\limits_{j\geq 0}h_j\gamma_j)\gamma_{n_1},\beta_{n_2}\rangle z_1+\langle S_{yz_2} (\sum\limits_{j\geq 0}h_j\gamma_j)\gamma_{n_1},\beta_{n_2}\rangle z_2\\+\langle S_{yz_3} (\sum\limits_{j\geq 0}h_j\gamma_j)\gamma_{n_2},\beta_{n_2}\rangle z_3+\langle S_{yz_4} (\sum\limits_{j\geq 0}h_j\gamma_j)\gamma_{n_2},\beta_{n_2}\rangle z_4 \\
\end{array}\right) \\
=\Psi(0)\left( \begin{array}{c}
\sum\limits_{j\geq 0}[b_{jn_1n_1} S_{yz_1}(h_j)z_1+b_{jn_1n_1}S_{yz_2}( h_j)z_2\\+b_{jn_2n_1}S_{yz_3}( h_j)z_3+b_{jn_2n_1}S_{yz_4}( h_j)z_4]  \\
\sum\limits_{j\geq 0} [ b_{jn_1n_2}S_{yz_1}( h_j)z_1+ b_{jn_1n_2}S_{yz_2}( h_j)z_2\\+b_{jn_2n_2}S_{yz_3}( h_j)z_3+b_{jn_2n_2}S_{yz_4}( h_j)z_4]  \\
\end{array}\right). \\
\end{array}\]

To give a clear expression of our derivation about the normal form, the rest part of derivation are divided into three cases: $n_1=n_2=0$,  $n_1= 0, n_2 \neq 0$, and $n_1\neq 0, n_2 \neq 0$.

\textbf{Case 1}
When $n_1=n_2=0$, in fact
\begin{equation*}
b_{jn_1n_1}=b_{jn_1n_2}=b_{jn_2n_1}=b_{jn_2n_2}=\int_0^{l\pi}\gamma_j(x)\gamma_0(x)\gamma_0(x)dx=\left\lbrace \begin{array}{ll}
\frac{1}{\sqrt{l\pi}}, &j=0,\\
0,&j \neq 0.\\
\end{array}\right.
\end{equation*}
Then, obviously
\[\begin{array}{l}
D_yf_2^1(z,0,0)(U_2^2(z,0))\\
=\frac{1}{\sqrt{l\pi}}\Psi(0)\left( \begin{array}{c}
 S_{yz_1}(h_0)z_1+S_{yz_2}( h_0)z_2+S_{yz_3}( h_0)z_3+S_{yz_4}( h_0)z_4 \\
 S_{yz_1}( h_0)z_1+ S_{yz_2}( h_0)z_2+S_{yz_3}( h_0)z_3+S_{yz_4}( h_0)z_4  \\
\end{array}\right) \\
=\frac{1}{\sqrt{l\pi}}\left(\begin{array}{c}
\psi_1(0)(S_{yz_1}(h_0)z_1+S_{yz_2}( h_0)z_2+S_{yz_3}( h_0)z_3+S_{yz_4}( h_0)z_4)\\
\psi_2(0)(S_{yz_1}(h_0)z_1+S_{yz_2}( h_0)z_2+S_{yz_3}( h_0)z_3+S_{yz_4}( h_0)z_4)\\
\psi_3(0)( S_{yz_1}( h_0)z_1+ S_{yz_2}( h_0)z_2+S_{yz_3}( h_0)z_3+S_{yz_4}( h_0)z_4)\\
\psi_4(0)( S_{yz_1}( h_0)z_1+ S_{yz_2}( h_0)z_2+S_{yz_3}( h_0)z_3+S_{yz_4}( h_0)z_4)\\
\end{array} \right) .
\end{array}\]
Thus,
\begin{equation*}
\frac{1}{3!}{\rm Proj} _{{\rm Im}(M_3^1)^c}(D_yf_2^1(z,0,0)(U_2^2(z,0)))\\
=\left( \begin{array}{c}
 E_{2100}z_1^2z_2+E_{1011}z_1z_3z_4 \\
 \overline{E_{2100}}z_1z_2^2+\overline{E_{1011}}z_2z_3z_4 \\
 E_{0021}z_3^2z_4+E_{1110}z_1z_2z_3  \\
 \overline{ E_{0021}}z_3z_4^2+\overline{E_{1110}}z_1z_2z_4  \\
  \end{array}\right) ,\\
\end{equation*}
where
\begin{equation}
\label{3!Dyf21z00U22(1)}
\begin{array}{l}
E_{2100}=\frac{1}{6\sqrt{l\pi}}\psi_1(0)\left[ S_{yz_1}(h_{01100})+S_{yz_2}( h_{02000})\right], \\
E_{1011}=\frac{1}{6\sqrt{l\pi}}\psi_1(0)\left[ S_{yz_1}(h_{00011})+S_{yz_3}( h_{01001})+S_{yz_4}( h_{01010})\right], \\
E_{0021}=\frac{1}{6\sqrt{l\pi}}\psi_3(0)\left[ S_{yz_3}( h_{00011})+S_{yz_4}( h_{00020})\right], \\
E_{1110}=\frac{1}{6\sqrt{l\pi}}\psi_3(0)\left[ S_{yz_1}( h_{00110})+S_{yz_2}( h_{01010})+S_{yz_3}( h_{01100})\right] .
\end{array}
\end{equation}
Clearly, we still need to calculate $h_{01100}$, $h_{02000}$, $h_{00011}$, $h_{01001}$, $h_{01010}$, $h_{00020}$, and $h_{00110}$.

From (\ref{mjU}),   (\ref{A}) and (\ref{zdoty}), we have
\[\begin{array}{l}
M_2^2U_2^2(z,0)(\theta)=M_2^2h(z)(\theta)=D_zh(z)Bz-A_{1}(h(z))\\
=\left\lbrace \begin{array}{ll}
D_zh(z)Bz-D_0\Delta h(0)-L_0(h(z)),&\theta=0,\\
D_zh(z)Bz-D_{\theta}h(z),&\theta\neq 0,
\end{array}\right.\\
=\left\lbrace \begin{array}{ll}
\sum\limits_{j\geq 0}[D_zh_j(z)\gamma_j(x)Bz-D_0\Delta h_j(0)\gamma_j(x)-L_0(h_j(z)\gamma_j(x))],&\theta=0,\\
\sum\limits_{j\geq 0}[D_zh_j(z)\gamma_j(x)Bz-D_{\theta}h_j(z)\gamma_j(x)],&\theta\neq 0.
\end{array}\right.\\
\end{array}\]

According to (\ref{fj1fj2}),
\[\begin{aligned}
& f_2^2(z,0,0)\\&=\left\lbrace \begin{array}{ll}
\begin{array}{l}
\widetilde{F}_2(z,0,0)-\phi_1(0)f_2^{1(1)}(z,0,0)\gamma_{n_1}-\phi_2(0)f_2^{1(2)}(z,0,0)\gamma_{n_1}\\~~~-\phi_3(0)f_2^{1(3)}(z,0,0)\gamma_{n_2}-\phi_4(0)f_2^{1(4)}(z,0,0)\gamma_{n_2},
\end{array}&\theta=0,\\
\begin{array}{l}
-\phi_1(\theta)f_2^{1(1)}(z,0,0)\gamma_{n_1}-\phi_2(\theta)f_2^{1(2)}(z,0,0)\gamma_{n_1}\\~~~-\phi_3(\theta)f_2^{1(3)}(z,0,0)\gamma_{n_2}-\phi_4(\theta)f_2^{1(4)}(z,0,0)\gamma_{n_2},
\end{array}&\theta\neq 0.
\end{array}\right.
\end{aligned}
\]

From (\ref{fj1fj2}), (\ref{zyalpha}), and (\ref{Ujzalpha}), we have
\begin{equation}\label{M22f22}
\left\langle M_2^2(U_2^2(z,0)),\beta_j\right\rangle =\left\langle f_2^2(z,0,0),\beta_j\right\rangle.
\end{equation}
Matching the coefficients of $z_1^{q_1}z_2^{q_2}z_3^{q_3}z_4^{q_4}$ in (\ref{M22f22})  when $j=0$,
we can get the results of $h_{01100}$, $h_{02000}$, $h_{00011}$, $h_{01001}$, $h_{01010}$, $h_{00020}$, and $h_{00110}$. We take $h_{02000}(\theta)$ as an example, and the others can be calculated in the same method.

When $\theta\neq 0$, solve the following equation
\[\begin{aligned}
2i\omega_1h_{02000}(\theta)-\dot{h}_{02000}(\theta)=& -\langle\phi_1(\theta)f_2^{1(1)}\gamma_{n_1},\beta_{0}\rangle-\langle\phi_2(\theta)f_2^{1(2)}\gamma_{n_1},\beta_{0}\rangle\\&-\langle\phi_3(\theta)f_2^{1(3)}\gamma_{n_2},\beta_{0}\rangle-\langle\phi_4(\theta)f_2^{1(4)}\gamma_{n_2},\beta_{0}\rangle,
\end{aligned}\]
and we get
\[\begin{aligned}
h_{02000}(\theta)&=e^{2i\omega_1\theta}h_{02000}(0)+\frac{1}{-i\omega_1}(e^{i\omega_1\theta}-e^{2i\omega_1\theta})\phi_1(0)f_{2000}^{1(1)}\\
&+\frac{1}{-3i\omega_1}(e^{-i\omega_1\theta}-e^{2i\omega_1\theta})\phi_2(0)f_{2000}^{1(2)}+\frac{1}{-2i\omega_1+i\omega_2}(e^{i\omega_2\theta}-e^{2i\omega_1\theta})\phi_3(0)f_{2000}^{1(3)}\\&+\frac{1}{-2i\omega_1-i\omega_2}(e^{-i\omega_2\theta}-e^{2i\omega_1\theta})\phi_4(0)f_{2000}^{1(4)}.
\end{aligned}\]
When $\theta=0$,
\[\begin{aligned}
2i\omega_1h_{02000}(0)-D_0\Delta h_{02000}(0)-L_0(h_{02000})
=\langle F_{2000}\gamma_{n_1}^2,\beta_{0}\rangle -\phi_1(0)f_2^{1(1)}\langle\gamma_{n_1},\beta_{0}\rangle\\-\phi_2(0)f_2^{1(2)}\langle\gamma_{n_1},\beta_{0}\rangle-\phi_3(0)f_2^{1(3)}\langle\gamma_{n_2},\beta_{0}\rangle-\phi_4(0)f_2^{1(4)}\langle\gamma_{n_2},\beta_{0}\rangle.
\end{aligned}\]
Combining with
\[[i\omega_1-\frac{n_1^2}{l^2}D_0-L_0(e^{i\omega_1\cdot}I_d)]\phi_1(0)=0,~~~[-i\omega_1-\frac{n_1^2}{l^2}D_0-L_0(e^{-i\omega_1\cdot}I_d)]\phi_2(0)=0,\]
\[[i\omega_2-\frac{n_2^2}{l^2}D_0-L_0(e^{i\omega_2\cdot}I_d)]\phi_3(0)=0,~~~[-i\omega_2-\frac{n_2^2}{l^2}D_0-L_0(e^{-i\omega_2\cdot}I_d)]\phi_4(0)=0,\]
we can get
\[\begin{aligned}
h_{02000}(0)&=\frac{1}{-i\omega_1}\phi_1(0)f_{2000}^{1(1)}+\frac{1}{-3i\omega_1}\phi_2(0)f_{2000}^{1(2)}+\frac{1}{-2i\omega_1+i\omega_2}\phi_3(0)f_{2000}^{1(3)}\\&
+\frac{1}{-2i\omega_1-i\omega_2}\phi_4(0)f_{2000}^{1(4)}-[-2i\omega_1+L_0(e^{2i\omega_1\cdot}I_d)]^{-1}\langle F_{2000}\gamma_{n_1}^2,\beta_0\rangle,
\end{aligned}\]
and
\begin{equation*}
\begin{aligned}
h_{02000}(\theta)&=\frac{1}{-i\omega_1}\phi_1(\theta)f_{2000}^{1(1)}+\frac{1}{-3i\omega_1}\phi_2(\theta)f_{2000}^{1(2)}
+\frac{1}{-2i\omega_1+i\omega_2}\phi_3(\theta)f_{2000}^{1(3)}\\&+\frac{1}{-2i\omega_1-i\omega_2}\phi_4(\theta)f_{2000}^{1(4)}- \frac{1}{\sqrt{l\pi}}e^{2i\omega_1\theta}[-2i\omega_1+L_0(e^{2i\omega_1\cdot}I_d)]^{-1}F_{2000}.
\end{aligned}
\end{equation*}

Using the same method mentioned above, we can work out
\begin{equation*}
\begin{aligned}
h_{01100}(\theta)&=\frac{1}{i\omega_1}\phi_1(\theta)f_{1100}^{1(1)}+\frac{1}{-i\omega_1}\phi_2(\theta)f_{1100}^{1(2)}\\&
+\frac{1}{i\omega_2}\phi_3(\theta)f_{1100}^{1(3)}+\frac{1}{-i\omega_2}\phi_4(\theta)f_{1100}^{1(4)}- \frac{1}{\sqrt{l\pi}}[L_0(I_d)]^{-1}F_{1100},
\end{aligned}
\end{equation*}
\begin{equation*}
\begin{aligned}
h_{00011}(\theta)&=\frac{1}{i\omega_1}\phi_1(\theta)f_{0011}^{1(1)}+\frac{1}{-i\omega_1}\phi_2(\theta)f_{0011}^{1(2)}\\&
+\frac{1}{i\omega_2}\phi_3(\theta)f_{0011}^{1(3)}+\frac{1}{-i\omega_2}\phi_4(\theta)f_{0011}^{1(4)}- \frac{1}{\sqrt{l\pi}}[L_0(I_d)]^{-1}F_{0011},
\end{aligned}
\end{equation*}
\begin{equation*}
\begin{aligned}
h_{01001}(\theta)&=\frac{1}{i\omega_2}\phi_1(\theta)f_{1001}^{1(1)}+\frac{1}{-2i\omega_1+i\omega_2}\phi_2(\theta)f_{1001}^{1(2)}\\&
+\frac{1}{-i\omega_1+2i\omega_2}\phi_3(\theta)f_{1001}^{1(3)}+\frac{1}{-i\omega_1}\phi_4(\theta)f_{1001}^{1(4)}\\&- \frac{1}{\sqrt{l\pi}}e^{(i\omega_1-i\omega_2)\theta}[-(i\omega_1-i\omega_2)+L_0(e^{(i\omega_1-i\omega_2)\cdot}I_d)]^{-1}F_{1001},
\end{aligned}
\end{equation*}
\begin{equation*}
\begin{aligned}
h_{01010}(\theta)&=\frac{1}{-i\omega_2}\phi_1(\theta)f_{1010}^{1(1)}+\frac{1}{-2i\omega_1-i\omega_2}\phi_2(\theta)f_{1010}^{1(2)}\\&
+\frac{1}{-i\omega_1}\phi_3(\theta)f_{1010}^{1(3)}+\frac{1}{-i\omega_1-2i\omega_2}\phi_4(\theta)f_{1010}^{1(4)}\\&- \frac{1}{\sqrt{l\pi}}e^{(i\omega_1+i\omega_2)\theta}[-(i\omega_1+i\omega_2)+L_0(e^{(i\omega_1+i\omega_2)\cdot}I_d)]^{-1}F_{1010},
\end{aligned}
\end{equation*}
\begin{equation*}
\begin{aligned}
h_{00110}(\theta)&=\frac{1}{2i\omega_1-i\omega_2}\phi_1(\theta)f_{0110}^{1(1)}+\frac{1}{-i\omega_2}\phi_2(\theta)f_{0110}^{1(2)}\\&
+\frac{1}{i\omega_1}\phi_3(\theta)f_{0110}^{1(3)}+\frac{1}{i\omega_1-2i\omega_2}\phi_4(\theta)f_{0110}^{1(4)}\\&- \frac{1}{\sqrt{l\pi}}e^{(-i\omega_1+i\omega_2)\theta}[-(-i\omega_1+i\omega_2)+L_0(e^{(-i\omega_1+i\omega_2)\cdot}I_d)]^{-1}F_{0110},
\end{aligned}
\end{equation*}
and
\begin{equation*}
\begin{aligned}
h_{00020}(\theta)&=\frac{1}{i\omega_1-2i\omega_2}\phi_1(\theta)f_{0020}^{1(1)}+\frac{1}{-i\omega_1-2i\omega_2}\phi_2(\theta)f_{0020}^{1(2)}
+\frac{1}{-i\omega_2}\phi_3(\theta)f_{0020}^{1(3)}\\&+\frac{1}{-3i\omega_2}\phi_4(\theta)f_{0020}^{1(4)}- \frac{1}{\sqrt{l\pi}}e^{2i\omega_2\theta}[-2i\omega_2+L_0(e^{2i\omega_2\cdot}I_d)]^{-1}F_{0020}.
\end{aligned}
\end{equation*}

\textbf{Case 2}
When $n_1=0,n_2\neq 0$, we have
\begin{equation*}
\begin{array}{l}
b_{jn_1n_1}=\int_0^{l\pi}\gamma_j(x)\gamma_0(x)\gamma_0(x)dx=\left\lbrace \begin{array}{ll}
\frac{1}{\sqrt{l\pi}}, &j=0,\\
0,&j \neq 0,\\
\end{array}\right.\\
b_{jn_1n_2}=b_{jn_2n_1}=\int_0^{l\pi}\gamma_j(x)\gamma_0(x)\gamma_{n_2}(x)dx=\left\lbrace \begin{array}{ll}
\frac{1}{\sqrt{l\pi}}, &j=n_2,\\
0, &j\neq n_2,
\end{array}\right.\\
b_{jn_2n_2}=\int_0^{l\pi}\gamma_j(x)\gamma_{n_2}(x)\gamma_{n_2}(x)dx=\left\lbrace \begin{array}{ll}
\frac{1}{\sqrt{l\pi}}, &j=0,\\
\frac{1}{\sqrt{2l\pi}},&j =2n_2, \\
0,&otherwise,
\end{array}\right.\\
\end{array}
\end{equation*}
and
\[\begin{array}{l}
D_yf_2^1(z,0,0)(U_2^2(z,0))\\
=\Psi(0)\left( \begin{array}{c}
\sum\limits_{j\geq 0}[b_{jn_1n_1} (S_{yz_1}(h_j)z_1+S_{yz_2}( h_j)z_2)+b_{jn_2n_1}(S_{yz_3}( h_j)z_3+S_{yz_4}( h_j)z_4)]  \\
\sum\limits_{j\geq 0} [ b_{jn_1n_2}(S_{yz_1}( h_j)z_1+ S_{yz_2}( h_j)z_2)+b_{jn_2n_2}(S_{yz_3}( h_j)z_3+S_{yz_4}( h_j)z_4)]  \\
\end{array}\right), \\
=\Psi(0)\left( \begin{array}{c}
\frac{1}{\sqrt{l\pi}}( S_{yz_1}(h_0)z_1+S_{yz_2}( h_0)z_2)+\frac{1}{\sqrt{l\pi}}(S_{yz_3}( h_{n_2})z_3+S_{yz_4}( h_{n_2})z_4)  \\
 \frac{1}{\sqrt{l\pi}}(S_{yz_1}( h_{n_2})z_1+ S_{yz_2}( h_{n_2})z_2)+\frac{1}{\sqrt{l\pi}}(S_{yz_3}( h_0)z_3+S_{yz_4}( h_0)z_4)\\+\frac{1}{\sqrt{2l\pi}}(S_{yz_3}( h_{2n_2})z_3+S_{yz_4}( h_{2n_2})z_4)  \\
\end{array}\right). \\
\end{array}\]
Thus, we obtain
\begin{equation*}
\frac{1}{3!}{\rm Proj} _{{\rm Ker}(M_3^1)^c}(D_yf_2^1(z,0,0)(U_2^2(z,0)))\\
=\left( \begin{array}{c}
 E_{2100}z_1^2z_2+E_{1011}z_1z_3z_4 \\
 \overline{E_{2100}}z_1z_2^2+\overline{E_{1011}}z_2z_3z_4 \\
 E_{0021}z_3^2z_4+E_{1110}z_1z_2z_3  \\
 \overline{ E_{0021}}z_3z_4^2+\overline{E_{1110}}z_1z_2z_4  \\
  \end{array}\right) ,\\
\end{equation*}
where
\begin{equation}
\begin{array}{l}
E_{2100}=\frac{1}{6\sqrt{l\pi}}\psi_1(0)\left[ S_{yz_1}(h_{01100})+S_{yz_2}( h_{02000})\right], \\
E_{1011}=\frac{1}{6\sqrt{l\pi}}\psi_1(0)\left[ S_{yz_1}h_{00011})+S_{yz_3}( h_{n_21001})+S_{yz_4}( h_{n_21010})\right], \\
E_{0021}=\frac{1}{6}\psi_3(0)\left[\frac{1}{\sqrt{l\pi}}( S_{yz_3}( h_{00011})+S_{yz_4}( h_{00020}))+\frac{1}{\sqrt{2l\pi}}( S_{yz_3}( h_{2n_20011})+S_{yz_4}( h_{2n_20020}))\right], \\
E_{1110}=\frac{1}{6}\psi_3(0)\left[\frac{1}{\sqrt{l\pi}}( S_{yz_1}( h_{n_20110})+S_{yz_2}( h_{n_21010})+S_{yz_3}( h_{01100}))+\frac{1}{\sqrt{2l\pi}}S_{yz_3}( h_{2n_21100})\right] .
\end{array}
\end{equation}
Clearly, we still need to calculate $h_{01100}$, $h_{02000}$, $h_{00011}$, $h_{n_21001}$, $h_{n_21010}$, $h_{00020}$, $h_{2n_20011}$, $h_{2n_20020}$, $h_{n_20110}$, and $h_{2n_21100}$.  Using the same method used in case 1, we can get the following results
\begin{equation*}
\begin{aligned}
h_{02000}(\theta)&=\frac{1}{-i\omega_1}\phi_1(\theta)f_{2000}^{1(1)}+\frac{1}{-3i\omega_1}\phi_2(\theta)f_{2000}^{1(2)}\\&-\frac{1}{\sqrt{l\pi}}e^{2i\omega_1\theta}[-2i\omega_1+L_0(e^{2i\omega_1\cdot}I_d)]^{-1}F_{2000},
\end{aligned}
\end{equation*}
\begin{equation*}
h_{01100}(\theta)=\frac{1}{i\omega_1}\phi_1(\theta)f_{1100}^{1(1)}+\frac{1}{-i\omega_1}\phi_2(\theta)f_{1100}^{1(2)}-\frac{1}{\sqrt{l\pi}}[L_0(I_d)]^{-1}F_{1100},
\end{equation*}
\begin{equation*}
h_{00011}(\theta)=\frac{1}{i\omega_1}\phi_1(\theta)f_{0011}^{1(1)}+\frac{1}{-i\omega_1}\phi_2(\theta)f_{0011}^{1(2)}- \frac{1}{\sqrt{l\pi}}[L_0(I_d)]^{-1}F_{0011},
\end{equation*}
\begin{equation*}
\begin{aligned}
h_{n_21001}(\theta)=&
\frac{1}{-i\omega_1+2i\omega_2}\phi_3(\theta)f_{1001}^{1(3)}+\frac{1}{-i\omega_1}\phi_4(\theta)f_{1001}^{1(4)}\\&- \frac{1}{\sqrt{l\pi}}e^{(i\omega_1-i\omega_2)\theta}[-(i\omega_1-i\omega_2)-\frac{n_2^2}{l^2}D_0+L_0(e^{(i\omega_1-i\omega_2)\cdot}I_d)]^{-1}F_{1001},
\end{aligned}
\end{equation*}
\begin{equation*}
\begin{aligned}
h_{n_21010}(\theta)&=\frac{1}{-i\omega_1}\phi_3(\theta)f_{1010}^{1(3)}+\frac{1}{-i\omega_1-2i\omega_2}\phi_4(\theta)f_{1010}^{1(4)}\\&- \frac{1}{\sqrt{l\pi}}e^{(i\omega_1+i\omega_2)\theta}[-(i\omega_1+i\omega_2)-\frac{n_2^2}{l^2}D_0+L_0(e^{(i\omega_1+i\omega_2)\cdot}I_d)]^{-1}F_{1010},
\end{aligned}
\end{equation*}
\begin{equation*}
\begin{aligned}
h_{00020}(\theta)&=\frac{1}{i\omega_1-2i\omega_2}\phi_1(\theta)f_{0020}^{1(1)}+\frac{1}{-i\omega_1-2i\omega_2}\phi_2(\theta)f_{0020}^{1(2)}\\
&-\frac{1}{\sqrt{l\pi}}e^{2i\omega_2\theta}[-2i\omega_2+L_0(e^{2i\omega_2\cdot}I_d)]^{-1}F_{0002},
\end{aligned}
\end{equation*}
\begin{equation*}
\begin{aligned}
h_{n_20110}(\theta)&=\frac{1}{i\omega_1}\phi_3(\theta)f_{0110}^{1(3)}+\frac{1}{i\omega_1-2i\omega_2}\phi_4(\theta)f_{0110}^{1(4)}\\&- \frac{1}{\sqrt{l\pi}}e^{(-i\omega_1+i\omega_2)\theta}[-(-i\omega_1+i\omega_2)-\frac{n_2^2}{l^2}D_0+L_0(e^{(-i\omega_1+i\omega_2)\cdot}I_d)]^{-1}F_{0110},
\end{aligned}
\end{equation*}
\begin{equation*}
h_{2n_20011}(\theta)=-\frac{1}{\sqrt{2l\pi}}[-\frac{(2n_2)^2}{l^2}D_0+L_0(I_d)]^{-1}F_{0011},
\end{equation*}
\begin{equation*}
h_{2n_20020}(\theta)=- \frac{1}{\sqrt{2l\pi}}e^{2i\omega_2\theta}[-2i\omega_2-\frac{(2n_2)^2}{l^2}D_0+L_0(e^{2i\omega_2\cdot}I_d)]^{-1}F_{0020},
\end{equation*}
and
\begin{equation*}
h_{2n_21100}(\theta)=0.
\end{equation*}

\textbf{Case 3}
When $n_1\neq 0,n_2\neq 0$, from
\begin{equation*}
\begin{array}{l}
b_{jn_kn_k}=\int_0^{l\pi}\gamma_j(x)\gamma_{n_k}(x)\gamma_{n_k}(x)dx=\left\lbrace \begin{array}{ll}
\frac{1}{\sqrt{l\pi}}, &j=0,\\
\frac{1}{\sqrt{2l\pi}},&j =2n_k ,\\
0,&otherwise,
\end{array} (k=1,2),\right.\\
b_{jn_1n_2}=b_{jn_2n_1}=\int_0^{l\pi}\gamma_j(x)\gamma_{n_1}(x)\gamma_{n_2}(x)dx=\left\lbrace \begin{array}{ll}
\frac{1}{\sqrt{2l\pi}}, &j=n_1+n_2,\\
\frac{1}{\sqrt{2l\pi}},&j =n_2-n_1\neq 0,\\
\frac{1}{\sqrt{l\pi}},&j =n_2-n_1=0,\\
0,&otherwise,
\end{array}\right.\\
\end{array}
\end{equation*}
we have
\[\begin{array}{l}
D_yf_2^1(z,0,0)(U_2^2(z,0))\\
=\Psi(0)\left( \begin{array}{c}
\sum\limits_{j\geq 0}[b_{jn_1n_1} (S_{yz_1}(h_j)z_1+S_{yz_2}( h_j)z_2)+b_{jn_2n_1}(S_{yz_3}( h_j)z_3+S_{yz_4}( h_j)z_4)]  \\
\sum\limits_{j\geq 0} [ b_{jn_1n_2}(S_{yz_1}( h_j)z_1+S_{yz_2}( h_j)z_2)+b_{jn_2n_2}(S_{yz_3}( h_j)z_3+S_{yz_4}( h_j)z_4)]  \\
\end{array}\right) \\
=\Psi(0)\left( \begin{array}{c}
b_{0n_1n_1} (S_{yz_1}(h_0)z_1+S_{yz_2}( h_0)z_2)+b_{2n_1n_1n_1}(S_{yz_1}(h_{2n_1})z_1+S_{yz_2}( h_{2n_1})z_2)\\
+b_{(n_1+n_2)n_2n_1}(S_{yz_3}( h_{n_1+n_2})z_3+S_{yz_4}( h_{n_1+n_2})z_4)\\+b_{n_2-n_1n_2n_1}(S_{yz_3}( h_{n_2-n_1})z_3+S_{yz_4}( h_{n_2-n_1})z_4)\\
b_{(n_1+n_2)n_1n_2} (S_{yz_1}(h_{n_1+n_2})z_1+S_{yz_2}( h_{n_1+n_2})z_2)\\+b_{(n_2-n_1)n_1n_2}(S_{yz_1}(h_{n_2-n_1})z_1+S_{yz_2}( h_{n_2-n_1})z_2)\\+b_{0n_2n_2}(S_{yz_3}( h_0)z_3+S_{yz_4}( h_0)z_4)+b_{2n_2n_2n_2}(S_{yz_3}( h_{2n_2})z_3+S_{yz_4}( h_{2n_2})z_4)\\
\end{array}\right) \\
=\Psi(0)\left( \begin{array}{c}
\frac{1}{\sqrt{l\pi}} (S_{yz_1}(h_0)z_1+S_{yz_2}( h_0)z_2)+\frac{1}{\sqrt{2l\pi}}(S_{yz_1}(h_{2n_1})z_1+S_{yz_2}( h_{2n_1})z_2)\\
+\frac{1}{\sqrt{2l\pi}}(S_{yz_3}( h_{n_1+n_2})z_3+S_{yz_4}( h_{n_1+n_2})z_4)\\+b_{n_2-n_1n_2n_1}(S_{yz_3}( h_{n_2-n_1})z_3+S_{yz_4}( h_{n_2-n_1})z_4)\\
\frac{1}{\sqrt{2l\pi}} (S_{yz_1}(h_{n_1+n_2})z_1+S_{yz_2}( h_{n_1+n_2})z_2)\\+b_{(n_2-n_1)n_1n_2}(S_{yz_1}(h_{n_2-n_1})z_1+S_{yz_2}( h_{n_2-n_1})z_2)\\+\frac{1}{\sqrt{l\pi}}(S_{yz_3}( h_0)z_3+S_{yz_4}( h_0)z_4)+\frac{1}{\sqrt{2l\pi}}(S_{yz_3}( h_{2n_2})z_3+S_{yz_4}( h_{2n_2})z_4)\\
\end{array}\right) ,
\end{array}\]
thus,
\begin{equation}
\label{3!Dyf21z00U22(2)}
{\rm Proj} _{{\rm Im}(M_3^1)^c}D_yf_2^1(z,0,0)(U_2^2(z,0))\\
=\left( \begin{array}{c}
 E_{2100}z_1^2z_2+E_{1011}z_1z_3z_4 \\
 \overline{E_{2100}}z_1z_2^2+\overline{E_{1011}}z_2z_3z_4 \\
 E_{0021}z_3^2z_4+E_{1110}z_1z_2z_3  \\
 \overline{ E_{0021}}z_3z_4^2+\overline{E_{1110}}z_1z_2z_4  \\
  \end{array}\right), \\
\end{equation}
where
\begin{equation}
\begin{aligned}
E_{2100}=&\frac{1}{6}\psi_1(0)[ \frac{1}{\sqrt{l\pi}}(S_{yz_1}(h_{01100})+S_{yz_2}( h_{02000}))+\frac{1}{\sqrt{2l\pi}}(S_{yz_1}(h_{2n_11100})+S_{yz_2}( h_{2n_12000}))], \\
E_{1011}=&\frac{1}{6}\psi_1(0)[\frac{1}{\sqrt{l\pi}}(S_{yz_1}(h_{00011})+\frac{1}{\sqrt{2l\pi}}(S_{yz_1}(h_{2n_10011})+S_{yz_3}( h_{(n_1+n_2)1001})\\&+S_{yz_4}( h_{(n_1+n_2)1010})) +b_{(n_2-n_1)n_2n_1}(S_{yz_3}( h_{(n_2-n_1)1001})+S_{yz_4}( h_{(n_2-n_1)1010}))], \\
E_{0021}=&\frac{1}{6}\psi_3(0)[ \frac{1}{\sqrt{l\pi}}(S_{yz_3}(h_{00011})+S_{yz_4}( h_{00020}))+\frac{1}{\sqrt{2l\pi}}(S_{yz_3}(h_{2n_20011})+S_{yz_4}( h_{2n_20020}))] ,\\
E_{1110}=&\frac{1}{6}\psi_3(0)[\frac{1}{\sqrt{l\pi}}(S_{yz_3}(h_{01100})+\frac{1}{\sqrt{2l\pi}}(S_{yz_3}(h_{2n_21100})+S_{yz_1}( h_{(n_1+n_2)0110})\\&+S_{yz_2}( h_{(n_1+n_2)1010})) +b_{(n_2-n_1)n_1n_2}(S_{yz_1}( h_{(n_2-n_1)0110})+S_{yz_2}( h_{(n_2-n_1)1010}))]. \\
\end{aligned}
\end{equation}
Clearly, we still need to calculate $h_{01100}$, $h_{02000}$, $h_{2n_11100}$, $h_{2n_12000}$, $h_{00011}$, $h_{2n_10011},$  $h_{(n_1+n_2)1001}$, $h_{(n_1+n_2)1010} $, $h_{(n_2-n_1)1001}$, $h_{(n_2-n_1)1010}$, $h_{00020},$ $h_{2n_20011}$, $h_{2n_20020}$, $h_{2n_21100}$, $h_{(n_1+n_2)0110}$, $h_{(n_1+n_2)1010}$, $h_{(n_2-n_1)0110}$, and $h_{(n_2-n_1)1010}.$
 In fact, we have
\begin{equation*}
h_{02000}(\theta)=- \frac{1}{\sqrt{l\pi}}e^{2i\omega_1\theta}[-2i\omega_1+L_0(e^{2i\omega_1\cdot}I_d)]^{-1}F_{2000},
\end{equation*}
\begin{equation*}
h_{01100}(\theta)=-\frac{1}{\sqrt{l\pi}}[L_0(I_d)]^{-1}F_{1100},
\end{equation*}
\begin{equation*}
h_{2n_12000}(\theta)=\left\lbrace
\begin{array}{ll}
-\frac{1}{\sqrt{2l\pi}}e^{2i\omega_1\theta}[-2i\omega_1-\frac{(2n_1)^2}{l^2}D_0+L_0(e^{2i\omega_1\cdot}I_d)]^{-1}F_{2000}, &n_2\neq 2n_1,\\
\frac{1}{-2i\omega_1+i\omega_2}\phi_3(\theta)f_{2000}^{1(3)}+\frac{1}{-2i\omega_1-i\omega_2}\phi_4(\theta)f_{2000}^{1(4)}\\-\frac{1}{\sqrt{2l\pi}}[-2i\omega_1-\frac{(2n_1)^2}{l^2}D_0+L_0(e^{2i\omega_1\cdot}I_d)]^{-1}F_{2000},&n_2= 2n_1,\\
\end{array}\right.
\end{equation*}
\begin{equation*}
h_{2n_11100}(\theta)=\left\lbrace \begin{array}{ll}
-\frac{1}{\sqrt{2l\pi}}[-D_0\frac{(2n_1)^2}{l^2}+L_0(I_d)]^{-1}F_{1100},& n_2\neq 2n_1,\\
-\frac{1}{\sqrt{2l\pi}}[-D_0\frac{(2n_1)^2}{l^2}+L_0(I_d)]^{-1}F_{1100}+\frac{1}{i\omega_2}\phi_3(\theta)f_{1100}^{1(3)}\\+\frac{1}{-i\omega_2}\phi_4(\theta)f_{1100}^{1(4)},& n_2=2n_1,\\
\end{array}\right.
\end{equation*}
\begin{equation*}
h_{00011}(\theta)=-\frac{1}{\sqrt{l\pi}}[L_0(I_d)]^{-1}F_{0011},
\end{equation*}
\begin{equation*}
h_{2n_10011}(\theta)=\left\lbrace \begin{array}{ll}
-\frac{1}{\sqrt{2l\pi}}[-D_0\frac{(2n_1)^2}{l^2}+L_0(I_d)]^{-1}F_{0011},& n_2=n_1,\\
\frac{1}{i\omega_2}\phi_3(\theta)f_{0011}^{1,3}+\frac{1}{-i\omega_2}\phi_4(\theta)f_{0011}^{1,4},&n_2=2n_1,\\
0,&otherwise,
\end{array}\right.
\end{equation*}
\begin{equation*}
\begin{aligned}
h_{n_1+n_21001}(\theta)&=- \frac{1}{\sqrt{2l\pi}}e^{(i\omega_1-i\omega_2)\theta}[-(i\omega_1-i\omega_2)\\&-\frac{(n_1+n_2)^2}{l^2}D_0+L_0(e^{(i\omega_1-i\omega_2)\cdot}I_d)]^{-1}F_{1001},
\end{aligned}
\end{equation*}
\begin{equation*}
\begin{aligned}
h_{n_1+n_21010}(\theta)&=- \frac{1}{\sqrt{2l\pi}}e^{(i\omega_1+i\omega_2)\theta}[-(i\omega_1+i\omega_2)\\&-\frac{(n_1+n_2)^2}{l^2}D_0+L_0(e^{(i\omega_1+i\omega_2)\cdot}I_d)]^{-1}F_{1010},
\end{aligned}
\end{equation*}
\begin{equation*}
h_{n_2-n_11001}(\theta)=\left\lbrace \begin{array}{ll}
- \frac{1}{\sqrt{l\pi}}e^{(i\omega_1-i\omega_2)\theta}[-(i\omega_1-i\omega_2)+L_0(e^{(i\omega_1-i\omega_2)\cdot}I_d)]^{-1}F_{1001},& n_2=n_1,\\\frac{1}{i\omega_2}\phi_1(\theta)f_{1001}^{1(1)}+\frac{1}{-2i\omega_1+i\omega_2}\phi_2(\theta)f_{1001}^{1(2)}
\\- \frac{1}{\sqrt{2l\pi}}e^{(i\omega_1-i\omega_2)\theta}[-(i\omega_1-i\omega_2)-\frac{(n_2-n_1)^2}{l^2}D_0\\~~~~~~~~~~~~~~~~~~~~~~~~~~~~+L_0(e^{(i\omega_1-i\omega_2)\cdot}I_d)]^{-1}F_{1001},& n_2= 2n_1,\\
- \frac{1}{\sqrt{2l\pi}}e^{(i\omega_1-i\omega_2)\theta}[-(i\omega_1-i\omega_2)-\frac{(n_2-n_1)^2}{l^2}D_0\\~~~~~~~~~~~~~~~~~~~~~~~~~~~~+L_0(e^{(i\omega_1-i\omega_2)\cdot}I_d)]^{-1}F_{1001},& otherwise,\\
\end{array}\right.
\end{equation*}
\begin{equation*}
h_{n_2-n_11010}(\theta)=\left\lbrace \begin{array}{ll}
- \frac{1}{\sqrt{l\pi}}e^{(i\omega_1+i\omega_2)\theta}[-(i\omega_1+i\omega_2)+L_0(e^{(i\omega_1+i\omega_2)\cdot}I_d)]^{-1}F_{1010},& n_2=n_1,\\
\frac{1}{-i\omega_2}\phi_1(\theta)f_{1010}^{1(1)}+\frac{1}{-2i\omega_1-i\omega_2}\phi_2(\theta)f_{1010}^{1(2)}
\\-\frac{1}{\sqrt{2l\pi}}e^{(i\omega_1+i\omega_2)\theta}[-(i\omega_1+i\omega_2)-\frac{(n_2-n_1)^2}{l^2}D_0\\~~~~~~~~~~~~~~~~~~~~~~~~~~~~+L_0(e^{(i\omega_1+i\omega_2)\cdot}I_d)]^{-1}F_{1010},& n_2= 2n_1,\\
- \frac{1}{\sqrt{2l\pi}}e^{(i\omega_1+i\omega_2)\theta}[-(i\omega_1+i\omega_2)-\frac{(n_2-n_1)^2}{l^2}D_0\\~~~~~~~~~~~~~~~~~~~~~~~~~~~~~+L_0(e^{(i\omega_1+i\omega_2)\cdot}I_d)]^{-1}F_{1010},&otherwise,\\
\end{array}\right.
\end{equation*}
\begin{equation*}
h_{00020}(\theta)=- \frac{1}{\sqrt{l\pi}}e^{2i\omega_2\theta}[-2i\omega_2+L_0(e^{2i\omega_2\cdot}I_d)]^{-1}F_{0020},
\end{equation*}
\begin{equation*}
h_{2n_20011}(\theta)=-\frac{1}{\sqrt{2l\pi}}[-\frac{(2n_2)^2}{l^2}D_0+L_0(I_d)]^{-1}F_{0011},
\end{equation*}
\begin{equation*}
h_{2n_20020}(\theta)=- \frac{1}{\sqrt{2l\pi}}e^{2i\omega_2\theta}[-2i\omega_2-\frac{(2n_2)^2}{l^2}D_0+L_0(e^{2i\omega_2\cdot}I_d)]^{-1}F_{0020},
\end{equation*}
\begin{equation*}
h_{2n_21100}(\theta)=\left\lbrace \begin{array}{cc}
-\frac{1}{\sqrt{2l\pi}}[-\frac{(2n_2)^2}{l^2}D_0+L_0(I_d)]^{-1}F_{1100}, & n_2=n_1,\\
0,&n_2\neq n_1,
\end{array}\right.
\end{equation*}
\begin{equation*}
\begin{aligned}
h_{n_1+n_20110}(\theta)&=- \frac{1}{\sqrt{2l\pi}}e^{(-i\omega_1+i\omega_2)\theta}[-(-i\omega_1+i\omega_2)\\&-\frac{(n_1+n_2)^2}{l^2}D_0+L_0(e^{(-i\omega_1+i\omega_2)\cdot}I_d)]^{-1}F_{0110},
\end{aligned}
\end{equation*}
\begin{equation*}
h_{n_2-n_10110}(\theta)=\left\lbrace \begin{array}{ll}- \frac{1}{\sqrt{l\pi}}e^{(-i\omega_1+i\omega_2)\theta}[-(-i\omega_1+i\omega_2)\\~~~~~~~~~~~~~~~~~~~~~~~~~~+L_0(e^{(-i\omega_1+i\omega_2)\cdot}I_d)]^{-1}F_{0110},& n_2=n_1,\\
\frac{1}{2i\omega_1-i\omega_2}\phi_1(\theta)f_{0110}^{1(1)}+\frac{1}{-i\omega_2}\phi_2(\theta)f_{0110}^{1(2)}
\\- \frac{1}{\sqrt{2l\pi}}e^{(-i\omega_1+i\omega_2)\theta}[-(-i\omega_1+i\omega_2)-\frac{(n_2-n_1)^2}{l^2}D_0\\~~~~~~~~~~~~~~~~~~~~~~~~~~+L_0(e^{(-i\omega_1+i\omega_2)\cdot}I_d)]^{-1}F_{0110},& n_2= 2n_1,\\- \frac{1}{\sqrt{2l\pi}}e^{(-i\omega_1+i\omega_2)\theta}[-(-i\omega_1+i\omega_2)-\frac{(n_2-n_1)^2}{l^2}D_0\\~~~~~~~~~~~~~~~~~~~~~~~~~~+L_0(e^{(-i\omega_1+i\omega_2)\cdot}I_d)]^{-1}F_{0110},& otherwise.\\
\end{array}\right.
\end{equation*}

Hence, by (\ref{f31bar}),  (\ref{3!Dzf21z00U21}), (\ref{3!Dyf21z00U22(1)}), and (\ref{3!Dyf21z00U22(2)}), we have
\begin{equation}\label{g31z00}
\frac{1}{3!}g_3^1(z,0,0)=\frac{1}{3!}{\rm Proj}_{{\rm Im }(M_3^1)^c}\overline{f}_3^1(z,0,0)=\left( \begin{array}{c}
 B_{2100}z_1^2z_2+B_{1011}z_1z_3z_4 \\
 \overline{B_{2100}}z_1z_2^2+\overline{B_{1011}}z_2z_3z_4 \\
 B_{0021}z_3^2z_4+B_{1110}z_1z_2z_3  \\
 \overline{ B_{0021}}z_3z_4^2+\overline{B_{1110}}z_1z_2z_4  \\
  \end{array}\right), \\
\end{equation}
with
\begin{equation}\label{B2100}
\begin{aligned}
B_{2100}=C_{2100}+\frac{3}{2}(D_{2100}+E_{2100}),~~B_{1011}=C_{1011}+\frac{3}{2}(D_{1011}+E_{1011}),\\
B_{0021}=C_{0021}+\frac{3}{2}(D_{0021}+E_{0021}),~~B_{1110}=C_{1110}+\frac{3}{2}(D_{1110}+E_{1110}).\\
\end{aligned}
\end{equation}
Therefore, by (\ref{2!g21zalpha}) and (\ref{g31z00}), the normal form truncated to the third order for double Hopf bifurcation reads as
\begin{equation}\label{normalform}
\begin{aligned}
&\dot{z}_1=&i\omega_1z_1&+B_{11}\alpha_1z_1+B_{21}\alpha_2z_1+ B_{2100}z_1^2z_2+B_{1011}z_1z_3z_4, \\
&\dot{z}_2=&-i\omega_1z_2&+\overline{B_{11}}\alpha_1z_2+\overline{B_{21}}\alpha_2z_2+ \overline{B_{2100}}z_1z_2^2+\overline{B_{1011}}z_2z_3z_4,\\
&\dot{z}_3=&i\omega_2z_3&+B_{13}\alpha_1z_3+B_{23}\alpha_2z_3+ B_{0021}z_3^2z_4+B_{1110}z_1z_2z_3 ,\\
&\dot{z}_4=&-i\omega_2z_4&+\overline{B_{13}}\alpha_1z_4+\overline{B_{23}}\alpha_2z_4+ \overline{B_{0021}}z_3z_4^2+\overline{B_{1110}}z_1z_2z_4 .
\end{aligned}
\end{equation}

\begin{remark}\label{calculationsteps}
To sum up, the whole calculation process above can be accomplished by following three  steps.

(1) Obtain the double Hopf bifurcation point by analyzing the associate characteristic equation, and find $n_1$, $n_2$.
Rewrite the original system into the form as (\ref{generalalpha}), and calculate $D_0$, $L_0$, $D_1^{(1,0)}$, $D_1^{(0,1)}$, $L_1^{(1,0)}$, and $L_1^{(0,1)}$.  Calculate the eigenfunctions $\phi_i$ and $\psi_i$ $(i=1,2,3,4)$.

(2) Calculate $B_{11}$, $B_{21}$, $B_{13}$, and $B_{23}$ from (\ref{B11}).

(3) Calculate $F_{mnij}$ $(m+n+i+j=3)$, and we can get  $C_{2100}$, $C_{1011}$, $C_{0021}$, and $C_{1110}$ from (\ref{C2100});  Calculate $F_{mnij}$ $(m+n+i+j=2)$, and we obtain  $D_{2100}$, $D_{1011}$, $D_{0021}$, and $D_{1110}$ from (\ref{D2100})-(\ref{D1110});  Calculate $f_{mnij}^{1(k)}$ $(m+n+i+j=2, k=1,2,3,4)$ by (\ref{f200011})-(\ref{f000211}), establish the linear operators $S_{yz_i} (i=1,2,3,4)$,  we can get  $E_{2100}$,  $E_{1011}$, $E_{0021}$, and $E_{1110}$ by the formulas in step3 in three different cases.  Finally, by (\ref{B2100}), we can get $B_{2100}$, $B_{1011}$, $B_{0021}$, and $B_{1110}$.
\end{remark}

Make double polar coordinates transformation by
\[\begin{array}{l}
z_1=r_1\cos \theta_1+ir_1\sin\theta_1,\\z_2=r_1\cos \theta_1-ir_1\sin\theta_1,\\z_3=r_2\cos \theta_2+ir_2\sin\theta_2,\\z_4=r_2\cos \theta_2-ir_2\sin\theta_2,
\end{array}\]
where $r_1,r_2>0$. Define $\epsilon_1={\rm Sign}({\rm Re}B_{2100})$,  $\epsilon_2={\rm Sign}({\rm Re}B_{0021})$, carry out the rescaling $\widehat{r}_1=r_1\sqrt{|{\rm}B_{2100}|}$, $\widehat{r}_2=r_2\sqrt{|{\rm}B_{0021}|}$, $\widehat{t}=t\epsilon_1$, and drop the hats, then system (\ref{normalform}) becomes
 \begin{equation}\label{normalformcylin}
 \begin{aligned}
 &\dot{r}_1=r_1(c_1+r_1^2+b_0r_2^2),\\
 &\dot{r}_2=r_2(c_2+c_0r_1^2+d_0r_2^2),
 \end{aligned}
 \end{equation}
 where
 \[\begin{aligned}
& c_1=\epsilon_1({\rm Re}B_{11}\alpha_1+{\rm Re}B_{21}\alpha_2)=\epsilon_1({\rm Re}B_{11}(\mu_1-\mu_{0,1})+{\rm Re}B_{21}(\mu_2-\mu_{0,2})),\\
& c_2=\epsilon_1({\rm Re}B_{13}\alpha_1+{\rm Re}B_{23}\alpha_2)=\epsilon_1({\rm Re}B_{13}(\mu_1-\mu_{0,1})+{\rm Re}B_{23}(\mu_2-\mu_{0,2}))),  \\&b_0=\frac{\epsilon_1\epsilon_2{\rm Re}B_{1011}}{{\rm Re}B_{0021}}, c_0=\frac{{\rm Re}B_{1110}}{{\rm Re}B_{2100}},d_0=\epsilon_1\epsilon_2.\\
 \end{aligned}\]
By chapter 7.5 in \cite{Guckenheimer}, Eq. (\ref{normalformcylin}) has twelve distinct kinds of unfoldings (see Table 1).
 \begin{table}[tbp]
 \label{twelve}
 \caption{The twelve unfoldings of system (\ref{normalformcylin}).}
 \centering
 \begin{tabular}{lcccccccccccc}
 \hline
 Case & {I}a & {I}b  & {II} & {III}& {IV}a & {IV}b  & {V}& {VI}a& {VI}b& {VII}a& {VII}b& {VIII}     \\ \hline  
 $d_0$ &  +1 &+1& +1&+1&+1&+1&--1&--1&--1&--1&--1&--1\\         
 $b_0 $ &  + &+& +&--&--&--&+&+&+&--&--&--\\        
 $c_0$ &  + &+& --&+&--&--&+&--&--&+&+&--\\
 $d_0-b_0c_0$ &  + &--& +&+&+&--&--&+&--&+&--&--\\
 \hline
 \end{tabular}
 \end{table}

\begin{figure}
                                \centering
                    \includegraphics[width=0.7\textwidth]{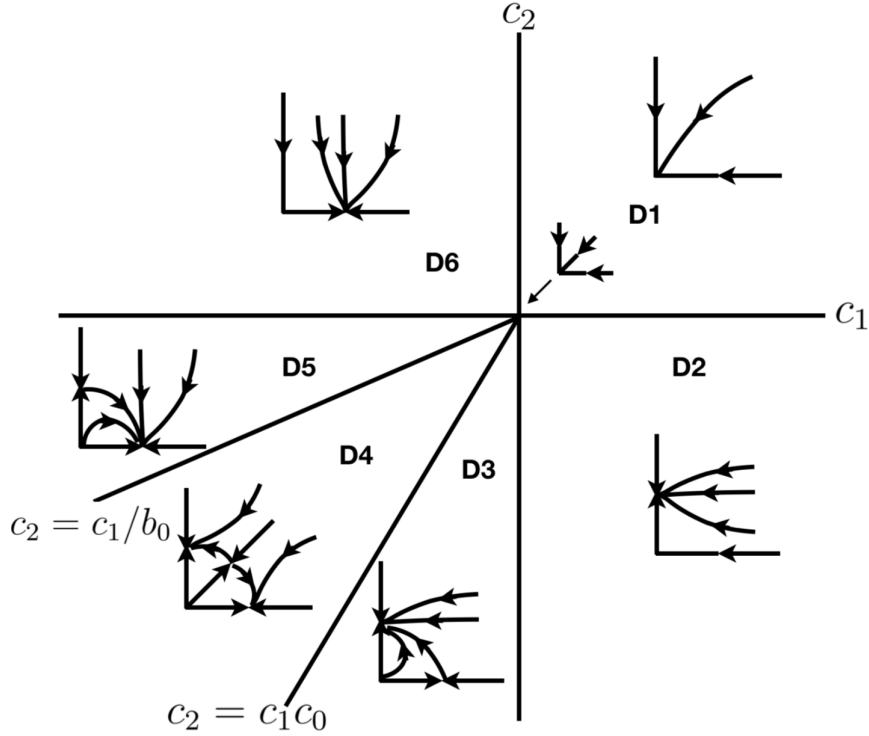}
                      \caption{Phase portraits for  the unfoldings of case Ib with $\epsilon_1=-1$.}
          \label{fig:Ib}
                        \end{figure}

\begin{remark} In  section  \ref{section epidemic},   case Ib appears, thus we draw bifurcation set and phase portraits for the unfolding of case  Ib in Figure \ref{fig:Ib}  \cite{Guckenheimer}.
For case Ib, near the double Hopf bifurcation point, the $(\alpha_1, \alpha_2)$ plane is divided into six regions D1-D6.  In region D1,  the  equilibrium is a sink; when the parameters vary and enter the region D2 (or D6), the stable equilibrium bifurcates into a stable periodic solution via supercritical Hopf bifurcations. For parameters in D3 (or D5), periodic solutions appear via secondary supercritical Hopf bifurcations, but they are saddle type and only stable on the center manifold; when the parameters cross the Hopf bifurcation curve   $c_2=c_1c_0$  (or  $c_2=c_1/b_0$) into the region D4,   there are two stable periodic solutions coexisting in D4.
\end{remark}

\begin{figure}
                                \centering
                    \includegraphics[width=0.7\textwidth]{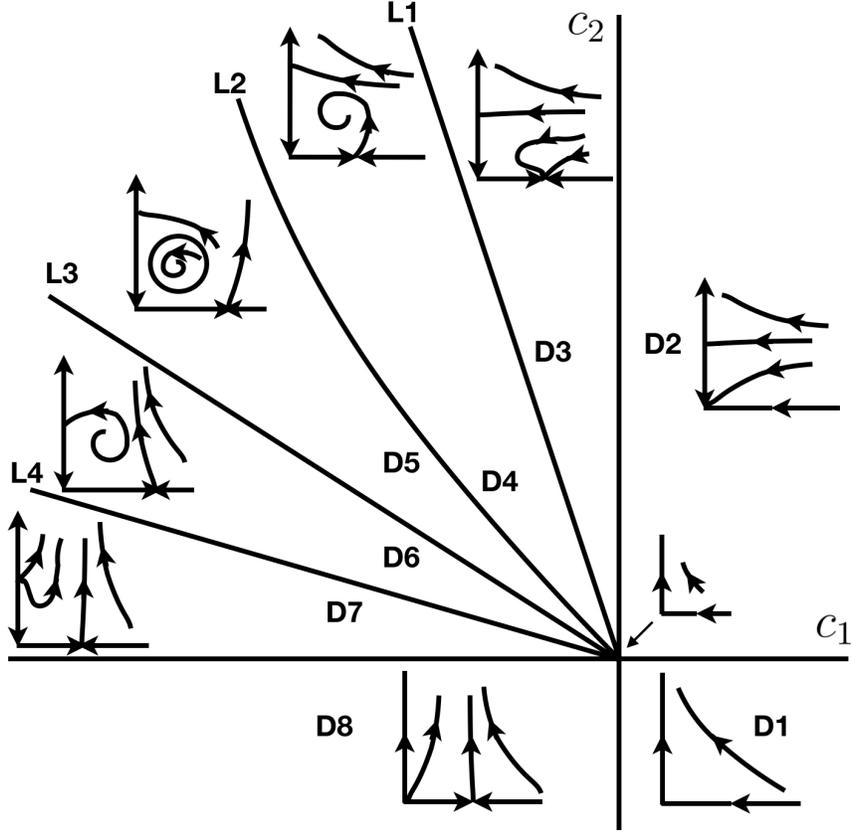}
                      \caption{Phase portraits for  the unfoldings of case VIa with $\epsilon_1=-1$, where $L_1:c_2=c_0c_1$, $L_2:c_2=\frac{c_0-1}{b_0+1}c_1+o(c_1)$, $L_3:c_2=\frac{c_0-1}{b_0+1}c_1$, $L_4:c_2=-\frac{1}{b_0}c_1$.}
          \label{fig:VIa}
                        \end{figure}

\begin{remark} In  section \ref{section predatorprey}  case VIa arises, thus we draw bifurcation set and phase portraits for the unfolding of case  VIa in Figure \ref{fig:VIa}  \cite{Guckenheimer}.
For case VIa, near the double Hopf bifurcation point, the $(\alpha_1, \alpha_2)$ plane is divided into eight regions D1-D8.  In region D2,  the  positive equilibrium is a sink; In region D3, there is a stable periodic solution. In D4, there is a quasi-periodic solution on the two-dimensional torus; In region D5,  there is a quasi-periodic solution on the
three-dimensional torus. When the parameters vary and enter D6, three-dimensional torus vanish.  Generally, a vanishing 3-torus might accompany the phenomenon of chaos \cite{P. Battelino,J.P. Eckmann,D. Ruelle}, so near the double Hopf bifurcation point, strange attractors may exist.
\end{remark}

 \section{Application to a Diffusive Epidemic Model with Delay and Stage Structure}
      \label{section epidemic}

 In this section, a  diffusive epidemic model with delay and stage structure is considered.
 Taking the time delay $\omega$ and a diffusive coefficient as bifurcation parameters, we show that double Hopf bifurcation can undergo with  two wave numbers in different cases, such as $n_1=0,n_2\neq 0$, and  $n_1\neq 0,n_2\neq 0$.  Following the steps of  Remark \ref{calculationsteps}, the normal form can be calculated, and the unfolding system can be got. Simulations illustrate that the spatio-temporal dynamics turn out to be very complicated near the double Hopf bifurcation point. In some regions, there are two stable nonhomogeneous periodic solutions or a homogeneous and a nonhomogeneous periodic solution coexisting.
\subsection{Model and the Existence of Double Hopf Bifurcation}
The  stage-structured epidemic model with the maturation delay and freely-moving delay is given by
\begin{equation}
\label{diffusion model}
\left\{
\begin{array}{l}
 \dfrac{\partial S(x,t)} {\partial t}= d_1\Delta S(x,t)+\alpha y(x,t)-dS(x,t)-\alpha e^{-d\tau}y(x,t-\tau)\\
~~~~~~~~~~~~~~-\mu S(x,t-\omega)I(x,t)+\gamma I(x,t), \\
\dfrac{\partial I(x,t)}{\partial t }= d_2\Delta I(x,t)+\mu S(x,t-\omega)I(x,t)-dI(x,t)-\gamma I(x,t),\\
\dfrac{\partial y(x,t)}{\partial t}= d_3\Delta y(x,t)+\alpha e^{-d\tau}y(x,t-\tau)-\beta y^2(x,t),~~~~~~~~~~~~~~~~~~~ x\in (0,l\pi),\\
\dfrac{\partial S(x,t)} {\partial x}= 0,~~\dfrac{\partial I(x,t)} {\partial x}=0,~~\dfrac{\partial y(x,t)} {\partial x}= 0, ~at~ x=0~and~l\pi,\\
\end{array}
  \right.
\end{equation}
where $S(x,t)$ and $I(x,t)$ represent the population of susceptible and infected immature individuals, and $y(x,t)$ represents the population of mature individuals.   $\tau$ is the time taken for the immature individuals to maturity,  $\omega$ is the time taken for the immature individuals from birth to moving freely.   $\alpha$ is the natural birth rate, $d$ is the natural death rate of the immature stage, $\beta$ is the death rate of the mature individuals  of logistic nature, $\mu$ is the disease transmission rate,  and $\gamma$ is the recovery rate.

 Denote the basic reproduction ratio by
  $R_0=\frac{\mu\alpha^2e^{-d\tau}(1-e^{-d\tau})}{d\beta(d+\gamma)}$.  We can easily verify that if $R_0>1$, system (\ref{diffusion model}) has a positive constant equilibrium  $E^*~(S^*,I^*,y^*)=(   \frac{d+\gamma}{\mu},\frac{(d+\gamma)}{\mu}(R_0-1), \frac{\alpha e^{-d\tau}}{\beta})$ \cite{duguo}.

Linearizing system (\ref{diffusion model}) at the positive equilibrium $E^* (S^*,I^*,y^*)$, we have
 \begin{equation*}
        \label{linear e2}
      \frac{\partial }{\partial t} \left( \begin{array}{l}
        S(x,t) \\
          I(x,t) \\
          y(x,t) \\
        \end{array}\right)
        =
       (D\Delta      +
       A) \left(  \begin{array}{l}
              S(x,t)\\
               I(x,t)\\
              y(x,t)\\
             \end{array}
             \right)
       + G_1 \left(  \begin{array}{l}
                   S(x,t-\omega)\\
                    I(x,t-\omega)\\
                   y(x,t-\omega)\\
                  \end{array}
                  \right)
       +G_2\left(  \begin{array}{l}
               S(x,t-\tau)\\
                 I(x,t-\tau)\\
            y(x,t-\tau)\\
        \end{array}
        \right),
             \end{equation*}
   where
   $D={\rm diag}\{d_1,d_2,d_3\}$,
     \begin{equation*}
     A=\left( \begin{array}{ccc}
      -d& -\mu S^*+\gamma& \alpha\\
      0& \mu S^*-d-\gamma& 0\\
      0&  0& -2\beta y^*\\
      \end{array}\right),~
      G_1=\left( \begin{array}{ccc}
       -\mu I^*& 0& 0\\
       \mu I^*& 0 & 0\\
           0& 0& 0\\
       \end{array}\right),~
       G_2=\left( \begin{array}{ccc}
        0&  0& -\alpha e^{-d\tau}\\
        0& 0 & 0\\
        0&  0& \alpha e^{-d\tau}\\
         \end{array}\right).
     \end{equation*}
   The corresponding characteristic equation is
    \begin{equation*}
    {\rm det}(\lambda I_3-M_n-A-G_1e^{-\lambda\omega}-G_2e^{-\lambda\tau})=0,
    \end{equation*}
    where $I_3$ is the $3\times 3$ identity matrix and $M_n=-\frac{n^2}{l^2} D$, $n\in \mathbb{N}_0$.  That is, each characteristic value $\lambda$ is a root of an equation
          \begin{equation}
      \label{charactershuang}
       (\lambda-\alpha e^{-d\tau}e^{-\lambda\tau}+2\beta y^*+d_3\frac{n^2}{l^2})\cdot\Delta_n(\lambda,\tau)=0,
      \end{equation}
     where
     \begin{equation}\label{charactershuangdisan}\begin{aligned}
     \Delta_n(\lambda,\tau)=  \lambda^2+(d+d_1\frac{n^2}{l^2}+d_2\frac{n^2}{l^2})\lambda+d_2\frac{n^2}{l^2}(d+d_1\frac{n^2}{l^2})+ e^{-\lambda\omega}(\mu I^*\lambda+\mu I^*d_2\frac{n^2}{l^2}+\mu I^*d),
     \end{aligned}
          \end{equation}
      with $n\in \mathbb{N}_0$.

  We can easily prove that the roots of $\lambda-\alpha e^{-d\tau}e^{-\lambda\tau}+2\beta y^*+d_3\frac{n^2}{l^2}=0$ have negative real part.  To investigate the location of the roots,   it remains to consider the roots of $\Delta_n(\lambda,\tau)=0$.

  When $\omega=0$,  Eq. (\ref{charactershuangdisan}) becomes
 \begin{equation}
 \label{characterdan}
  \lambda^2+T_n\lambda+J_n=0,~~~~~~n\in \mathbb{N}_0.
 \end{equation}
Since $T_n=d+d_1\frac{n^2}{l^2}+d_2\frac{n^2}{l^2}+\mu I^*>0, ~J_n=d_2\frac{n^2}{l^2}(d+d_1\frac{n^2}{l^2})+\mu I^*d_2\frac{n^2}{l^2}+\mu I^*d>0$,
 we know that all roots of Eq. (\ref{characterdan}) have negative real part, and thus so do  the roots of Eq. (\ref{charactershuang}) for  $\omega=0$ when $R_0>1$.

By substituting   $iz~(z>0)$ into  Eq. (\ref{charactershuangdisan}), and separating the real and imaginary parts, we have
      \begin{eqnarray}
   \label{shixu}
   \left\{
   \begin{array}{l}
   -z^2+d_2\frac{n^2}{l^2}(d+d_1\frac{n^2}{l^2})= -\mu I^*z\sin \omega z-(\mu I^*d_2\frac{n^2}{l^2}+\mu I^*d)\cos \omega z ,\\
   (d+d_1\frac{n^2}{l^2}+d_2\frac{n^2}{l^2})z=-\mu I^*z\cos \omega z+(\mu I^*d_2\frac{n^2}{l^2}+\mu I^*d)\sin \omega z ,
   \end{array}
   \right.
   \end{eqnarray}
which is solved by
 \begin{equation*}\label{SnCn}
            \begin{array}{l}
    \sin z_n\omega=\dfrac{(d+d_1\frac{n^2}{l^2}+d_2\frac{n^2}{l^2})(\mu I^*d_2\frac{n^2}{l^2}+\mu I^*d)z_n-[d_2\frac{n^2}{l^2}(d+d_1\frac{n^2}{l^2})-z_n^2]\mu I^*z_n}{(\mu I^*z_n)^2+(\mu I^*d_2\frac{n^2}{l^2}+\mu I^*d)^2}\\\stackrel{\vartriangle}{=}S_n(z_n),\\ \cos z_n\omega =-\dfrac{(d+d_1\frac{n^2}{l^2}+d_2\frac{n^2}{l^2})\mu I^*z_n^2+[d_2\frac{n^2}{l^2}(d+d_1\frac{n^2}{l^2})-z_n^2](\mu I^*d_2\frac{n^2}{l^2}+\mu I^*d)}{(\mu I^*d_2\frac{n^2}{l^2}+\mu I^*d)^2+(\mu I^*z_n)^2}\\\stackrel{\vartriangle}{=}C_n(z_n).
             \end{array}
             \end{equation*}
  Then, we obtain
    \begin{equation}
    \label{fz}
   G(z)=z^4+P_nz^2+Q_n=0,
    \end{equation}
      where
           \begin{equation}
           \label{ABC}
           \begin{array}{l}
            P_n=(d_2\frac{n^2}{l^2})^2+(d+d_1\frac{n^2}{l^2}+\mu I^*)(d+d_1\frac{n^2}{l^2}-\mu I^*),\\
            Q_n=J_nK_n=J_n(d_2\frac{n^2}{l^2}(d+d_1\frac{n^2}{l^2})-\mu I^*d_2\frac{n^2}{l^2}-\mu I^*d).
           \end{array}
                      \end{equation}
    Noticing that  $J_n>0$, the sign of $Q_n$ coincides with that of $K_n$.
    Since $K_0=-\mu I^*d<0$, and $K_n$ is a quadratic polynomial with respect to $n^2$,
    we can conclude that there exists $k_1\in\mathbb{N}_0$, such that
    \begin{equation}
    \label{BDdayuxiaoyu}
    \begin{array}{l}
    K_n<0~~~{\rm for} ~0\leq n \leq k_1,\\
            K_n>0~~~{\rm for} ~ n \geq k_1+1,~n\in\mathbb{N}_0.\\
    \end{array}
       \end{equation}

    Denote the positive real root of the equation $K_n=0$ by $n_2$ $(k_1<n_2<k_1+1)$, then we have
  $          K_{n_2}
           =d_1d_2\frac{1}{l^4}n_2^4+(d-\mu I^*)d_2\frac{1}{l^2}n_2^2-\mu I^*d=0.$
    Since $-\mu I^*d<0$, we have $d_1d_2\frac{1}{l^4}n_2^4+(d-\mu I^*)d_2\frac{1}{l^2}n_2^2=( d+d_1\frac{n_2^2}{l^2}-\mu I^*)d_2\frac{1}{l^2}n_2^2>0$.
    It means that
    \begin{equation*}
    d+d_1\frac{n_2^2}{l^2}-\mu I^*>0.
    \end{equation*}
    By (\ref{ABC}), we have
     \begin{equation}
     \label{ABCno}
   P_{n_2}=(d_2\frac{n_2^2}{l^2})^2+(d+d_1\frac{n_2^2}{l^2}+\mu I^*)(d+d_1\frac{n_2^2}{l^2}-\mu I^*)>0.
               \end{equation}
   Noticing $k_1+1>n_2$ and by (\ref{BDdayuxiaoyu}) and   (\ref{ABCno}), we get
   \begin{equation}
   \label{ABCn}
   P_n>0,~~{\rm for} ~n\geq k_1+1, ~k_1\in\mathbb{N}_0.
   \end{equation}

  From (\ref{BDdayuxiaoyu}) and (\ref{ABCn}), we can conclude that for $n\in\mathbb{N}_0$ and $n\leq k_1$, (\ref{fz}) has only one positive real root
    \begin{equation}
    \label{zk}
    z_n= \sqrt{\dfrac{-P_n+\sqrt{P_n^2-4Q_n}}{2}}.   \end{equation}
           For $n\in\mathbb{N}_0$ and $n\geq k_1+1$, (\ref{fz}) has no positive real roots.

In fact,
\begin{equation*}
\label{Sn}
S_n(z_n)=\dfrac{z_n\{(d+d_1\frac{n^2}{l^2}+d_2\frac{n^2}{l^2})(\mu I^*d_2\frac{n^2}{l^2}+\mu I^*d)-[d_2\frac{n^2}{l^2}(d+d_1\frac{n^2}{l^2})]\mu I^*+z_n^2\mu I^*\}}{(\mu I^*z_n)^2+(\mu I^*d_2\frac{n^2}{l^2}+\mu I^*d)^2},
\end{equation*}
where $(d+d_1\frac{n^2}{l^2}+d_2\frac{n^2}{l^2})(\mu I^*d_2\frac{n^2}{l^2}+\mu I^*d)-[d_2\frac{n^2}{l^2}(d+d_1\frac{n^2}{l^2})]\mu I^*
=\mu I^*[d(d+d_1\frac{n^2}{l^2}+d_2\frac{n^2}{l^2})+(d_2\frac{n^2}{l^2})^2]
>0.$
Thus, when  $n\in\{0,1,\cdots,k_1\}$, $S_n\geq 0$, define
\begin{equation}\label{omeganj}
\omega_n^j=
\frac{ \arccos C_n(z_n)+2j\pi}{z_n}.\\
\end{equation}
 
  Differentiating the two sides of Eq. (\ref{charactershuangdisan}) with respective to $\omega$,    Using (\ref{charactershuangdisan}) and (\ref{shixu}), we obtain
  \begin{equation}\label{transver}
  \frac{\mathrm{d}\mathrm{Re}  \lambda (\omega)}{\mathrm{d} \omega}\mid _{\omega=\omega_n^j}=\frac{\sqrt{P_n^2-4Q_n}}{(\mu I^*)^2z_n^2+(\mu I^*d_2\frac{n^2}{l^2}+\mu I^*d)^2}>0.
  \end{equation}
  From (\ref{omeganj}),  we have  the very first critical value as
      \begin{equation*}
      \omega^*=\omega_{n_0}^0=\min_{n\in \{0,1,\cdots,k_1\}}\{\omega_n^0\}, ~~~~z^*=z_{n_0}.
      \end{equation*}

  Due to the general Hopf bifurcation theorem \cite {Faria,JWu}, we have the following theorem.
 \begin{theorem}
    \label{Hopf bifurcation}
   Suppose $R_0>1$. \\
  (1) The equilibrium $E^*$ of system (\ref{diffusion model}) is locally asymptotically stable for $0\leq \omega<  \omega_{n_0}^0$ and is unstable for $ \omega> \omega_{n_0}^0$.\\
  (2) System (\ref{diffusion model}) undergoes  a Hopf bifurcation at the equilibrium $E^*$ when $\omega=\omega_n^j$, for $j\in \mathbb{N}_0$ and $n\in\{0,1,\cdots,k_1\}$.
  \end{theorem}

  Now we need to give a condition under which a double Hopf bifurcation occurs. From (\ref{omeganj}) and (\ref{transver}), we have that  when  $n\in\{0,1,\cdots,k_1\}$, $
  \omega_n^0=
  \frac{ \arccos C_n(z_n)}{z_n}$,   and
  $   \frac{\mathrm{d}\mathrm{Re}  \lambda (\omega)}{\mathrm{d} \omega}\mid _{\omega=\omega_n^0}>0 $, which means that if we vary the coefficients $d_2$ and $\omega$, and fix other coefficients in Eq. (\ref{diffusion model}), there are $k_1+1$ Hopf bifurcation curves on the $d_2-\omega$ plane. On every Hopf bifurcation curve $\omega_n^0$, the characteristic equation (\ref{charactershuang})  always has one pair of eigenvalues $\pm i z_n $, which  crosses  transversely the imaginary axis when the parameters cross the Hopf bifurcation curve, and  the rest eigenvalues  have non-zero real part. If we can   find the intersection of  two certain Hopf bifurcation curves $\omega_{n_1}^0$ and $\omega_{n_2}^0$, there must exist two pairs of  eigenvalues $\pm i z_{n_1} $ and  $\pm i z_{n_2} $ at the intersection point, and all the other eigenvalues have non-zero real part.  Thus, a double Hopf singularity  can be found by searching for the intersection of the Hopf bifurcation curves, which can be done
   by the following process. Firstly, for $n_1,n_2\in\{0,1,\cdots,k_1\}$, we regard $z_{n_1}$ and $z_{n_2}$  as  functions of $d_2$  from Eq. (\ref{zk}).
   Secondly,  the expression of $\omega_{n_1}^0$ and $\omega_{n_2}^0$ are obtained from (\ref{omeganj}). Finally,  from  $\omega_{n_1}^0=\omega_{n_2}^0$, we can solve the value of $d_2$, denoted by $d_2^*$, such that $\omega_{n_1}^0=\omega_{n_2}^0$. Thus, we have that when $d_2=d_2^*$, $\omega=\omega_{n_1}^0=\omega_{n_2}^0$, the Hopf bifurcation curves $\omega_{n_1}^0$  and $\omega_{n_2}^0$ intersect. Thus, system (\ref{diffusion model}) undergoes double Hopf bifurcation at the intersection.

    \begin{theorem}
       \label{bifurcation}
  Suppose that $R_0>1$ and there exists $d_2^*$, $n_1$, $n_2$, such that when $d_2=d_2^*$, $\omega_{n_1}^0=\omega_{n_2}^0$. Then system (\ref{diffusion model}) undergoes  a double Hopf bifurcation at $E^*$ when  $d_2=d_2^*$, $\omega=\omega_{n_1}^0=\omega_{n_2}^0\stackrel{\triangle}{=}\omega^*$.
  \end{theorem}

      \subsection{Normal Form of Double Hopf Bifurcation}
     Let $u_1(x,t)=S(x,\omega t)-S^*$,      $u_2(x,t)=I(x,\omega t)-I^*$,     $u_3(x,t)=y(x,\omega t)-y^*$. Thus, these transformations not only  transform the equilibrium $(S^*,I^*,y^*)$ into $(0,0,0)$, but also normalize  the delay  $\omega$ to 1, and transform the other delay  $\tau$ into $\tau/\omega$. Denote  $U(t)= (u_1(x,t),u_2(x,t),$ $u_3(x,t))^T$,  and  $\omega=\omega^*+\alpha_1$, $d_2=d_2^*+\alpha_2$, then  system (\ref{diffusion model}) can be transformed into
     \begin{equation}\label{dUdt1}
     \begin{aligned}
   \dfrac{{\rm d}U(t)}{{\rm d}t}=D(\omega^*+\alpha_1,d_2^*+\alpha_2)\Delta U(t)+L(\omega^*+\alpha_1,d_2^*+\alpha_2)(U^t)\\+F(\omega^*+\alpha_1,d_2^*+\alpha_2,U^t),
     \end{aligned}
          \end{equation}
   where
    \begin{equation*}
      \begin{array}{l}
     D(\omega^*+\alpha_1,d_2^*+\alpha_2)= (\omega^*+\alpha_1)\left( \begin{array}{ccc}
          d_1&0&0\\
          0&d_2^*+\alpha_2&0\\
          0&0&d_3
          \end{array}\right)=D_0+\alpha_1D_1^{(1,0)}+\alpha_2D_1^{(0,1)},\\
    L(\omega^*+\alpha_1,d_2^*+\alpha_2)U^t=(L_0+\alpha_1L_1^{(1,0)}+\alpha_2L_1^{(0,1)})U^t,\\
    F(\omega^*+\alpha_1,d_2^*+\alpha_2,U^t)=(\omega^*+\alpha_1)\left( \begin{array}{c}
    -\mu U_1^t(-1)U_2^t(0)\\
    \mu U_1^t(-1)U_2^t(0)\\
    -\beta U_3^{t2}(0)
    \end{array}\right),
      \end{array}
         \end{equation*}
 with $$D_0=\omega^*\left( \begin{array}{ccc}
                     d_1&0&0\\
                     0&d_2^*&0\\
                     0&0&d_3
                     \end{array}\right), D_1^{(1,0)}=\left( \begin{array}{ccc}
                                d_1&0&0\\
                              0&d_2^*&0\\
                           0&0&d_3
                 \end{array}\right), D_1^{(0,1)}=\omega^*\left( \begin{array}{ccc}
                         0&0&0\\
                     0&1&0\\
                       0&0&0
                  \end{array}\right),$$
 $$L_0U^t=\omega^*\left[AU^t(0)+ G_1U^t(-1)+ G_2U^t(-\tau^*)\right],$$    $$L_1^{(1,0)}U^t=\left[AU^t(0)+ G_1U^t(-1)+ G_2U^t(-\tau^*)\right],$$ $$L_1^{(0,1)}U^t=0,$$ and  $$\tau^*=\tau/\omega.$$

 Eq. (\ref{dUdt1}) can be rewritten as
 \begin{equation}
      \label{dUdt2}
      \dfrac{dU}{dt}= D_0\Delta U(t)+L_0U^t+\widetilde{F}(\alpha,U^t),
      \end{equation}
 where
\begin{equation}
\widetilde{F}(\alpha,U^t)=\alpha_1D_1^{(1,0)} \Delta U+\alpha_2D_1^{(0,1)} \Delta U
+\alpha_1L_1^{(1,0)}U^t+F( \omega^*+\alpha_1,d_2^*+\alpha_2,U^t).
\end{equation}

 Consider the linearized system of (\ref{dUdt2})
  \begin{equation}
            \label{dUdt3}
            \dfrac{dU}{dt}=D_0\Delta U(t)+L_0U^t.
            \end{equation}
 From the previous discussion, we know that system (\ref{dUdt3}) has pure imaginary eigenvalues  $\{\pm iz_{n_1}\omega^*,\pm iz_{n_2}\omega^*\}$ at the double Hopf bifurcation point and the other eigenvalues with non-zero real part.

Assume that   the non-resonant condition holds true and use the algorithm we give in Section \ref{section3}.
 After a few calculations, we have that  the bases of $P_\varLambda$ and $P^*$, respectively, are
$\varPhi(\theta)=(\phi_1(\theta),\bar\phi_1(\theta),\phi_3(\theta),\bar\phi_3(\theta)) $,  and  $\varPsi(s)=(\psi_1(s),\bar\psi_1(s),\psi_3(s),\bar\psi_3(s))^T, $ with  \begin{equation}\label{p12}
\begin{aligned}
\phi_1(\theta)=(1,p_{12},p_{13})^Te^{iz_{n_1}\omega^*\theta},\phi_3(\theta)=(1,p_{32},p_{33})^Te^{iz_{n_2}\omega^*\theta},\\
\psi_1^*(s)=D_1(1,q_{12}^*,q_{13}^*)e^{-iz_{n_1}\omega^*s}, \psi_3^*(s)=D_2(1,q_{32}^*,q_{33}^*)e^{-iz_{n_2}\omega^*s},
\end{aligned}\end{equation}
   where
   \begin{equation*}
         \begin{array}{l}
 p_{12}=\dfrac{\mu I^*e^{-iz_{n_1}\omega^*}}{d_2\frac{n^2}{l^2}+iz_{n_1}},  ~p_{13}=0,   ~
 p_{32}=\dfrac{\mu I^*e^{-iz_{n_2}\omega^*}}{d_2\frac{n^2}{l^2}+iz_{n_2}}, ~p_{33}=0, \end{array}
      \end{equation*}\begin{equation*}
         \begin{array}{l}
 q_{12}^*=-\dfrac{-\mu S^*+\gamma}{-d_2\frac{n^2}{l^2}-iz_{n_1}},~ q_{32}^*=-\dfrac{-\mu S^*+\gamma}{-d_2\frac{n^2}{l^2}-iz_{n_2}}, \end{array}
      \end{equation*}\begin{equation*}
         \begin{array}{l}q_{13}^*=-\dfrac{\alpha-\alpha e^{-d\tau}e^{-iz_{n_1}\omega^*\tau^*}}{-2\beta y^*-d_3\frac{n^2}{l^2}+\alpha e^{-d\tau}e^{-iz_{n_1}\omega^*\tau^*}-iz_{n_1}}, \end{array}
      \end{equation*}
         \begin{equation*}
         \begin{array}{l}q_{33}^*=-\dfrac{\alpha-\alpha e^{-d\tau}e^{-iz_{n_2}\omega^*\tau^*}}{-2\beta y^*-d_3\frac{n^2}{l^2}+\alpha e^{-d\tau}e^{-iz_{n_2}\omega^*\tau^*}-iz_{n_2}}, \end{array}
      \end{equation*}\begin{equation*}
         \begin{array}{l}
 D_1=\frac{1 }{1+q_{12}^*p_{12}+\omega^*\mu I^*e^{-iz_{n_1}\omega^*}(q_{12}^*-1) },   D_2=\frac{1 }{1+q_{32}^*p_{32}+\omega^*\mu I^*e^{-iz_{n_2}\omega^*}(q_{32}^*-1) }.
           \end{array}
      \end{equation*}

Consider the Taylor expansion
\begin{equation*}
\widetilde{F}(\alpha,U^t)=\frac{1}{2!}\widetilde{F}_2(\alpha,U^t)+\frac{1}{3!}\widetilde{F}_3(\alpha,U^t),
\end{equation*}
where
$\widetilde{F}_2(\alpha,U^t)=2\alpha_1D_1^{(1,0)} \Delta U+2\alpha_2D_1^{(0,1)} \Delta U
+2\alpha_1L_1^{(1,0)}U^t+F_2( \omega^*+\alpha_1,d_2^*+\alpha_2,U^t)$, and  $\frac{1}{3!}\widetilde{F}_3(\alpha,U^t)=0$.
By a few calculations, we have
\begin{equation*}
\begin{array}{l}
F_{2000}=2\omega^*\left( \begin{array}{c}
-\mu e^{-iz_{n_1}\omega^*}p_{12}\\\mu e^{-iz_{n_1}\omega^*}p_{12}\\-\beta p_{13}^2
\end{array}\right),\\ F_{1100}=2\omega^*\left( \begin{array}{c}
-\mu( e^{-iz_{n_1}\omega^*}\overline{p}_{12}+ e^{iz_{n_1}\omega^*}p_{12})\\\mu( e^{-iz_{n_1}\omega^*}\overline{p}_{12}+ e^{iz_{n_1}\omega^*}p_{12})\\-2\beta p_{13}\overline{p}_{13}
\end{array}\right),\end{array}
\end{equation*}
\begin{equation*}
\begin{array}{l}
F_{1010}=2\omega^*\left( \begin{array}{c}
-\mu( e^{-iz_{n_1}\omega^*}p_{32}+ e^{-iz_{n_2}\omega^*}p_{12})\\\mu( e^{-iz_{n_1}\omega^*}p_{32}+ e^{-iz_{n_2}\omega^*}p_{12})\\-2\beta p_{13}p_{33}
\end{array}\right), \\ F_{1001}=2\omega^*\left( \begin{array}{c}
-\mu( e^{-iz_{n_1}\omega^*}\overline{p}_{32}+ e^{iz_{n_2}\omega^*}p_{12})\\\mu( e^{-iz_{n_1}\omega^*}\overline{p}_{32}+ e^{iz_{n_2}\omega^*}p_{12})\\-2\beta p_{13}\overline{p}_{33}
\end{array}\right),\end{array}
\end{equation*}
\begin{equation*}
\begin{array}{l}
F_{0200}=2\omega^*\left( \begin{array}{c}
-\mu e^{iz_{n_1}\omega^*}\overline{p}_{12}\\\mu e^{iz_{n_1}\omega^*}\overline{p}_{12}\\-\beta \overline{p}_{13}^2
\end{array}\right),\\ F_{0110}=2\omega^*\left( \begin{array}{c}
-\mu( e^{iz_{n_1}\omega^*}p_{32}+ e^{-iz_{n_2}\omega^*}\overline{p}_{12})\\\mu( e^{iz_{n_1}\omega^*}p_{32}+ e^{-iz_{n_2}\omega^*}\overline{p}_{12})\\-2\beta \overline{p}_{13}p_{33}
\end{array}\right),\end{array}
\end{equation*}
\begin{equation*}
\begin{array}{l}
F_{0101}=2\omega^*\left( \begin{array}{c}
-\mu( e^{iz_{n_1}\omega^*}\overline{p}_{32}+ e^{iz_{n_2}\omega^*}\overline{p}_{12})\\\mu( e^{iz_{n_1}\omega^*}\overline{p}_{32}+ e^{iz_{n_2}\omega^*}\overline{p}_{12})\\-2\beta \overline{p}_{13}\overline{p}_{33}
\end{array}\right),\\ F_{0020}=2\omega^*\left( \begin{array}{c}
-\mu e^{-iz_{n_2}\omega^*}p_{32}\\\mu e^{-iz_{n_2}\omega^*}p_{32}\\-\beta p_{33}^2
\end{array}\right), \end{array}
\end{equation*}
\begin{equation*}
\begin{array}{l}

F_{0011}=2\omega^*\left( \begin{array}{c}
-\mu( e^{-iz_{n_2}\omega^*}\overline{p}_{32}+ e^{iz_{n_2}\omega^*}p_{32})\\\mu( e^{-iz_{n_2}\omega^*}\overline{p}_{32}+ e^{iz_{n_2}\omega^*}p_{32})\\-2\beta p_{33}\overline{p}_{33}
\end{array}\right),\\ F_{0002}=2\omega^*\left( \begin{array}{c}
-\mu e^{iz_{n_2}\omega^*}\overline{p}_{32}\\\mu e^{iz_{n_2}\omega^*}\overline{p}_{32}\\-\beta \overline{p}_{33}^2
\end{array}\right).
\end{array}
\end{equation*}
By Remark \ref{S} in Section \ref{section3}, define the linear operators
\begin{equation*}
\begin{aligned}
S_{yz_1}(y)=F_{y(0)z_1}y(0)+F_{y(-1)z_1}y(-1)+F_{y(-\frac{\tau}{\omega^*})z_1}y(-\frac{\tau}{\omega^*}),\\
S_{yz_2}(y)=F_{y(0)z_2}y(0)+F_{y(-1)z_2}y(-1)+F_{y(-\frac{\tau}{\omega^*})z_2}y(-\frac{\tau}{\omega^*}),\\
S_{yz_3}(y)=F_{y(0)z_3}y(0)+F_{y(-1)z_3}y(-1)+F_{y(-\frac{\tau}{\omega^*})z_3}y(-\frac{\tau}{\omega^*}),\\
S_{yz_4}(y)=F_{y(0)z_4}y(0)+F_{y(-1)z_4}y(-1)+F_{y(-\frac{\tau}{\omega^*})z_4}y(-\frac{\tau}{\omega^*}),\\
\end{aligned}
\end{equation*}
where
\begin{equation*}
\begin{array}{ccc}
F_{y(0)z_1}=2\omega^*\left(\begin{array}{ccc}
0 &-\mu e^{-iz_{n_1}\omega^*}& 0\\0 &\mu e^{-iz_{n_1}\omega^*}&0\\0&0&-2\beta p_{13}
\end{array} \right), & F_{y(-1)z_1}=2\omega^*\left(\begin{array}{ccc}
-\mu p_{12}&0&0\\\mu p_{12}&0&0\\0&0&0
\end{array} \right),& F_{y(-\frac{\tau}{\omega^*})z_1}=0,\\
F_{y(0)z_2}=2\omega^*\left(\begin{array}{ccc}
0 &-\mu e^{iz_{n_1}\omega^*}& 0\\0 &\mu e^{iz_{n_1}\omega^*}&0\\0&0&-2\beta \overline{p}_{13}
\end{array} \right), & F_{y(-1)z_2}=2\omega^*\left(\begin{array}{ccc}
-\mu \overline{p}_{12}&0&0\\\mu \overline{p}_{12}&0&0\\0&0&0
\end{array} \right),& F_{y(-\frac{\tau}{\omega^*})z_2}=0,\end{array}
\end{equation*}
\begin{equation*}
\begin{array}{ccc}
F_{y(0)z_3}=2\omega^*\left(\begin{array}{ccc}
0 &-\mu e^{-iz_{n_2}\omega^*}& 0\\0 &\mu e^{-iz_{n_2}\omega^*}&0\\0&0&-2\beta p_{33}
\end{array} \right), & F_{y(-1)z_3}=2\omega^*\left(\begin{array}{ccc}
-\mu p_{32}&0&0\\\mu p_{32}&0&0\\0&0&0
\end{array} \right),& F_{y(-\frac{\tau}{\omega^*})z_3}=0,\\
F_{y(0)z_4}=2\omega^*\left(\begin{array}{ccc}
0 &-\mu e^{iz_{n_2}\omega^*}& 0\\0 &\mu e^{iz_{n_2}\omega^*}&0\\0&0&-2\beta \overline{p}_{33}
\end{array} \right), & F_{y(-1)z_4}=2\omega^*\left(\begin{array}{ccc}
-\mu \overline{p}_{32}&0&0\\\mu \overline{p}_{32}&0&0\\0&0&0
\end{array} \right),& F_{y(-\frac{\tau}{\omega^*})z_4}=0.\\
\end{array}
\end{equation*}
Following the steps in Remark \ref{calculationsteps}, we can get all the coefficients in (\ref{normalform}), and then obtain the normal forms up to third order.

         \subsection{Simulations}
In the following, we take $\alpha=2.1$, $d=0.5$, $\mu=0.5$, $\gamma=0.1$, $\beta=0.3$, $\tau=1$, $d_1=0.05$, $d_3=0.06$, and let $\omega$ and $d_2$ be the bifurcation parameters.

\begin{figure}
             \centering
                    a)\includegraphics[width=0.45\textwidth,height=0.35\textwidth]{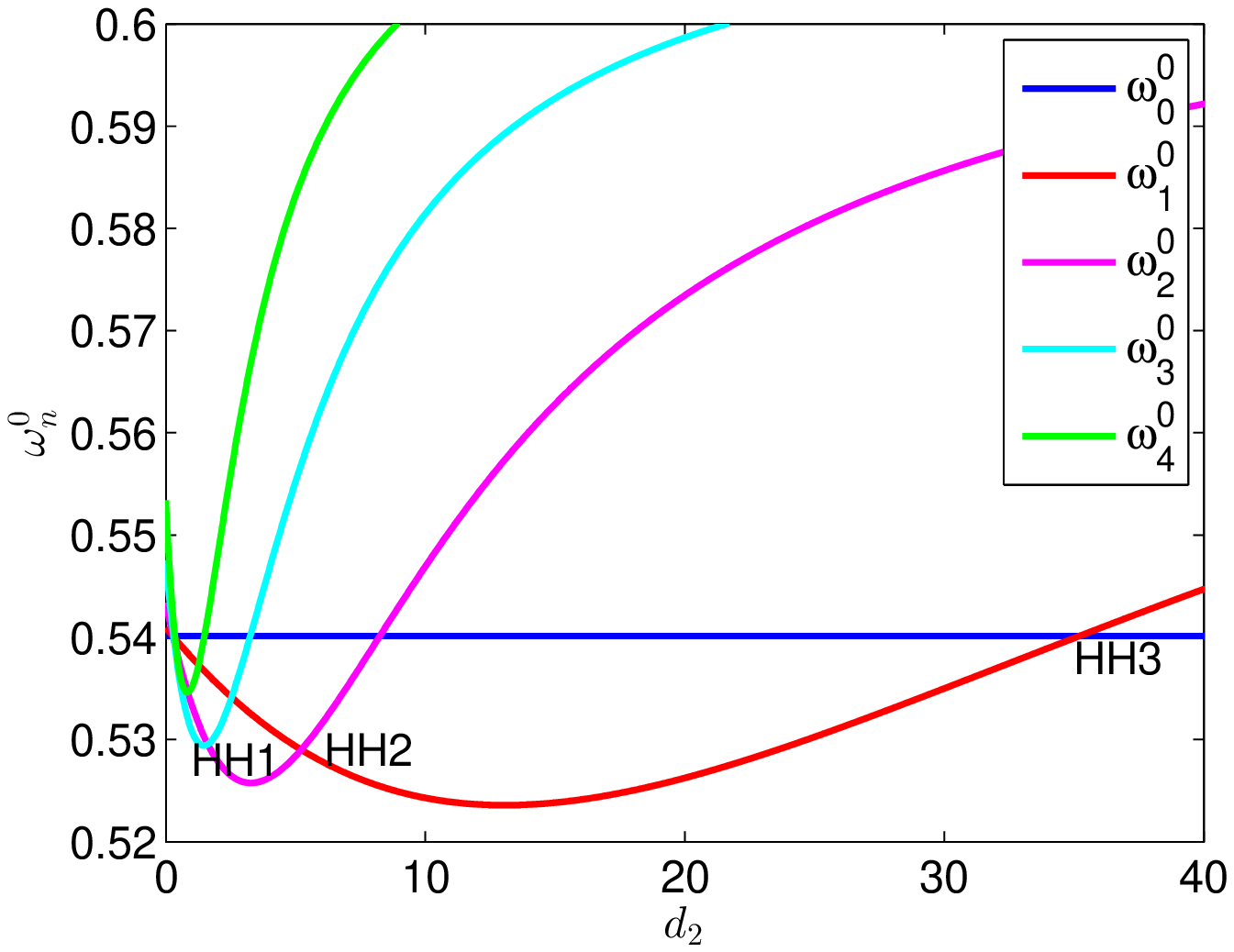}
                     b)\includegraphics[width=0.49\textwidth,height=0.35\textwidth]{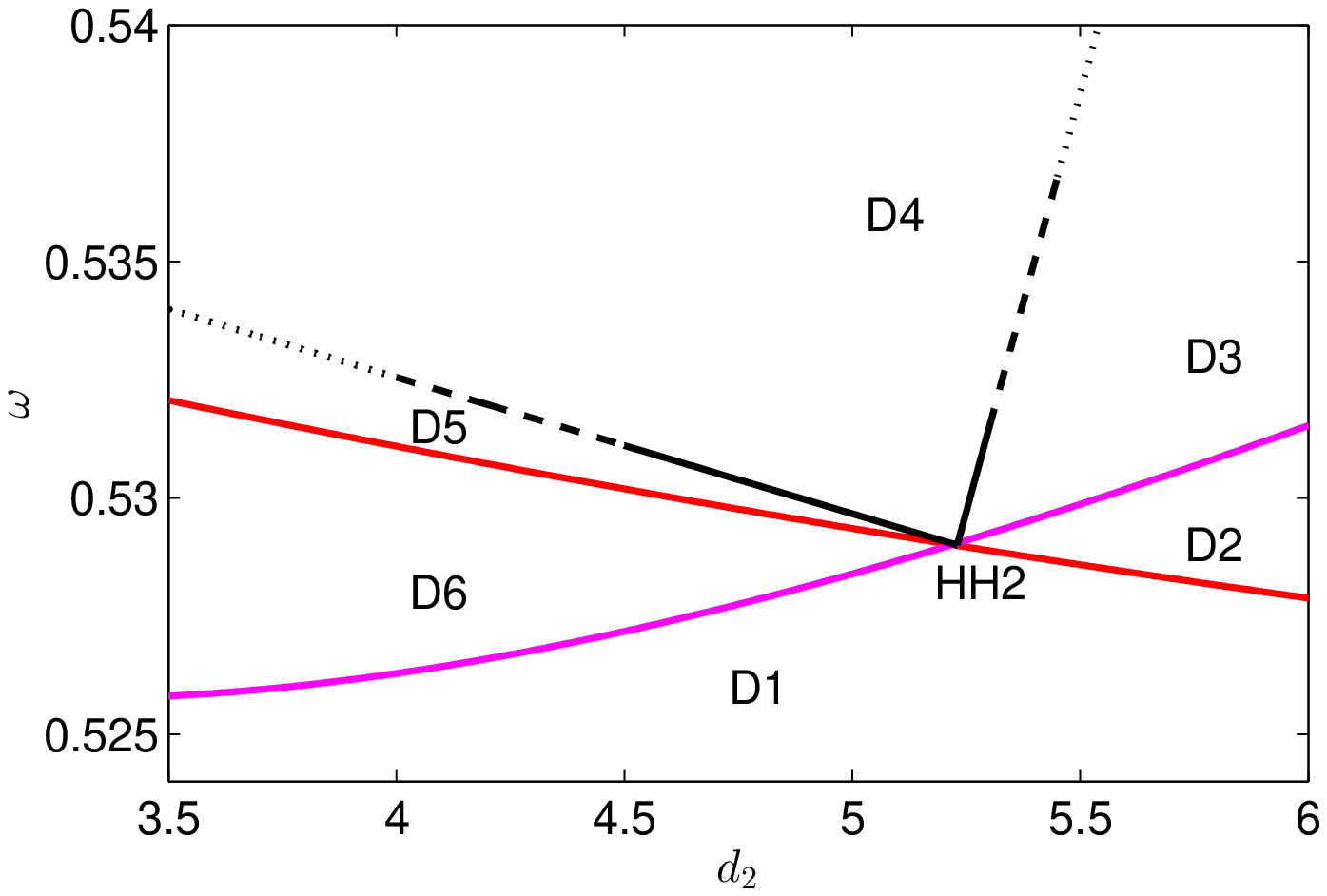}                                   
            \caption  {a) When $\alpha=2.1$, $d=0.5$, $\mu=0.5$, $\gamma=0.1$, $\beta=0.3$, $\tau=1$, $d_1=0.05$, $d_3=0.06$, the  bifurcation set on the $d_2-\omega$ plane is drawn, and  double
            Hopf bifurcation points are marked. b) Complete bifurcation set
                        near the double Hopf point HH2.}
           \label{fig:d2omega}
             \end{figure}

\begin{figure}
             \centering
     \includegraphics[width=0.95\textwidth]{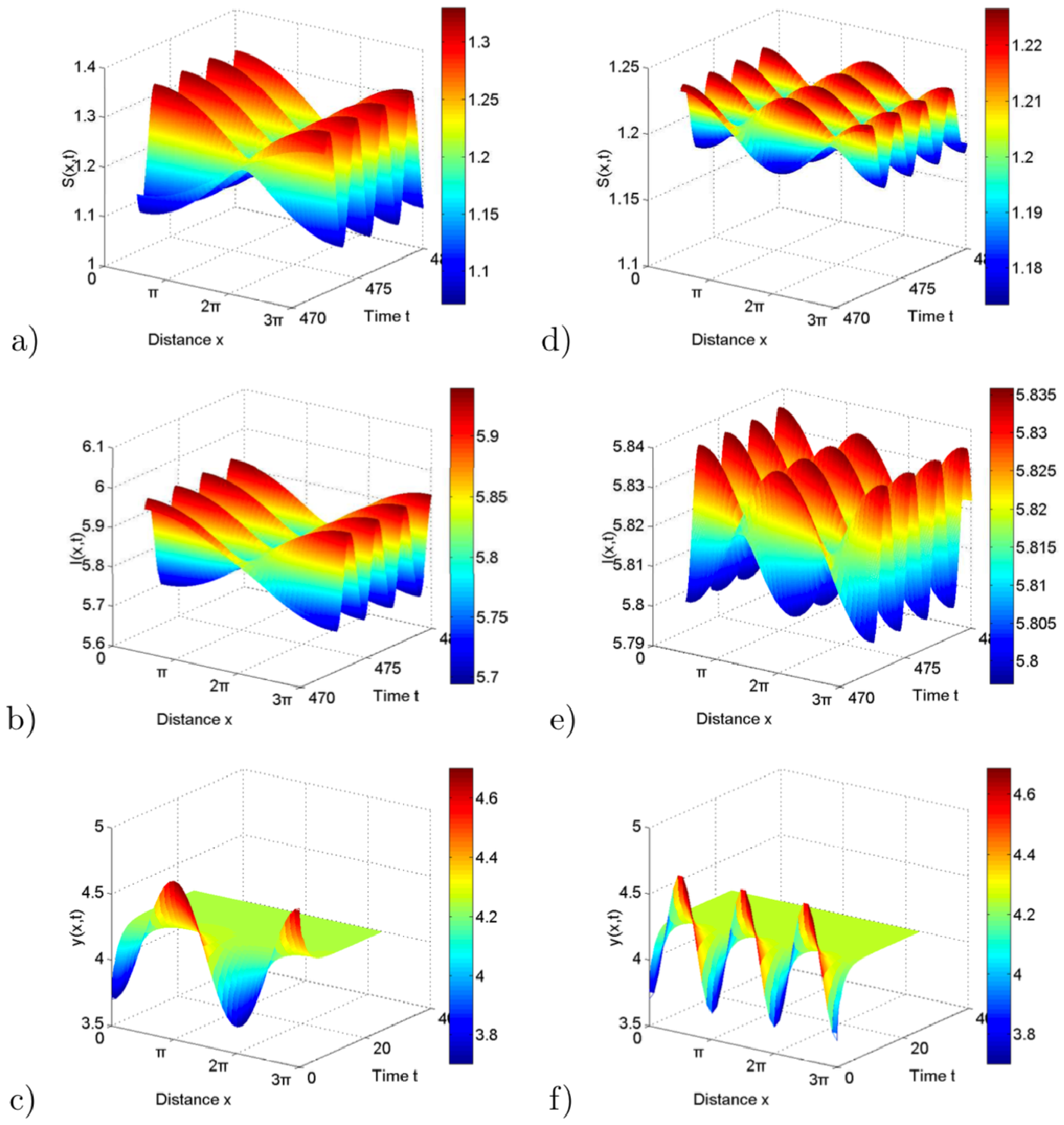}       
          \caption{ When  $d_2=5.23$, and $\omega=0.53$ in $D4$, two stable nonhomogeneous periodic solutions coexist when we choose two different initial conditions. For (a-c)  initial conditions are $S(x,t)=1.2+0.01\cos x$, $I(x,t)=5.8-0.06\cos x$, $y(x,t)=4.2-0.05\cos x$, $t\in [-\tau,0]$. For (d-f)  initial conditions are $S(x,t)=1.2+0.01\cos 2x$, $I(x,t)=5.8-0.06\cos 2x$, $y(x,t)=4.2-0.05\cos 2x$, $t\in [-\tau,0]$.  }
         \label{fig:hh2D4}
        \end{figure}

By (\ref{omeganj}), we can draw  the curves of Hopf bifurcation values when   $d_2$ varies.
 As shown clearly in Fig. \ref{fig:d2omega} a), every two Hopf bifurcation curves intersect at a double Hopf bifurcation point such as HH1, HH2 and HH3.
   When  $d_2=5.23$, $\omega_1^0$ intersects $\omega_2^0$, and we denote the double Hopf bifurcation point HH2. It means that the wave numbers are the case of $n_1=1,n_2=2$.  For HH2, we have $z_{n_1}=z_1=2.9930$, $z_{n_2}=z_2=3.1037$, and $\omega_{1}^0=\omega_{2}^0=0.5290$. Using the algorithm we established, we get that $B_{11} =0.7184 + 0.5138i$, $B_{21} =0.0021 + 0.0042i$,
   $B_{13} =0.6431 + 0.5805i$, $B_{23} =-0.0037 + 0.0101i$,
$B_{2100} =  -0.0001 - 0.1942i$, $B_{1011} =  -0.0010 - 0.3398i$, $B_{0021} =  -0.00071 + 0.00055i$, and $B_{1110 }=  -0.0851 - 0.4349i$. Thus, $\epsilon_1 =   -1$, $\epsilon_2 =  -1$, $b_0 =   1.4401$, $c_0 =  0.0016$,
$d_0 =    1$, $d_0-b_0c_0 = -0.0023$.    It means that the unfolding system is of type Ib. Moreover, we have $c_1=-0.7184\alpha_1-0.0021\alpha_2 $, $c_2=-0.6431\alpha_1+0.0037\alpha_2 $. Thus,    the bifurcation set  near HH2 is shown in Figure \ref{fig:d2omega} b), in which the two black lines are two pitchfork bifurcation curves   $\omega=(d_2-5.23)/(    28.2318)+0.5290$,  $\omega=(d_2-5.23)/(   -346.6577)+0.5290$.   The parameter plane near the bifurcation point is divided into six regions. In D1, the positive equilibrium is asymptotically stable. In D2 or D6, there are stable periodic solutions. In region D4, there are two stable nonhomogeneous periodic solutions coexisting, which are corresponding to two different eigenfunctions.
 From Theorem 2.2 in Section 6.2 of \cite{JWu}, Hopf bifurcating periodic solutions of system (\ref{dUdt1}) can be parameterized by a small parameter $\epsilon$. Thus, when $\nu=\nu(\epsilon)$ and $\epsilon$ are near 0, (i.e. when $\omega$ is near $\omega^*$),  two periodic solutions  have the following representations respectively
$$ \begin{aligned}
 U_t(\nu,\theta)(x)=\epsilon Re\phi_1(\theta)e^{iz_{n_1}\omega^*t}\cos\frac{n_1}{l}x+O(\epsilon^2),\\
    U_t(\nu,\theta)(x)=\epsilon Re\phi_3(\theta)e^{iz_{n_2}\omega^*t}\cos\frac{n_2}{l}x+O(\epsilon^2),
 \end{aligned}$$
where $\phi_1(\theta)$ and $\phi_3(\theta)$ is defined in (\ref{p12}). Figure \ref{fig:hh2D4} illustrates that the two solutions have totally different spatial shape.

For  double Hopf bifurcation points HH1 and HH3, we can calculate the parameters in Table 2. Obviously, these two points are also of type Ib.  We notice that during the calculation of step 3 in third order normal form of the three double Hopf bifurcation points,  the cases $n_1=0, n_2\neq 0$ (e.g. HH3), and $ n_1\neq 0,  n_2\neq 0$ (e.g. HH1, HH2)   both  occur.

 \begin{table}[tbp]
  \label{twopairs}
  \caption{Parameter values at double Hopf bifurcation points.}
  \centering
  \begin{tabular}{ccccccccccc}
  \hline
  $\begin{array}{c}
  Point
  \end{array}$ & $d_2^*$&$\omega^*$ & $n_1$ & $z_{n_1}$  & $n_2$& $z_{n_2}$ &$b_0$ &$c_0$ &$d_0$& $d_0-b_0c_0$     \\ \hline  
  HH1 & $1.62$& $0.5295$& 2& 3.0071& 3& 3.0763&    1.8774&    1.1651 & 1&  -1.1872\\
               HH2& 5.23 & 0.5290& 1& 2.9930 & 2& 3.1037&    1.4401&    0.0016& 1&    -0.0023\\
              HH3& 35.2& 0.5401 & 0 & 2.9082 & 1&  3.1007&   3.1290
              &            1.4043
               & 1&  -3.3939

       \\ \hline
  \end{tabular}
  \end{table}

\section{Application to a Diffusive Predator-Prey Model with Delay} \label{section predatorprey}
In this section,  a diffusive predator-prey model with delay is considered. Taking $r_1$ and $\tau$ as bifurcation parameters, we can find the double Hopf bifurcation point.  By the steps in Remark \ref{calculationsteps}, we can obtain the normal form and get the unfolding system. We will show that the system will exhibit complex dynamical behavior near the double Hopf bifurcation point: the existence of quasi-periodic solution on a 2-torus,  quasi-periodic solution on a 3-torus, and even  chaos.
  \subsection{Model and the Existence of Double Hopf Bifurcation}
 Let us consider the predator-prey system
 \begin{equation}\label{predator1}\left\{
 \begin{aligned}
 &\frac{\partial X(x,t)}{\partial t}=d_1\Delta X(x,t)+X(x,t)[r_1-a_{11}X(x,t-\tau)-a_{12}Y(x,t)],\\
 &\frac{\partial Y(x,t)}{\partial t}=d_2\Delta Y(x,t)+Y(x,t)[-r_2+a_{21}X(x,t)-a_{22}Y(x,t)],
 \end{aligned}\right.x\in(0,l\pi),
 \end{equation}
 equipping with homogeneous Neumann boundary condition, where  $X(x,t)$ and $Y(x,t)$ represent the densities of  prey and predator populations at time $t$ and location $x$, respectively, $d_1, d_2>0$ denote the diffusion coefficients of prey and predator,   $\tau$ is the generation time of the prey species, and $r_i, a_{ij}$ $ (i,j=1,2)$ are positive constants.

 It is obvious that if $r_1a_{21}-r_2a_{11}>0$,  (\ref{predator1}) has a unique positive equilibrium $E^*=(X^*,Y^*)$ \cite{Y SongHopf}
 where
 \[X^*=\frac{r_1a_{22}+r_2a_{12}}{a_{11}a_{22}+a_{12}a_{21}},~Y^*=\frac{r_1a_{21}-r_2a_{11}}{a_{11}a_{22}+a_{12}a_{21}}.\]
 By the tanslation $u=X-X^*,v=Y-Y^*$, (\ref{predator1}) can be written as
 \begin{equation}\label{predatorpingyi}
 \begin{aligned}
 &\frac{\partial u(x,t)}{\partial t}=d_1\Delta u(x,t)+(u(x,t)+X^*)[-a_{11}u(x,t-\tau)-a_{12}v(x,t)],\\
 &\frac{\partial v(x,t)}{\partial t}=d_2 \Delta v(x,t)+(v(x,t)+Y^*)[a_{21}u(x,t)-a_{22}v(x,t)].
 \end{aligned}
 \end{equation}
 The linearization of (\ref{predatorpingyi}) at the origin is
 \begin{equation*}\label{predatorlinear}
 \begin{aligned}
 &\frac{\partial u(x,t)}{\partial t}=d_1\Delta u(x,t)-a_{11}X^*u(x,t-\tau)-a_{12}X^*v(x,t)],\\
 &\frac{\partial v(x,t)}{\partial t}=d_2\Delta v(x,t)+a_{21}Y^*u(x,t)-a_{22}Y^*v(x,t),
 \end{aligned}
 \end{equation*}
 whose characteristic equation is  \begin{equation}
       \label{characterpredator}
          \lambda^2+A_n\lambda+B_n+ e^{-\lambda\tau}(C\lambda+D_n) =0,
       \end{equation}
       with $n\in \mathbb{N}_0$,
       \begin{equation*}
       \begin{array}{ll}
      A_n=d_1\frac{n^2}{l^2}+d_2\frac{n^2}{l^2}+a_{22}Y^*, \\
         B_n=d_1\frac{n^2}{l^2}(a_{22}Y^*+d_2\frac{n^2}{l^2})+a_{12}a_{21}X^*Y^*, \\       C=a_{11}X^*,          \\     D_n=(a_{22}Y^*+d_2\frac{n^2}{l^2})a_{11}X^*.
                \end{array}
    \end{equation*}

    Clearly, $\lambda=0$ is not a root of (\ref{characterpredator}), which excludes the existence of Turing bifurcation.
    When $\tau=0$,  Eq. (\ref{characterpredator}) becomes the following sequence of quadratic polynomial equations
  \begin{equation}
  \label{characterpredatoromega0}
   \lambda^2+(A_n+C)\lambda+(B_n+D_n)=0,~~~~~~n\in \mathbb{N}_0,
  \end{equation}
 where $A_n+C=d_1\frac{n^2}{l^2}+d_2\frac{n^2}{l^2}+a_{22}Y^*+a_{11}X^* >0, ~B_n+D_n=d_1\frac{n^2}{l^2}(a_{22}Y^*+d_2\frac{n^2}{l^2})+a_{12}a_{21}X^*Y^*+(a_{22}Y^*+d_2\frac{n^2}{l^2})a_{11}X^*>0.$
  We know that all roots of Eq. (\ref{characterpredatoromega0}) have negative real part.

  We would like to seek critical values of $\tau$ such that there exists a pair of imaginary eigenvalues. Let $\pm i\omega (\omega>0)$ be the solution of Eq. (\ref{characterpredator}), then we have
       \begin{eqnarray*}
      \label{shixupredator}
      \left\{
      \begin{array}{l}
      -\omega^2+B_n= -C\omega\sin \omega\tau -D_n\cos \omega\tau ,\\
      A_n\omega=-C\omega\cos  \omega\tau+D_n\sin  \omega\tau ,
      \end{array}
      \right.
      \end{eqnarray*}
   which is solved by
    \begin{equation*}\label{SnCnpredator}
               \begin{array}{l}
       \sin  \omega\tau=\dfrac{A_n\omega D_n-(B_n-\omega^2)C\omega}{(C\omega)^2+D_n^2}=S_n(\omega),\\ \cos  \omega\tau =-\dfrac{A_nC\omega^2+(B_n-\omega^2)D_n}{D_n^2+(C\omega)^2}=C_n(\omega).
                \end{array}
                \end{equation*}
     Then, we obtian
       \begin{equation}
       \label{fomega}
      G(\omega)=\omega^4+(A_n^2-2B_n-C^2)\omega^2+B_n^2-D_n^2=0.
       \end{equation}

 Suppose that
 \begin{equation*}
 (H_5)~ A_n^2-2B_n-C^2<0, ~B_n^2-D_n^2>0, ~and~  (A_n^2-2B_n-C^2)^2-4(B_n^2-D_n^2)>0
 \end{equation*}
 holds, Eq. (\ref{fomega}) has two positive roots given by
 \begin{equation*}
 \omega_n^{\pm}=\sqrt{\frac{-(A_n^2-2B_n-C^2)\pm\sqrt{(A_n^2-2B_n-C^2)^2-4(B_n^2-D_n^2)}}{2}}.
 \end{equation*}
 \begin{equation}\label{taunj}
 \tau_n^{j\pm}=\left\lbrace \begin{array}{ll}
 \frac{\arccos C_n(\omega)+2j\pi}{\omega_n^{\pm}},&\rm{if} ~ S_n(\omega_n^{\pm})>0,\\ \frac{2\pi-\arccos C_n(\omega)+2j\pi}{\omega_n^{\pm}},&\rm{if} ~ S_n(\omega_n^{\pm})<0.\\
 \end{array}\right.
 \end{equation}

  By calculation and the results in \cite{Y SongHopf}, we can verify the  transversality condition     $$\mathrm{Sign}\mathrm{Re}\frac{d\lambda}{d\tau}\mid_{\lambda=i\omega_n^+}>0,$$ and $$\mathrm{Sign}\mathrm{Re}\frac{d\lambda}{d\tau}\mid_{\lambda=i\omega_n^-}<0.$$
 \begin{theorem}
  Suppose that $(H_5)$ holds, system (\ref{predator1}) undergoes a Hopf bifurcation at the origin when $\tau=\tau_{n}^{j-}$ or $\tau_{n}^{j+}$.
 \end{theorem}
 We fix the other parameters, and choose $r_1$ and $\tau$ as double Hopf bifurcation parameters.
 \begin{theorem}
  Suppose that $(H_5)$ holds,  and there exists $r_1^*$ such that $\tau_{n_1}^{j_1-}=\tau_{n_2}^{j_2+}$. Then system (\ref{predator1}) undergoes a double Hopf bifurcation at the origin when $r_1=r_1^*$, $\tau=\tau_{n_1}^{j_1-}=\tau_{n_2}^{j_2+}\stackrel{\bigtriangleup}{=}\tau^*$.
 \end{theorem}
  \subsection{Normal Form of Double Hopf Bifurcation}
 Let   $u(x,t)=X(x,\tau t)-X^*,v(x,t)=Y(x,\tau t)-Y^*$,  $\tau=\tau^*+\alpha_1, r_1=r_1^*+\alpha_2$,  then (\ref{predator1}) can be written as
 \begin{equation*}\label{predator}
  \begin{aligned}
  &\frac{\partial u(x,t)}{\partial t}=\tau^*d_1\Delta u(x,t)+\tau^*(u(x,t)+X^*)[-a_{11}u(x,t-1)-a_{12}v(x,t)],\\
  &\frac{\partial v(x,t)}{\partial t}=\tau^*d_2\Delta v(x,t)+\tau^*(v(x,t)+Y^*)[a_{21}u(x,t)-a_{22}v(x,t)].
  \end{aligned}
  \end{equation*}
  Denote $U(t)=(u(x,t),v(x,t)^T$, then we have
   \begin{equation}
      \label{dUdtpredator}
      \dfrac{{\rm d} U}{{\rm d} t}=D(\tau^*+\alpha_1,r_1^*+\alpha_2)\Delta U(t)+L(\tau^*+\alpha_1,r_1^*+\alpha_2)U^t+F(\tau^*+\alpha_1,r_1^*+\alpha_2,U^t).
      \end{equation}
   Here
     \begin{equation*}
       \begin{array}{l}
      D(\tau^*+\alpha_1,r_1^*+\alpha_2)= (\tau^*+\alpha_1)\left( \begin{array}{cc}
           d_1&0\\
           0&d_2\end{array}\right)=D_0+\alpha_1D_1^{(1,0)}+\alpha_2D_1^{(0,1)},\\
     L(\tau^*+\alpha_1,r_1^*+\alpha_2)U^t=(L_0+\alpha_1L_1^{(1,0)}+\alpha_2L_1^{(0,1)})U^t,\\
     F(\omega^*+\alpha_1,r_1^*+\alpha_2,U^t)=(\tau^*+\alpha_1)\left( \begin{array}{c}
     -a_{11} U_1^t(-1)U_1^t(0)-a_{12}U_1^t(0)U_2^t(0)\\
     a_{21}U_1^t(0)U_2^t(0)-a_{22}U_2^{t2}     \end{array}\right),
       \end{array}
          \end{equation*}
 with
 \[D_0=\tau^*\left( \begin{array}{cc}
            d_1&0\\
            0&d_2\end{array}\right), D_1^{(1,0)}=\left( \begin{array}{cc}
                        d_1&0\\
                        0&d_2\end{array}\right),D_1^{(0,1)}=0,\]
\[L_0U^t=\tau^*(AU^t(0)+G_1U^t(-1)), L_1^{(1,0)}U^t=AU^t(0)+G_1U^t(-1),\]
\[L_1^{(0,1)}U^t=\tau^*\frac{1}{a_{11}a_{22}+a_{12}a_{21}}\left[ \left( \begin{array}{cc}
0&-a_{12}a_{22}\\a_{21}a_{21}&-a_{22}a_{21}
\end{array}\right)U^t(0)+\left(\begin{array}{cc}
-a_{11}a_{22}&0\\0&0
\end{array}\right) U^t(-1) \right] , \]
where
\[A=\left(\begin{array}{cc}
0&-a_{12}X_{r_1^*}^*\\a_{21}Y_{r_1^*}^*&-a_{22}Y_{r_1^*}^*
\end{array} \right),
G_1=\left( \begin{array}{cc}
-a_{11}X_{r_1^*}^*&0\\0&0
\end{array}\right) , \]  and $X_{r_1^*}^*=\frac{r_1^*a_{22}+r_2a_{12}}{a_{11}a_{22}+a_{12}a_{21}}, Y_{r_1^*}^*=\frac{r_1^*a_{21}-r_2a_{11}}{a_{11}a_{22}+a_{12}a_{21}}$.

        Eq. (\ref{dUdtpredator}) can be rewritten as
           \begin{equation}
                \label{dUdt2predator}
                \dfrac{{\rm d} U}{{\rm d}t}= D_0\Delta U(t)+L_0U^t+\widetilde{F}(\alpha,U^t),
                \end{equation}
           where
          \begin{equation*}
          \begin{aligned}
           \widetilde{F}(\alpha,U^t)=(\alpha_1D_1^{(1,0)}+ \alpha_2D_1^{(0,1)}) \Delta U
                    +(\alpha_1L_1^{(1,0)}+\alpha_2L_1^{(0,1)})U^t+F( \tau^*+\alpha_1,r_1^*+\alpha_2,U^t).
          \end{aligned}
                   \end{equation*}

           Consider the linearized system of (\ref{dUdt2predator})
            \begin{equation}
                      \label{dUdt3predator}
                      \dfrac{{\rm d} U}{{\rm d}t}=D_0\Delta U(t)+L_0U^t.
                      \end{equation}
  From the previous discussion, we know that system (\ref{dUdt3predator}) has pure imaginary eigenvalues  $\pm i\omega_{n_2}^+\tau^*,\pm i\omega_{n_1}^-\tau^*$ at the double Hopf bifurcation point and the other eigenvalues with non-zero real part. 
 After a few calculations, we have that  the bases of $P$ and $P^*$, respectively, are
   $\varPhi(\theta)=(\phi_1(\theta),\bar\phi_1(\theta),\phi_3(\theta),\bar\phi_3(\theta)) $,  and  $\varPsi(s)=( \psi_1(s),\bar\psi_1(s),$ $\psi_3(s),\bar\psi_3(s)) ^T, $ with
   \[\phi_1(\theta)=(1,p_{12})^Te^{i\omega_{n_1}^-\tau^*\theta},~~~
   \phi_3(\theta)=(1,p_{32})^Te^{i\omega_{n_2}^+\tau^*\theta},\]
  \[
 \psi_1^*(s)=D_1(1,q_{12})e^{-i\omega_{n_1}^-\tau^*s},~~ \psi_3^*(s)=D_2(1,q_{32})e^{-i\omega_{n_2}^+\tau^*s},  \]
    where
    \[p_{12}=\frac{-a_{11}X_{r_1}^*e^{-i\omega_{n_1}^-\tau^*}-i\omega_{n_1}^--d_1\frac{n^2}{l^2}}{a_{12}X_{r_1^*}^*}, p_{32}=\frac{-a_{11}X_{r_1}^*e^{-i\omega_{n_2}^+\tau^*}-i\omega_{n_2}^+-d_1\frac{n^2}{l^2}}{a_{12}X_{r_1^*}^*},\]
    \[q_{12}=-\frac{-a_{11}X_{r_1}^*e^{-i\omega_{n_1}^-\tau^*}-i\omega_{n_1}^--d_1\frac{n^2}{l^2}}{a_{21}Y_{r_1^*}^*},q_{32}=-\frac{-a_{11}X_{r_1}^*e^{-i\omega_{n_2}^+\tau^*}-i\omega_{n_2}^+-d_1\frac{n^2}{l^2}}{a_{21}Y_{r_1^*}^*},\]
   \[  D_1=\frac{1 }{1+p_{12}q_{12}-a_{11}X_{r_1^*}^*\tau^*e^{-i\omega_{n_1}^-\tau^*}}, D_2=\frac{1 }{1+p_{32}q_{32}-a_{11}X_{r_1^*}^*\tau^*e^{-i\omega_{n_2}^+\tau^*}}.  \]

   Consider the Taylor expansion
 \begin{equation*}
 \widetilde{F}(\alpha,U^t)=\frac{1}{2!}\widetilde{F}_2(\alpha,U^t)+\frac{1}{3!}\widetilde{F}_3(\alpha,U^t),
 \end{equation*}
 where
 $\widetilde{F}_2(\alpha,U^t)=2\alpha_1D_1^{(1,0)} \Delta U+2\alpha_2D_1^{(0,1)} \Delta U
 +(2\alpha_1L_1^{(1,0)}+2\alpha_2L_1^{(0,1)})U^t+2F( \omega^*+\alpha_1,d_2^*+\alpha_2,U^t)$, and  $\frac{1}{3!}\widetilde{F}_3(\alpha,U^t)=0$.

By a few calculations, we have
\begin{equation*}
\begin{array}{ll}
F_{2000}=2\tau^*\left( \begin{array}{c}
-a_{11} e^{-i\omega_{n_1}^-\tau^*}-a_{12}p_{12}\\a_{21}p_{12}-a_{22}p_{12}^2
\end{array}\right),\\ F_{1100}=2\tau^*\left( \begin{array}{c}
-a_{11}(e^{i\omega_{n_1}^-\tau^*}+e^{-i\omega_{n_1}^-\tau^*})-a_{12}(p_{12}+\overline{p}_{12})\\a_{21}(p_{12}+\overline{p}_{12})-a_{22}2p_{12}\overline{p}_{12}
\end{array}\right),\end{array}
\end{equation*}
\begin{equation*}
\begin{array}{ll}
F_{1010}=2\tau^*\left( \begin{array}{c}
-a_{11}(e^{-i\omega_{n_2}^+\tau^*}+e^{-i\omega_{n_1}^-\tau^*})-a_{12}(p_{12}+p_{32})\\a_{21}(p_{12}+p_{32})-a_{22}2p_{12}p_{32}
\end{array}\right),\\ F_{1001}=2\tau^*\left( \begin{array}{c}
-a_{11}(e^{i\omega_{n_2}^+\tau^*}+e^{-i\omega_{n_1}^-\tau^*})-a_{12}(p_{12}+\overline{p}_{32})\\a_{21}(p_{12}+\overline{p}_{32})-a_{22}2p_{12}\overline{p}_{32}
\end{array}\right),\end{array}
\end{equation*}
\begin{equation*}
\begin{array}{ll}
F_{0200}=2\tau^*\left( \begin{array}{c}
-a_{11} e^{i\omega_{n_1}^-\tau^*}-a_{12}\overline{p}_{12}\\a_{21}\overline{p}_{12}-a_{22}\overline{p}_{12}^2
\end{array}\right),\\ F_{0110}=2\tau^*\left( \begin{array}{c}
-a_{11}(e^{-i\omega_{n_2}^+\tau^*}+e^{i\omega_{n_1}^-\tau^*})-a_{12}(\overline{p}_{12}+p_{32})\\a_{21}(\overline{p}_{12}+p_{32})-a_{22}2\overline{p}_{12}p_{32}
\end{array}\right),\end{array}
\end{equation*}
\begin{equation*}
\begin{array}{ll}
F_{0101}=2\tau^*\left( \begin{array}{c}
-a_{11}(e^{i\omega_{n_2}^+\tau^*}+e^{i\omega_{n_1}^-\tau^*})-a_{12}(\overline{p}_{12}+\overline{p}_{32})\\a_{21}(\overline{p}_{12}+\overline{p}_{32})-a_{22}2\overline{p}_{12}\overline{p}_{32}
\end{array}\right),\\
F_{0020}=2\tau^*\left( \begin{array}{c}
-a_{11} e^{-i\omega_{n_2}^+\tau^*}-a_{12}p_{32}\\a_{21}p_{32}-a_{22}p_{32}^2
\end{array}\right),\end{array}
\end{equation*}
\begin{equation*}
\begin{array}{ll}
F_{0011}=2\tau^*\left( \begin{array}{c}
-a_{11}(e^{i\omega_{n_2}^+\tau^*}+e^{-i\omega_{n_2}^+\tau^*})-a_{12}(p_{32}+\overline{p}_{32})\\a_{21}(p_{32}+\overline{p}_{32})-a_{22}2p_{32}\overline{p}_{32}
\end{array}\right),\\ F_{0002}=2\tau^*\left( \begin{array}{c}
-a_{11} e^{i\omega_{n_2}^+\tau^*}-a_{12}\overline{p}_{32}\\a_{21}\overline{p}_{32}-a_{22}\overline{p}_{32}^2
\end{array}\right).
\end{array}
\end{equation*}
By Remark \ref{S} in Section \ref{section3}, establish the linear operators
\begin{equation*}
\begin{aligned}
S_{yz_1}(y)=F_{y(0)z_1}y(0)+F_{y(-1)z_1}y(-1),\\
S_{yz_2}(y)=F_{y(0)z_2}y(0)+F_{y(-1)z_2}y(-1),\\
S_{yz_3}(y)=F_{y(0)z_3}y(0)+F_{y(-1)z_3}y(-1),\\
S_{yz_4}(y)=F_{y(0)z_4}y(0)+F_{y(-1)z_4}y(-1),\\
\end{aligned}
\end{equation*}
where
\begin{equation*}
\begin{array}{cc}
F_{y(0)z_1}=2\tau^*\left(\begin{array}{cc}
-a_{11} e^{-i\omega_{n_1}^-\tau^*}-a_{12}p_{12}& -a_{12}\\a_{21}p_{12} &a_{21}-2a_{22}p_{12}
\end{array} \right), & F_{y(-1)z_1}=2\tau^*\left(\begin{array}{ccc}
-a_{11}&0\\0&0\end{array} \right),\\
F_{y(0)z_2}=2\tau^*\left(\begin{array}{cc}
-a_{11} e^{i\omega_{n_1}^-\tau^*}-a_{12}\overline{p}_{12}& -a_{12}\\a_{21}\overline{p}_{12} &a_{21}-2a_{22}\overline{p}_{12}
\end{array} \right), & F_{y(-1)z_2}=2\tau^*\left(\begin{array}{ccc}
-a_{11}&0\\0&0\end{array} \right),\\
F_{y(0)z_3}=2\tau^*\left(\begin{array}{cc}
-a_{11} e^{-i\omega_{n_2}^+\tau^*}-a_{12}p_{32}& -a_{12}\\a_{21}p_{32} &a_{21}-2a_{22}p_{32}
\end{array} \right), & F_{y(-1)z_3}=2\tau^*\left(\begin{array}{ccc}
-a_{11}&0\\0&0\end{array} \right),\\
F_{y(0)z_4}=2\tau^*\left(\begin{array}{cc}
-a_{11} e^{i\omega_{n_2}^+\tau^*}-a_{12}\overline{p}_{32}& -a_{12}\\a_{21}\overline{p}_{32} &a_{21}-2a_{22}\overline{p}_{32}
\end{array} \right),& F_{y(-1)z_4}=2\tau^*\left(\begin{array}{ccc}
-a_{11}&0\\0&0\end{array} \right).\\
\end{array}
\end{equation*}
Following the steps in Remark \ref{calculationsteps}, we can get all the coefficients in (\ref{normalform}), and then obtain the normal forms up to the third order.

              \begin{figure}
               \centering
       a)\includegraphics[width=0.47\textwidth]{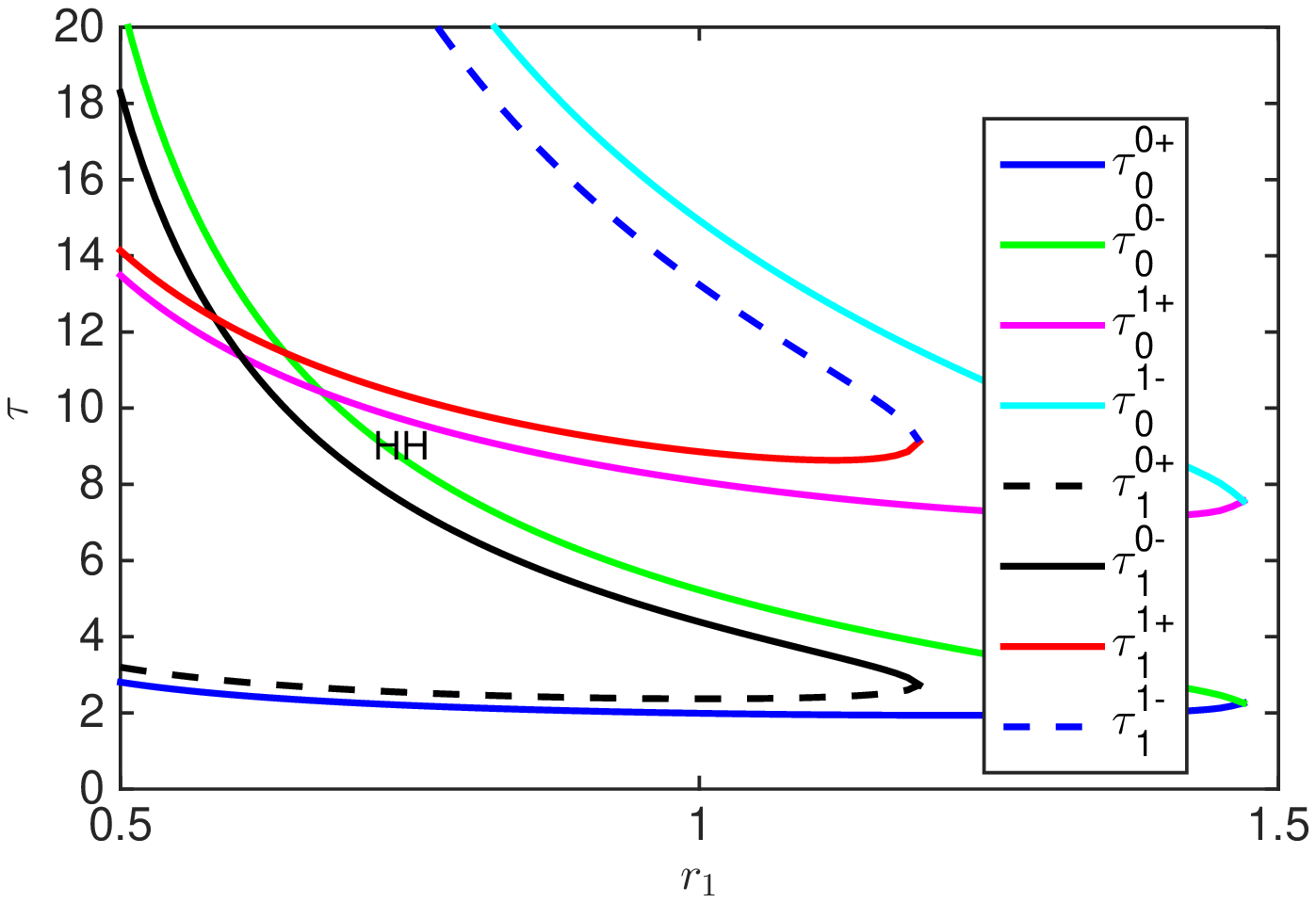}b)\includegraphics[width=0.47\textwidth]{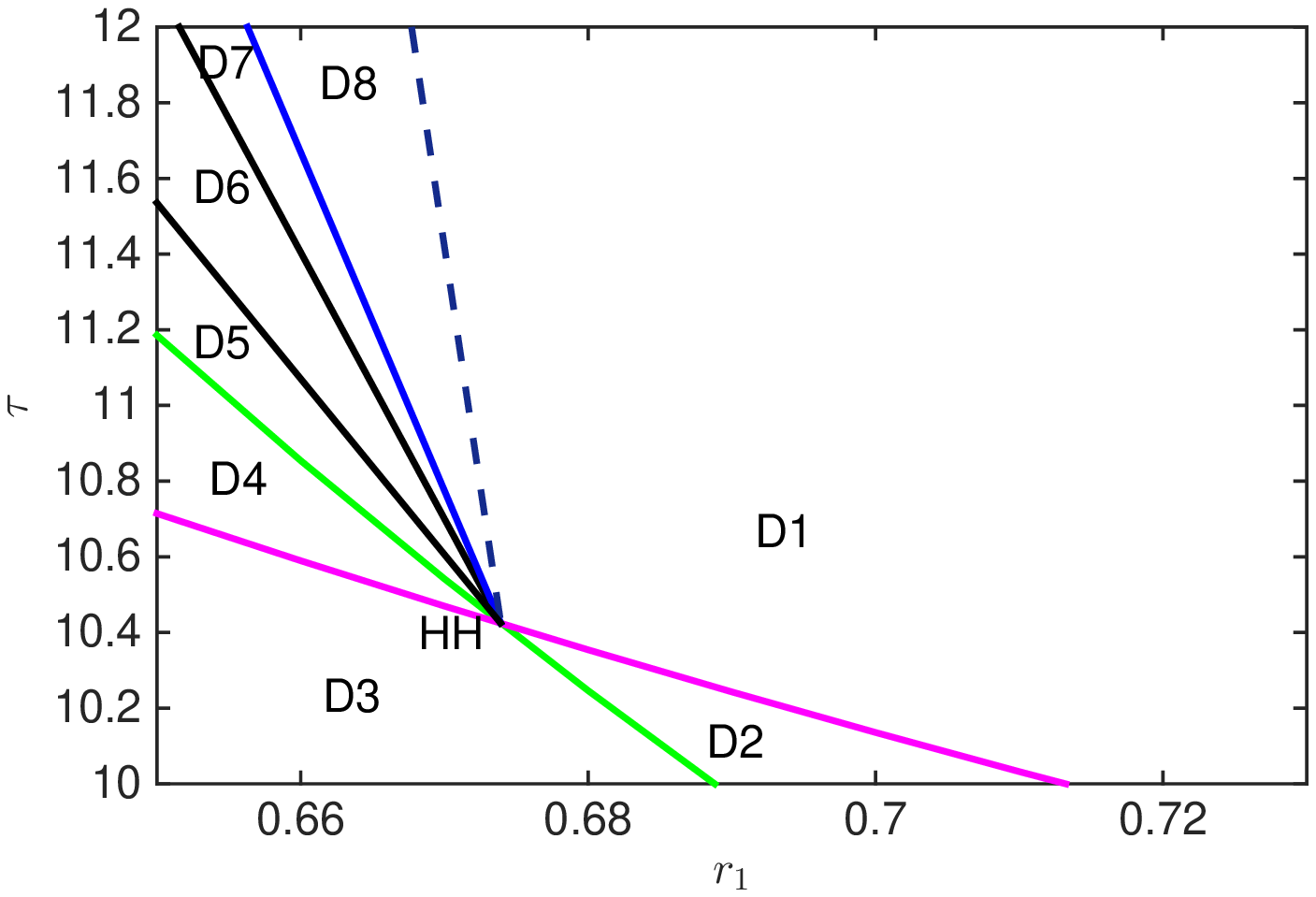}
    \caption  {  When $r_2=1,a_{11}=1,a_{12}=1.2,a_{21}=2.8,a_{22}=1,d_1=0.1,
    d_2=0.2,l=3$, (a) the partial  bifurcation set on the $r_1-\tau$ plane and (b) the complete bifurcation set around HH are shown. }
     \label{fig:predatorr1tau}
    \end{figure}

  \subsection{Simulations}
\begin{figure}
               \centering
        \includegraphics[width=0.99\textwidth]{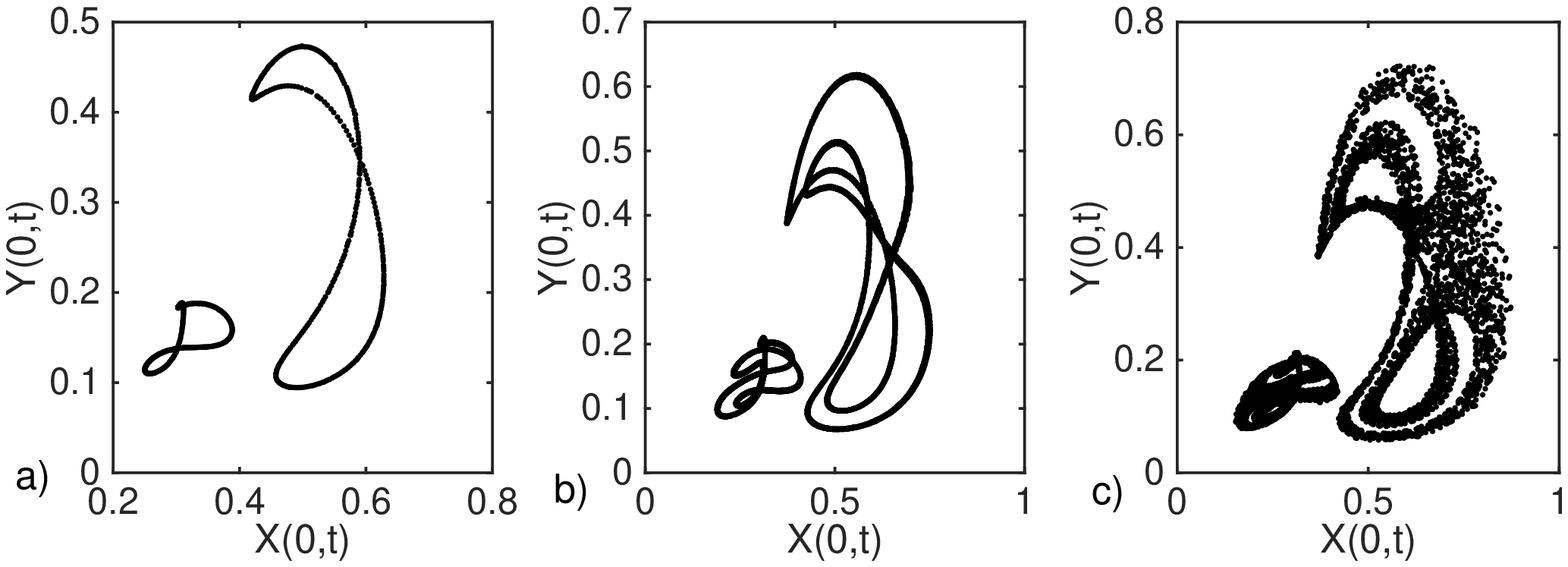}
    \caption  {  (a) When $\tau=10.8$, $r_1=0.69$, there exists a stable quasi-periodic solution on a 2-torus of system (\ref{predator1}); (b) when  $\tau=10.8$, $r_1=0.71$, there exists a quasi-periodic     solution on a 3-torus of system (\ref{predator1}); (c) $\tau=10.8$, $r_1=0.726$, there is a strange attractor of system (\ref{predator1}), and the system exhibits chaotic behavior.}
     \label{fig:predator3figure}
    \end{figure}

In this section, we fix  $r_2=1$, $a_{11}=1$, $a_{12}=1.2$, $a_{21}=2.8$, $a_{22}=1$, $d_1=0.1$, $d_2=0.2$, $l=3$, and take $r_1$ and $\tau$ as bifurcation parameters. We can find that when $0.316<r_1<1.47$, $(H_5)$ holds for all $n=0$, and when $0.316<r_1<1.19$, $(H_5)$ holds for all $n=1$. By (\ref{taunj}),  we can draw the curves of Hopf bifurcation values when $r_1$ varies, which is shown in Figure \ref{fig:predatorr1tau}(a). When $r_1=0.6739271475$,   $\tau_0^{0-}$ intersects $\tau_0^{1+}$, and we denote the double Hopf bifurcation point HH. For HH, we have $\omega_0^+=  0.77444$, $\omega_0^-=   0.362170$, and $\tau_0^{0-}=\tau_0^{1+}= 10.4238045$. We get that $B_{11} =  0.17069 + 0.12592i$,
 $B_{21 }=  1.97884 + 2.26811i$, $B_{13} =
 -0.11732 + 0.09868i$, $B_{23 }=
 -3.48553 + 1.52722i$,
$B_{2100} = -0.58923 - 0.57368i$, $B_{1011} =
 -7.57432 +11.46887i$, $B_{0021} =
  3.603569 - 7.55242i$, $B_{1110} =
  0.374678+ 2.89123i$, $\epsilon_1 =-1$, $\epsilon_2 =  1$, $b_0 =  2.10189$, $c_0 =
  -0.63587$, $d_0 =    -1$, and $d_0-b_0c_0 =    0.336547$.  It means that the unfolding system is of type VIa, and the bifurcation set is shown in Figure \ref{fig:predatorr1tau}(b) in which the phase portraits in D1-D8  have been shown in Figure \ref{fig:VIa}.  System  (\ref{predator1}) has a quasi-periodic solution on a 2-torus and a quasi-periodic solution on a 3-torus, which will be eliminated by a saddle connection orbit. According to  the ``Ruelle-Takens-Newhouse" scenario to chaos \cite{P. Battelino,J.P. Eckmann,D. Ruelle},  a vanishing 3-torus might accompany strange attractors, and lead a system into chaos. Thus, we know that near point HH, there may exist a strange attractor, which will be shown numerically in Figure  \ref{fig:predator3figure}.

In  Figure  \ref{fig:predator3figure},
we  show the existence of quasi-periodic solutions with multiple frequencies, and show the results through the Poincar\'e map on a  Poincar\'e section.  Since all periodic or quasi-periodic solutions are  spatially homogeneous, we choose the solution curve  of $(X(0,t),Y(0,t))$, i.e., $x=0$.  Moreover, because the periodic solutions oscillate in an infinite dimensional phase space \cite{JWu}, we give simulations near the double Hopf bifurcation point on the Poincar\'e section $X(0,t-\tau)=X^*$. We can see that the system exhibits rich dynamical behavior near the bifurcation point HH in Figure \ref{fig:predator3figure}.  System  (\ref{predator1}) has a quasi-periodic solution on a 2-torus which becomes a quasi-periodic solution on a 3-torus and then vanishes through a heteroclinic orbit, which is shown in Figure \ref{fig:predator3figure} (a) and (b) respectively.   After the vanishing of 3-torus, system  (\ref{predator1}) has a strange attractor, and exhibits chaotic behavior, which is shown in \ref{fig:predator3figure} (c).

     \section{ Concluding Remarks}
      \label{Conclusion}
   In this paper, we have extended the center manifold reduction and normal form method to analyze the dynamical behavior near the double Hopf bifurcation point in a delayed reaction-diffusion system. The method appears to be quite complicated as proceeded in this paper, but it is still an explicit algorithm and is not difficult to be implanted into a computer program. What we should care more about in this method is to calculate the double Hopf point and all the rest part can be exactly obtained by following the procedure we have given.

    After normal form derivation with respect to a delayed reaction-diffusion system, the dynamics near the bifurcation point, like the case in ordinary differential equations, is also governed by twelve distinct kinds of unfolding systems, and the bifurcation set about ordinary differential equations for each of the twelve types of unfoldings still applies to our model.  To show that our method is a powerful method for analyzing local behavior around a double Hopf bifurcation point, we give two examples: in the stage-structured epidemic model, we show that two spatially inhomogeneous periodic oscillations coexist near the singular point; in the predator-prey system, we show  the existence of  quasi-periodic solution on a 2-torus,  quasi-periodic solution on a 3-torus, and  strange attractor appears near the bifurcation point.

\begin{acknowledgements}
The authors are grateful to the handling editor and anonymous  referees for their careful reading of the manuscript and  valuable comments, which improve the exposition of the paper very much.  This research is supported by  National Natural Science Foundation of China (11701120, 11771109) and   Shaanxi Provincial Education Department grant (18JK0123).
\end{acknowledgements}


\begin{thebibliography}{}
%
%

 \bibitem{Q. An} An, Q.,  Jiang, W:  Spatiotemporal attractors generated by the Turing-Hopf bifurcation in a time-delayed reaction-diffusion system. Disctete Cont. Dyn-B  doi: 10.3934/dcdsb.2018183 (2018)

\bibitem{Andronov} Andronov, A.A.: Application of Poincar\'{e}  theorem on  bifurcation points and change in stability to simple auto-oscillatory systems.  C. R. Acad. Sci. Paris  189, 559-561 (1929)


\bibitem{AK Bajaj} Bajaj, A.K., Sethna, P.R.: Bifurcations in three-dimensional motions of articulated tubes. I - Linear systems and symmetry. II - Nonlinear analysis. J. Appl. Mech. 49, 606-618 (1982)


\bibitem{P. Battelino}Battelino, P.M., Grebogi, C., Ott, E., Yorke, J.A.: Chaotic attractors on a 3-torus, and torus break-up. Physica D 39, 299-314 (1989)


 \bibitem{M. Baurmann} Baurmann, M., Gross, T., Feudel, U.: Instabilities in spatially extended predator-prey systems: spatio-temporal patterns in the neighborhood of Turing-Hopf bifurcations. J. Theor. Bio. 245, 220-229 (2007)

\bibitem{J Belair}Belair, J., Campbell, S.A., Driessche, P.V.D.: Frustration, stability, and delay-induced oscillations in a neural network model. SIAM. J. Appl. Math. 56, 245-255 (1996)

\bibitem{Ping BiTumor}Bi, P., Ruan, S.: Bifurcations in delay differential equations and applications to tumor and
immune system interaction models. SIAM J. Appl. Dyn. Syst. 12, 1847-1888 (2013)

\bibitem{P.L. Buono}Buono, P.L., B\'{e}lair, J.:  Restrictions and unfolding of double Hopf bifurcation in functional differential equations. J. Differ. Equ. 189,
234-266 (2003)

\bibitem{S.A. Campell} Campell, S.A., B\'{e}lair, J.:  Analytical and symbolically-assisted investigation of Hopf bifurcations in delay-differential equations. Can. Appl. Math. Q. 3, 137-154 (1995)

\bibitem{S.A. Campellfeedback} Campell, S.A., B\'{e}lair, J., Ohira, T., Milton, J.: Limit cycles, tori, and complex dynamics in a second-order differential equations with delayed negative feedback. J. Dyn. Differ. Equ. 7, 213-236 (1995)



\bibitem{S. A. CampbellResonant} Campell, S.A.,  LeBlanc, V.G.: Resonant Hopf-Hopf interaction in delay differential equations.
J. Dyn. Differ. Equ. 10, 327-346 (1998)

 \bibitem{S. Chen} Chen, S., Shi, J., Wei, J.: Global stability and Hopf bifurcation in a delayed diffusive Leslie-Gower predator-prey system, Int. J. Bifurcat. Chaos  22, 331-517 (2012)

\bibitem{S. Chennonlocal} Chen, S., Yu, J.: Stability and bifurcations in a nonlocal delayed reaction-diffusion population model. J. Differ. Equations 260, 218-240 (2016)

 \bibitem{A. De Wit}  De Wit, A., Dewel, G., Borckmans, P.: Chaotic Turing-Hopf mixed mode. Phys. Rev. E 48, R4191-R4194 (1993)

\bibitem{duguo} Du, Y., Guo, Y., Xiao, P.: Freely-moving delay induces periodic oscillations in a structured SEIR model. Int. J. Bifurcat. Chaos  27, 1750122 (2017)




\bibitem{J.P. Eckmann}Eckmann, J.P.: Roads to turbulence in dissipative dynamical systems. Rev. Modern Phys. 53, 643-654 (1981)

\bibitem{C. Elphick} Elphick, C., Tiraopegui, E., Brachet, M.E., Coullet, P., Iooss, G.: A simple global characterization for normal forms of singular vector fields. Physica D  29, 95-127 (1987)

\bibitem{Faria}  Faria, T.:  Normal forms and Hopf bifurcation for partial differential equations with delays. Trans. Amer. Math. Soc. 352, 2217-2238 (2000)

\bibitem{T. Fariapredator-preydiffusion} Faria, T., Stability and bifurcation for a delayed predator-prey model and the effect of diffusion. J. Math. Anal. Appl. 254,  433-463 (2001)

 \bibitem{Faria W. Huang}Faria, T., Huang, W.: Stability of periodic solutions arising from Hopf bifurcation for a reaction-diffusion equation with time delay. Fields Inst. Comm. 31, 125-141  (2002)

\bibitem{FariaJDE} Faria, T., Magalh\~aes, L.T.: Normal forms for retarded functional differential equations with parameters and applications to Hopf bifurcation. J.  Differ. Equations  122, 181-200 (1995)


\bibitem{T. FariaBT}Faria, T., Magalh\~aes, L.T.: Normal form for retarded functional differential equations and applications to Bogdanov-Takens singularity. J. Differ. Equations  122, 201-224 (1995)





\bibitem{S. VAN GILS}Gils, S.A.V., Krupa, M., Langford, W.F.:  Hopf bifurcation with non-semisimple 1:1 resonance.
Nonlinearity 3, 825-850 (1990)

\bibitem{W. GOVAERTS} Govaerts, W.,  Guckenheimer, J., Khibnik, A.: Defining functions for multiple Hopf bifurcations. SIAM J. Numer. Anal. 34, 1269-1288 (1997)

\bibitem{Guckenheimer}
Guckenheimer, J., Holmes, P.: Nonlinear Oscillations, Dynamical Systems, and Bifurcations of Vector Fields.
Springer, New York (1983)

\bibitem{S. Guononlocal} Guo, S.: Stability and bifurcation in a reaction-diffusion model with nonlocal delay effect. J.  Differ. Equations   259, 1409-1448 (2015)

\bibitem{S. GuoDirichlet} Guo, S., Ma, L.: Stability and bifurcation in a delayed reaction-diffusion equation with Dirichlet boundary condition. J. Nonl. Sci. 26, 545-580 (2016)

\bibitem{Hale} Hale, J.K., Kocak,  H.: Dynamics and Bifurcations. Springer, New York  (1991)

\bibitem{J. Halefde} Hale, J.K., Lunel, S.M.V.:  Introduction to Functional Differential Equations. Springer, New York (1993)
  \bibitem{B. Hassard}Hassard, B.D., Kazarinoff, N.D., Wan, Y.H.: Theory and Applications of Hopf Bifurcation. Cambridge Univ. Press, New York  (1981)

\bibitem{HW Hethcote}Hethcote, H.W., Lewis, M.A.,  Driessche, P.V.D.: An epidemiological model with a delay and a nonlinear incidence rate. J. Math. Biol. 27, 49-64 (1989)

\bibitem{Hopf}Hopf, E.: Abzweigung einer periodischen l\"{o}sung eines Differential Systems. Berichen Math.
Phys. Kl. S\"{a}ch. Akad. Wiss. Leipzig  94, 1-22 (1942)

 \bibitem{S. B. HSUpredator}Hsu, S.B., Huang, T.W.:  Global stability for a class of predator-prey systems. SIAM J. Appl. Math. 55, 763-783 (1995)

\bibitem{Ji} Ji, J., Li,  X., Luo, Z.: Two-to-one resonant Hopf bifurcations in a quadratically nonlinear oscillator involving time delay. Int. J. Bifurcat. Chaos  22, 1250060 (2012)

\bibitem{H. Kielhfer}Kielh\"{o}fer, H.: Bifurcation Theory: An Introduction with Applications to Partial Differential Equations. Springer, New York (2011)

\bibitem{Kuznetsov} Kuznetsov, Y.A.:  Elements of Applied Bifurcation Theory. Springer,  New York (2011)

       \bibitem{GREGORY}Lewis, G.M., Nagata, W.: Double Hopf bifurcations in the differentially heated rotating annulus. SIAM J. Appl. Math. 63, 1029-1055 (2003)

   \bibitem{Lin}Lin, X., So, J.W.H., Wu, J.: Centre manifolds for partial differential equations with delays. P. Roy. Soc. Edinb. A 122, 237-254 (1992)

\bibitem{Luongo}Luongo, A., Paolone, A.: Perturbation methods for bifurcation analysis from multiple nonresonant complex eigenvalues. Nonlinear Dynam.  14, 193-210 (1997)
\bibitem{Suqi Ma}Ma, S., Lu, Q., Feng, Z.: Double Hopf bifurcation for van der Pol-Duffing oscillator with parametric delay feedback control. J. Math. Anal. Appl. 338, 993-1007 (2008)
 \bibitem{M. Meixner} Meixner, M., De Wit, A., Bose, S., Sch\"{o}ll, E.: Generic spatiotemporal dynamics near codimension-two Turing-Hopf bifurcations. Phys. Rev. E 55, 6690-6697 (1997)

\bibitem{Poincare} Poincar\'{e}, H.:  Les M\'{e}thodes Nouvelles de la M\'{e}canique C\'{e}leste. Cauthier-Villars, Paris  (1892)

\bibitem{DVR Reddy} Reddy, D.V.R., Sen, A., Johnston, G.L.: Time delay effects on coupled limit cycle oscillators at Hopf bifurcation. Physica D  129, 15-34  (1999)

\bibitem{Revel}Revel, G., Alonso, D.M., Moiola, J.L.: Interactions between oscillatory modes near a 2:3 resonant Hopf-Hopf bifurcation. Chaos 20, 113-129  (2010)

\bibitem{Revel1:2}Revel, G., Alonso, D.M., Moiola, J.L.: Numerical semi-global analysis of a 1:2 resonant Hopf-Hopf bifurcation. Physica D 247, 40-53  (2013)

 \bibitem{S. RUANpredator} Ruan, S., Xiao, D.: Global analysis in a predator-prey system with nonmonotonic functional response. SIAM J.  Appl. Math. 61, 1445-1472 (2000)

\bibitem{D. Ruelle}Ruelle, D., Takens, F.: On the nature of turbulence. Comm. Math. Phys. 20, 167-192  (1971)

 \bibitem{Y SongHopf}Song, Y., Wei, J.: Local Hopf bifurcation and global periodic solutions in a delayed predator-prey system. J. Math. Anal. Appl. 301,  1-21 (2005)


\bibitem{YongliSongTuringHopf}Song, Y., Zhang, T., Peng, Y.: Turing-Hopf bifurcation in the reaction-diffusion equations and its applications. Commun. Nonlinear Sci. Numer. Simulat. 33, 229-258 (2016)

\bibitem{PH Steen}Steen, P.H., Davis, S.H.:
Quasiperiodic bifurcation in nonlinearly-coupled oscillators near a point of strong resonance. SIAM J. Appl. Math. 42, 1345-1368 (1982)


\bibitem{Y. SuDirichlet} Su, Y., Wei, J., Shi, J.: Hopf bifurcations in a reaction-diffusion population model with delay effect. J. Differ. Equations  247, 1156-1184 (2009)

\bibitem{Wiggins}Wiggins, S.: Introduction to Applied Nonlinear Dynamical Systems and Chaos. Springer, New York  (2003).

\bibitem{JWu} Wu, J.: Theory and Applications of Partial Functional-Differential Equations. Springer, New York (1996)

\bibitem{D. XIAOpredator}Xiao, D.:  Bifurcations of a ratio-dependent predator-prey system with constant rate harvesting. SIAM  J. Appl. Math. 65, 737-753 (2005)

 \bibitem{YXiao} Xiao, Y., Chen, L.:  An SIS epidemic model with stage structure and a delay. Acta Math. Appl. Sin. E. 18, 607-618  (2002)

 \bibitem{X. Xu}Xu, X., Wei, J.: Turing-Hopf bifurcation of a class of modified Leslie-Gower model with diffusion. Disc. Continu. Dyn. Sys. B 23, 765-783 (2018)

\bibitem{X.-P. YanDirichlet} Yan, X., Li, W.: Stability of bifurcating periodic solutions in a delayed reaction-diffusion population model. Nonlinearity 23, 1413-1431 (2010)
\bibitem{F. Yipredator-prey}Yi, F., Wei, J., Shi, J.: Bifurcation and spatiotemporal patterns in a homogenous diffusive predator-prey system. J. Differ. Equations  246, 1944-1977 (2009)

\bibitem{PEI YU}Yu, P.: Analysis on double Hopf bifurcation using computer algebra
with the aid of multiple scales. Nonlinear Dynamics 27, 19-53 (2002)

\bibitem{Yu}Yu, P., Bi, Q.: Analysis of non-linear dynamics and bifurcations of a double pendulum. J. Sound Vib. 217, 691-736  (1998)


\bibitem{P. Yu} Yu, P., Yuan, Y., Xu, J.: Study of double Hopf bifurcation and chaos for oscillator with time delay feedback. Commun. Nonlinear Sci. Numer.
Simul. 7, 69-91 (2002)



\bibitem{Y.Y. Zhang} Zhang, Y., Xu, J.: Classification and computation of non-resonant double Hopf bifurcations and solutions in delayed van der Pol-Duffing system.  Int. J. Nonlinear Sci. Numer. Simul. 6, 67-74 (2005)




\end{thebibliography}


\end{document}